\documentclass{amsart}
\usepackage{amsmath, amssymb, latexsym}
\title{Flag f-vectors of three-colored complexes}
\author{Andrew Frohmader}
\address{Department of Mathematics, 581 Malott Hall, Cornell University, Ithaca, NY 14853-4201}
\email{froh@math.cornell.edu}

\newtheorem{theorem}{Theorem}[section]

\newtheorem{lemma}[theorem]{Lemma}
\newtheorem{definition}[theorem]{Definition}

\newtheorem{example}[theorem]{Example}

\def\proof{\smallskip\noindent {\it Proof: \ }}
\def\endproof{\hfill\ensuremath{\square}\medskip}

\begin{document}

\maketitle

\begin{abstract}
The flag f-vectors of three-colored complexes are characterized.  This also characterizes the flag h-vectors of balanced Cohen-Macaulay complexes of dimension two, as well as the flag h-vectors of balanced shellable complexes of dimension two.
\end{abstract}

\section{Introduction}

In the late 1970s, Stanley \cite{stanley} showed that two particular classes of simplicial complexes have equivalent characterizations of their flag f-vectors.  Several years later, Bj\"{o}rner, Frankl, and Stanley \cite{bfs} showed that two additional classes of simplicial complexes shared this same characterization.  Unfortunately, no one has a characterization for any of these classes of simplicial complexes, but we only know that characterizing one would suffice for all four.  There are some already known cases that are trivial.  In this paper, we solve one of the simplest non-trivial cases by characterizing the flag f-vectors of 3-colored simplicial complexes.

Recall that a \textit{simplicial complex} $\Delta$ on a vertex set $W$ is a collection of subsets of $W$ such that (i) for every $v \in W$, $\{v\} \in \Delta$ and (ii) for every $B \in \Delta$, if $A \subset B$, then $A \in \Delta$.  The elements of $\Delta$ are called \textit{faces}.  A face on $i$ vertices is said to have \textit{dimension} $i-1$, while the dimension of a complex is the maximum dimension of a face of the complex.

The \textit{$i$-th f-number} of a simplicial complex $\Delta$, $f^{i-1}(\Delta)$ is the number of faces of $\Delta$ on $i$ vertices.  The \textit{f-vector} of $\Delta$ lists the f-numbers of $\Delta$.  One interesting question to ask is which integer vectors can arise as f-vectors of simplicial complexes.  Much work has been done toward answering this for various classes of simplicial complexes.  For example, the Kruskal-Katona theorem \cite{kruskal, katona} characterizes the f-vectors of all simplicial complexes.

In this paper, we wish to deal with colored complexes, where the coloring provides additional data.  A \textit{coloring} of a simplicial complex is a labeling of the vertices of the complex with colors such that no two vertices in the same face are the same color.  Because any two vertices in a face are connected by an edge, this is equivalent to merely requiring that any two adjacent vertices be assigned different colors.  If the set of colors has $n$ colors, we refer to the colors as $1, 2, \dots , n$.  The set of colors is denoted by $[n] = \{1, 2, \dots , n \}$.  The \textit{color set} of a face is the subset of $[n]$ consisting of the colors of the vertices of the face.  The Frankl-F\"{u}redi-Kalai \cite{ffk} theorem characterizes the f-vectors of all simplicial complexes that can be colored with $n$ colors.

We wish to use a refinement of the usual notion of f-vectors. The \textit{flag f-numbers} of a colored simplicial complex $\Delta$ on a color set $[n]$ are defined by, for any subset $S \subset [n]$, $f_S(\Delta)$ is the number of faces of $\Delta$ whose color set is $S$.  The \textit{flag f-vector} of $\Delta$, $f(\Delta)$ is the collection of the flag f-numbers of $\Delta$ for all subsets $S \subset [n]$.

For simplicity, we sometimes refer to colors by their numbers and drop the brackets when we do so.  For example, $f_3(\Delta)$ is the number of vertices of $\Delta$ of color 3. Similarly, $f_{12}(\Delta)$ is the number of edges of $\Delta$ with one vertex of color 1 and one vertex of color 2.

This is a refinement of the usual notion of f-numbers and the f-vector of a complex.  The relation between the f-numbers and the flag f-numbers is that the former ignores the colors.  It can be computed from the latter as
$$f^{i-1}(\Delta) = \sum_{|S| = i} f_S(\Delta).$$

The f-numbers of a complex are usually written with the number as a subscript, not a superscript.  We do not do this because this paper is mainly interested in flag f-numbers, and we wish to be able to drop the brackets and write, for example, $f_{12}(\Delta)$ rather than $f_{\{1, 2\}}(\Delta)$ without it being mistaken for the number of twelve dimensional faces of $\Delta$.

One can ask which nonnegative integer vectors can arise as the flag f-vectors of colored simplicial complexes.  It can help to define the flag h-vector of a complex by
$$h_S(\Delta) = \sum_{T \subset S} f_T(\Delta) (-1)^{|S| - |T|}.$$
The flag h-vector of a complex contains the same information as the flag f-vector, and is easier to work with in some cases.  If given the flag h-vector, we can recover the flag f-vector by
$$f_S(\Delta) = \sum_{T \subset S} h_T(\Delta).$$

Stanley \cite{stanley} showed that the flag h-vector of a balanced Cohen-Macaulay complex is the flag f-vector of a simplicial complex and vice versa.  That is, if a Cohen-Macaulay complex has dimension $n-1$ and can be colored with $n$ colors, then there is a simplicial complex that can be colored with $n$ colors whose flag f-vector is the flag h-vector of the Cohen-Macaulay complex.  Bj\"{o}rner, Frankl, and Stanley \cite{bfs} showed that both of these are also equivalent to a vector being the flag h-vector of a balanced shellable complex, and also to being the flag f-vector of a color-shifted complex.  The problem here is that while four different classes of complexes have equivalent characterizations, none of them have a known characterization.

One might hope that stronger local restrictions than what Bj\"{o}rner, Frankl, and Stanley found could be placed upon the complexes without changing the characterization of the flag f-vectors, and work toward a solution that way.  For example, the Kruskal-Katona theorem says that to characterize the f-vectors of simplicial complexes, we can restrict to the ``rev-lex" complexes.  As there is only one possible rev-lex complex for a given f-vector, this effectively solves the problem.  Frankl, F\"{u}redi, and Kalai \cite{ffk} did something similar to characterize the f-vectors of colored complexes.

Indeed, the paper of Bj\"{o}rner, Frankl, and Stanley already did impose stronger restrictions to some extent.  Every color-shifted colored complex is, in particular, a colored complex, so they showed that in order to characterize the flag f-vectors of colored complexes, it sufficed to consider only color-shifted complexes.  Likewise, every shellable complex is Cohen-Macaulay, so they showed that to characterize the flag h-vectors of balanced Cohen-Macaulay complexes, it suffices to consider only the balanced shellable complexes.  However, another paper of the author \cite{nosharper} showed that one cannot impose stronger local conditions than color-shifting in a certain sense.

Another approach is to try to impose some bounds.  Walker \cite{walker} showed that the only linear inequalities on the flag f-numbers of simplicial complexes are the trivial ones, namely, that all flag f-numbers are non-negative.  In the same paper, he computed all linear inequalities on the logarithms of the flag f-numbers of a simplicial complex.  These give inequalities on the products of flag f-numbers.  For example, the most trivial case is that $f_1(\Delta) f_2(\Delta) \geq f_{12}(\Delta)$, as any edge whose color set is $\{1, 2\}$ must use a vertex of color 1 and a vertex of color 2, and the ways to pick these vertices are $f_1(\Delta)$ and $f_2(\Delta)$, respectively.  While an interesting result, Walker's result is shy of a full characterization of the flag f-vectors of simplicial complexes in multiple ways.

First, it avoids dealing with discreteness issues.  For example, if $f_{12}(\Delta) = f_{13}(\Delta) = f_{23}(\Delta) = 5$, then Walker's result shows that $f_{123}(\Delta) \leq 5\sqrt{5} \approx 11.18$.  Since flag f-numbers must be integers, this immediately gives that $f_{123}(\Delta) \leq 11$.  By checking the possible cases, it is possible to obtain an upper bound of $f_{123}(\Delta) \leq 9$.

What happens here is that discreteness gets in the way to make Walker's bounds not sharp.  If vertices were not discrete, countable things, then we could take a complete tripartite graph on $\sqrt{5}$ vertices of each color and hit his bound exactly.

Another issue is that, even if given sharp bounds on flag f-numbers of certain color sets in terms of their flag f-numbers on proper subsets, these can sometimes conflict with each other.  For example, suppose that $f_{12}(\Delta) = 1$, $f_{13}(\Delta) = 2$, $f_{23}(\Delta) = 2$, $f_{24}(\Delta) = 2$, and $f_{34}(\Delta) = 1$.  Walker's bounds give that $f_{123}(\Delta) \leq 2$ and $f_{234}(\Delta) \leq 2$.  It is possible to obtain either one of these.  The former is a complete tripartite graph on two vertices of color 3 and one vertex each of colors 1 and 2.  The latter is a complete tripartite graph on two vertices of color 2 and one vertex each of colors 3 and 4.

It is not possible to make $f_{123}(\Delta) = 2$ and $f_{234}(\Delta) = 2$ simultaneously, however.  The former requires that the edges with color set $\{2, 3\}$ have a common vertex of color 2, while the latter requires that the edges have a common vertex of color 3.

Walker's bounds are enough to settle the case of two colors.  A proposed nonnegative integer flag f-vector corresponds to a non-empty two-colored simplical complex if and only if $f_{\emptyset}(\Delta) = 1$ and $f_1(\Delta)f_2(\Delta) \geq f_{12}(\Delta)$.  The problem remains open for more colors, however.

In this paper, we give a characterization of the flag f-vectors of 3-colored complexes in Theorem~\ref{maintheorem}.  Discreteness issues are accounted for, so we do produce sharp bounds.  Given any prospective flag f-vector for a 3-colored complex, we can either give a complex that has the desired flag f-vector or show that no such complex exists.  The problem of different higher dimensional flag f-numbers forcing different configurations on lower dimensional faces only appears when there are at least four colors.

Our solution consists of greatly restricting the class of simplicial complexes and then checking those that remain to see whether any give enough facets.  The number of complexes we must check depends on the number of vertices and edges of each color set allotted.  If the number of edges of each color set is chosen independently and uniformly at randomly from $[n]$, then the expected number of complexes we must check to find the one that maximizes the number of facets is less than 15, independent of $n$ and regardless of how many vertices of each color are allowed.  In the worst possible case, we check on the order of $n^{{1 \over 4}}$ complexes.

In Section~2, we give our characterization of the flag f-vectors of 3-colored complexes.  In Section~3, we give some examples of computations to characterize the flag f-vectors of 3-colored complexes.  These computations demonstrate why there shouldn't be a much nicer characterization.  In Section~4, we discuss the analogous problem for more than three colors.

\section{Three colors}

As we have seen, getting a complete characterization of the flag f-vectors of colored complexes would be quite difficult.  The case of only two colors is trivial, however.  In this section, we give a solution to the case of three colors.  The case of three colors is enough that discreteness issues matter, unlike how they can be ignored with only two colors.  Three colors are still few enough that the issue of different sets of flag f-numbers trying to force different and contradictory sets of faces on the same color set cannot appear.

We begin by dispensing with a couple of trivialities.  First, if $f_S(\Delta) = 0$ for some $S \subset [3]$, then $\Delta$ has no two-dimensional faces, leaving only trivial obstructions to a prospective flag f-vector corresponding to an actual complex.  As such, for the rest of this section and the next one, we assume that $f_S(\Delta) > 0$ for all $S \subset [3]$.  This also allows us to refer to the two-dimensional faces of a three-colored complex as \textit{facets}.

Next, there is one empty set of vertices.  The empty set is a face of every non-empty colored complex because it is a subset of the vertex set of some other face.  Therefore, for every non-empty colored complex $\Delta$, we get $f_{\emptyset}(\Delta) = 1$.  We thus assume that any prospective flag f-vector has $f_{\emptyset}(\Delta) = 1$ and do not further bother with empty faces in the rest of this paper.

Our characterization is to take given $f_1(\Delta), f_2(\Delta), f_3(\Delta), f_{12}(\Delta), f_{13}(\Delta)$, and $f_{23}(\Delta)$ satisfying $f_{12}(\Delta) \leq f_1(\Delta)f_2(\Delta)$, $f_{13}(\Delta) \leq f_1(\Delta)f_3(\Delta)$, and $f_{23}(\Delta) \leq f_2(\Delta)f_3(\Delta)$, and compute the largest possible value for $f_{123}(\Delta)$.  We can restrict to the situation where these inequalities hold, as otherwise it is trivial that no $\Delta$ with the desired flag f-vector can exist.  Apart from this, if we have a proposed flag f-vector and wish to know whether it corresponds to any 3-colored complex, we can answer the question by comparing the proposed value of $f_{123}(\Delta)$ to the computed maximum value.

The basic idea of the characterization is to put increasingly stronger restrictions on the class of complexes to be considered.  Eventually the class of complexes is small enough that we can check it by brute force in a reasonable amount of time, as it typically ends up only being around six or ten complexes to check.

A first step in this direction is color-shifting.  We can place an arbitrary order on the vertices of each color.  We label the $j$-th vertex of color $i$ as $v_j^i$, so that the vertices of color $i$ are $v_1^i, v_2^i, \dots, v_{f_i(\Delta)}^i$.

\begin{definition}
\textup{Let $\Delta$ be an $n$-colored simplicial complex.  We say that $\Delta$ is \textit{color-shifted} if, for all $b_1 \leq a_1, b_2 \leq a_2, \dots , b_j \leq a_j$, $\{v_{a_1}^{i_1}, v_{a_2}^{i_2}, \dots v_{a_j}^{i_j}\} \in \Delta$ implies $\{v_{b_1}^{i_1}, v_{b_2}^{i_2}, \dots v_{b_j}^{i_j}\} \in \Delta$.}
\end{definition}

\begin{theorem} \label{colorshift}
Let $\Delta$ be an $n$-colored simplicial complex.  Then there is an $n$-colored, color-shifted simplicial complex $\Gamma$ such that $f_S(\Delta) = f_S(\Gamma)$ for all $S \subset [n]$.
\end{theorem}

This was proven by Bj\"{o}rner, Frankl, and Stanley as part of \cite[Theorem 1]{bfs}, in which they proved the equivalence of characterizing the flag f-vectors of four different classes of complexes.  They called the concept ``compressed" rather than color-shifted.  Furthermore, their proof allowed for a more general notion of coloring where, for example, one could have three colors, but allow a face to have up to 3 vertices of color 1, up to 5 vertices of color 2, and up to 2 vertices of color 3.  In this paper, we focus on the case where only one vertex of each color is allowed in a face.

In proofs throughout this section, we define complexes by specifying the edges.  The facets are assumed to be all possible facets for which all three edges are present in the complex.  We can do this because we are trying to find the complex that maximizes the number of facets, and discarding some facets from a complex obviously cannot increase the number of facets.

The intuition behind our characterization is fairly simple.  Start by considering what would happen if vertices were not discrete things, that is, if we could have non-integer numbers of vertices.  For example, in the f-vector case, if we wanted to have as many edges as possible on a bipartite graph with exactly 5 vertices, a standard exercise in elementary calculus would show that the optimal solution is to have two and a half vertices in each part, yielding six and a quarter edges total.

We then step back and say, we can't have half of a vertex or a quarter of an edge, but the real solution is probably close to this.  Indeed, in this simple example, the exact solution is very close:  a bipartite graph with three vertices in one part and two in the other has six edges.  We check a handful of 3-colored complexes that are close to the continuous construction and show that at least one of them must maximize $f_{123}(\Delta)$.

The most complicated case is what happens if we have plenty of vertices, so that the important restrictions are the numbers of edges.  If the numbers of vertices are significant restrictions, then often they make the problem trivial.  At the very least, they restrict what it is necessary to check.

If the main restrictions come from the numbers of edges, one might guess that the optimal complex is the complete tripartite complex on $c_1$ vertices of color 1, $c_2$ vertices of color 2, and $c_3$ vertices of color 3, for suitable constants $c_1, c_2$, and $c_3$.  The relevant equations are $f_{12}(\Delta) = c_1c_2$, $f_{13}(\Delta) = c_1c_3$, and $f_{23}(\Delta) = c_2c_3$.  We can solve for the constants to get $c_1 = \sqrt{{f_{12}(\Delta)f_{13}(\Delta) \over f_{23}(\Delta)}}$, $c_2 = \sqrt{{f_{12}(\Delta)f_{23}(\Delta) \over f_{13}(\Delta)}}$, and $c_3 = \sqrt{{f_{13}(\Delta)f_{23}(\Delta) \over f_{12}(\Delta)}}$.  With this in mind, we make the following definitions.

\begin{definition}
\textup{Let $\Delta$ be a 3-colored simplicial complex with flag f-vector $f(\Delta)$.  Define}
\begin{eqnarray*}
b_1(\Delta) & = & \Bigg\lfloor \sqrt{{f_{12}(\Delta)f_{13}(\Delta) \over f_{23}(\Delta)}} \Bigg\rfloor, \\ b_2(\Delta) & = & \Bigg\lfloor \sqrt{{f_{12}(\Delta)f_{23}(\Delta) \over f_{13}(\Delta)}} \Bigg\rfloor, \qquad \textup{and} \\ b_3(\Delta) & = & \Bigg\lfloor \sqrt{{f_{13}(\Delta)f_{23}(\Delta) \over f_{12}(\Delta)}} \Bigg\rfloor.
\end{eqnarray*}
\end{definition}

For simplicity, we refer to an edge of color set $\{1, 2\}$ as being an edge of color 12, and similarly for edges of color 13 or 23.  With only three colors, this cannot lead to ambiguity about the color set intended.

The next several lemmas allow us to put much stronger conditions on a 3-colored complex than merely requiring it to be color-shifted.

\begin{lemma} \label{onevert}
Let $\Delta$ be a 3-colored simplicial complex.  Then there is a simplicial complex $\Gamma$ and positive integers $c_1, c_2, c_3, c_4, c_5$, and $c_6$ such that
\begin{enumerate}
\item $\Gamma$ is 3-colored,
\item $\Gamma$ is color-shifted,
\item $f_S(\Gamma) = f_S(\Delta)$ for all $S \subset [3]$ with $S \not = [3]$,
\item $f_{123}(\Gamma) \geq f_{123}(\Delta)$,
\item $\{v_{a_1}^1, v_{a_2}^2 \} \in \Gamma$ for all $a_1 \leq c_1$ and $a_2 \leq c_2$,
\item $\Gamma$ has at most one vertex other than $\{v_1^1, \dots, v_{c_1}^1, v_1^2, \dots, v_{c_2}^2 \}$ contained in an edge of $\Gamma$ of color 12.
\item $\{v_{a_3}^1, v_{a_4}^3 \} \in \Gamma$ for all $a_3 \leq c_3$ and $a_4 \leq c_4$,
\item $\Gamma$ has at most one vertex other than $\{v_1^1, \dots, v_{c_3}^1, v_1^3, \dots, v_{c_4}^3 \}$ contained in an edge of $\Gamma$ of color 13.
\item $\{v_{a_5}^2, v_{a_6}^3 \} \in \Gamma$ for all $a_5 \leq c_5$ and $a_6 \leq c_6$, and
\item $\Gamma$ has at most one vertex other than $\{v_1^2, \dots, v_{c_5}^2, v_1^3, \dots, v_{c_6}^3 \}$ contained in an edge of $\Gamma$ of color 23.
\end{enumerate}
\end{lemma}

\proof  Theorem~\ref{colorshift} ensures that there is a complex $\Delta_0$ that satisfies properties 1 through 4, though it could fail the rest.  We start with this complex and rearrange edges of one color at a time such that after each rearrangement, the complex satisfies the two properties relevant to that color, while retaining any numbered properties that it held before the rearrangement.  After rearranging edges of all three colors, we have the complex $\Gamma$.

Suppose that it is not possible to choose $c_1$ and $c_2$ and create a complex $\Delta_1$ satisfying properties 1 through 6 with the edges of $\Delta_1$ exactly the same as the edges of $\Delta_0$ except for those of color 12.  Taking $\Delta_1 = \Delta_0$ and $c_1 = c_2 = 1$ would leave the necessary faces unchanged and satisfy properties 1 through 5, so the only obstruction here is property 6.

Let $\Sigma$ be a complex that minimizes the number of vertices contained in an edge of color 12 among all rearrangements of the edges of $\Delta_0$ of color 12 that satisfy conditions 1 through 5.  Let the vertices of $\Sigma$ contained in an edge of color 12 be $\{v_1^1, \dots, v_{d_1}^1, v_1^2, \dots, v_{d_2}^2\}$.  Let the number of vertices of color 1 adjacent to the vertex $v_i^3$ be $p_i$ and the number of vertices of color 2 adjacent to $v_i^3$ be $q_i$.  Since $\Sigma$ is color-shifted, we must have $p_1 \geq p_2 \geq \dots \geq p_{f_3(\Delta)}$ and $q_1 \geq q_2 \geq \dots \geq q_{f_3(\Delta)}$.

Suppose first that $p_{f_3(\Delta)} \geq d_1$ and $q_{f_3(\Delta)} \geq d_2$.  In this case, every edge of color 12 together with a vertex of color 3 forms a facet of $\Sigma$.  Furthermore, rearranging edges of $\Sigma$ of color 12 does not change the number of facets provided that all edges only use the vertices $\{v_1^1, \dots, v_{d_1}^1, v_1^2, \dots, v_{d_2}^2\}$.  Let $w = \frac{f_{12}(\Delta)}{d_2}$.  Rearrange the edges of $\Sigma$ of color 12 such that the vertices $\{v_1^1, \dots, v_{\lfloor w \rfloor}^1\}$ are each adjacent to all of $\{v_1^2, \dots, v_{d_2}^2\}$ and the vertex $v_{\lfloor w \rfloor +1}^1$ is adjacent to $\{v_1^2, \dots, v_{d_2(w - \lfloor w \rfloor)}^2\}$.  Take $c_1 = \lfloor w \rfloor$ and $c_2 = d_2$.  It is easy to check that this new $\Sigma$ satisfies all of the conditions necessary for $\Delta_1$.

Otherwise, let $m = \min\{i\ |\ \textup{either } p_i < d_1 \textup{ or } q_i < d_2\}$.  Suppose without loss of generality that $p_m < d_1$; exactly the same argument applies if $q_m < d_2$.  If $f_{12}(\Delta) \leq (d_1 - 1)d_2$, then we can move the edges of $\Sigma$ of color 12 that contain $v_{d_1}^1$ to instead use one vertex in $\{v_1^1, \dots, v_{d_1-1}^1\}$ (without $v_{d_1}^1$) and one vertex in $\{v_1^2, \dots, v_{d_2}^2\}$, with the new edges chosen so as to keep the new complex color-shifted.  For $i < m$, this does not change the number of facets of $\Sigma$ containing $v_i^3$ because every edge of $\Sigma$ of color 12 together with $v_i^3$ forms a facet of $\Sigma$ both before and after the rearrangement.  It also cannot decrease the number of facets of $\Sigma$ containing $v_i^3$ for $i \geq m$, as $v_i^3$ together with the vertices of a removed edge did not form a facet as $\{v_i^3, v_{d_1}^1\} \not \in \Sigma$.  As such, the new $\Sigma$ satisfies conditions 1 through 5 while having one vertex fewer contained in an edge of color 12.  This contradicts the choice of $\Sigma$.

The other possibility is that $f_{12}(\Delta) > (d_1 - 1)d_2$. In this case, rearrange the edges of $\Sigma$ of color 12 such that the edges are every possible combination of one vertex in $\{v_1^1, \dots, v_{d_1-1}^1\}$ (without $v_{d_1}^1$) and one vertex in $\{v_1^2, \dots, v_{d_2}^2\}$, as well as that $v_{d_1}^1$ forms an edge with each vertex of $\{v_1^2, \dots, v_{f_{12}(\Delta) - (d_1 - 1)d_2}^2\}$.  As in the previous paragraph, this cannot decrease the number of facets of $\Sigma$ containing vertex $v_i^3$ with $i < m$ because $v_i^3$ would form a facet with the edge both before and after it is moved.  It cannot decrease the number of facets of $\Sigma$ containing $v_i^3$ for $i \geq m$, as $v_i^3$ together with the vertices of a removed edge did not form a facet as $\{v_i^3, v_{d_1}^1\} \not \in \Sigma$.  We can take $c_1 = d_1 - 1$, $c_2 = d_2$, and $\Delta^1 = \Sigma$ and satisfy conditions 1 through 6, as $v_{d_1}^1$ is the only extra vertex.

Now we repeat the process by rearranging the edges of other colors.  We can create $\Delta_2$ from $\Delta_1$ by rearranging the edges with color 13 in the same manner as how $\Delta_1$ was created.  This retains properties 1 through 4 for the same reasons that $\Delta_1$ did and makes properties 7 and 8 hold if $\Gamma = \Delta_2$ for the same reasons that properties 5 and 6 hold if $\Gamma = \Delta_1$.  Since the edges of color 12 are unchanged, properties 5 and 6 hold for $\Delta_2$ because they hold for $\Delta_1$.

Finally, we create $\Delta_3$ from $\Delta_2$ by rearranging edges of color 23.  This inherits properties 5 thought 8 from $\Delta_2$.  It retains properties 1 through 4 and adds properties 9 and 10 for the same reasons as happened analogously with $\Delta_1$ and $\Delta_2$.  Taking $\Gamma = \Delta_3$ completes the proof.  \endproof

\begin{lemma} \label{twovert}
Let $\Delta$ be a 3-colored simplicial complex.  It is possible to choose positive integers $g_1$, $g_2$, and $g_3$ and construct a 3-colored simplicial complex $\Gamma$ such that
\begin{enumerate}
\item $f_S(\Gamma) \leq f_S(\Delta)$ for all $S \subset [3]$ with $S \not = [3]$,
\item $f_{123}(\Gamma) \geq f_{123}(\Delta)$,
\item $\{v_{a_1}^1, v_{a_2}^2, v_{a_3}^3 \} \in \Gamma$ for all $a_1 \leq g_1$, $a_2 \leq g_2$, and $a_3 \leq g_3$, and
\item $\Gamma$ has at most two vertices other than $\{v_1^1, \dots, v_{g_1}^1, v_1^2, \dots, v_{g_2}^2, v_1^3, \dots, v_{g_3}^3 \}$ contained in a facet of $\Gamma$.
\end{enumerate}
Furthermore, it is possible to construct $\Gamma$ in the following manner for suitable choices of $p, q \in [3]$ with $p \not = q$.  Start with $g_1$ vertices of color 1, $g_2$ vertices of color 2, and $g_3$ vertices of color 3.  Make any two of these vertices not of the same color adjacent to each other.  If $f_p(\Delta) > g_p$, then add another vertex of color $p$ and make it adjacent to as many vertices of other colors as possible without violating condition 2, with the vertices of other colors chosen the earliest of their colors.  Similarly, if $f_q(\Delta) > g_q$, then add another vertex of color $q$ and make it adjacent to as many vertices of other colors as possible without violating condition 2.
\end{lemma}

\proof  By Lemma~\ref{onevert}, there is a simplicial complex $\Delta_1$ that satisfies all of the conditions of that lemma. Suppose first that some two of the extra vertices of conditions 6, 8, and 10 from Lemma~\ref{onevert} are the same color.  Since those conditions say at most one extra vertex, if there isn't an extra vertex, we can choose its color arbitrarily, and hence choose to make it match the color of one of the other vertices and be in this case.

We construct $\Gamma$ from $\Delta_1$ as follows.  Let $g_1 = \min\{c_1, c_3\}$, $g_2 = \min\{c_2, c_5\}$, and $g_3 = \min\{c_4, c_6\}$.  Let $p$ be the color of the two extra vertices and let $q$ be the color of the other extra vertex.  With these constants, construct $\Gamma$ as described in the statement of the lemma.  That conditions 1, 3 and 4 hold is immediate from the construction.  Condition 2 also holds because all of the facets of $\Delta_1$ are also in $\Gamma$, as the only edges of $\Delta_1$ not in $\Gamma$ are those for which at least one vertex is not adjacent to any vertex of a particular other color, so the edge is not contained in any facet.

Otherwise, the extra vertices are forced upon us by $\Delta_1$ and one extra vertex is of each color.  One possibility is that the extra vertex for an edge of color set $\{i,i+1\}$ (modulo 3) is of color $i$ for all $i \in [3]$; the other possibility is that the extra vertex is always of color $i+1$.  By symmetry, it suffices to consider only the former case.

Suppose that $c_1 \geq c_3$.  The extra vertex of color 1 is then not adjacent to any vertices of color 3.  As such, if we remove all edges containing it, we do not lose any facets.  This makes there no longer an extra vertex corresponding to the color 12, so we are back in a previous case.  The same analysis applies if $c_5 \geq c_2$ or $c_4 \geq c_6$.  This leaves only the case where $c_1 < c_3$, $c_5 < c_2$, and $c_4 < c_6$.  These conditions imply that no two of the extra vertices are adjacent.

Construct a simplicial complex $\Delta_2$ from $\Delta_1$ by discarding all vertices and edges not contained in a facet of $\Delta_1$.  In the new complex $\Delta_2$, let $v_{g_1+1}^1$ be adjacent to $d_{12}$ vertices of color 2 and $d_{13}$ vertices of color 3, let $v_{g_2+1}^2$ be adjacent to $d_{21}$ vertices of color 1 and $d_{23}$ vertices of color 3, and let $v_{g_3+1}^3$ be adjacent to $d_{31}$ vertices of color 1 and $d_{32}$ vertices of color 2.  Because no two extra vertices are adjacent, any two vertices of distinct colors adjacent to one of these extra vertices must be adjacent to each other.  As such, $v_{g_1+1}^1$ is contained in $d_{12}d_{13}$ facets, $v_{g_2+1}^2$ is contained in $d_{21}d_{23}$ facets, and $v_{g_3+1}^3$ is contained in $d_{31}d_{32}$ facets.  Our construction gives that $d_{21} = g_1$, $d_{32} = g_2$, and $d_{13} = g_3$, which we can substitute into the above formulas.  No two of the extra vertices are adjacent to each other, so if we add the facets not containing any of the extra vertices, we can compute
$$f_{123}(\Delta_2) = g_1g_2g_3 + d_{12}d_{13} + d_{21}d_{23} + d_{31}d_{32} = g_1g_2g_3 + d_{12}g_3 + d_{23}g_1 + d_{31}g_2.$$

Suppose that $g_1 \geq g_2$.  We wish to compute $f_{12}(\Delta_2)$.  There are $g_1g_2$ edges containing neither of the extra vertices.  There are $d_{12}$ edges containing the extra vertex of color 1.  There are $d_{21} = g_1$ edges containing the extra vertex of color 2.  Add the edges containing neither of the extra vertices and we get $$f_{12}(\Delta_2) = g_1g_2 + d_{12} + d_{21} = g_1g_2 + d_{12} + g_1.$$

We create $\Gamma$ from $\Delta_2$ by rearranging the edges of color 12 as follows.  Remove the edges containing $v_{g_1+1}^1$ or $v_{g_2+1}^2$.  In their place, make $v_{g_1+1}^1$ adjacent to the first $g_2$ vertices of color 2. Use the remaining edges to make $v_{g_2+1}^2$ adjacent to the first $g_1 + d_{12} - g_2$ edges of color 1.  We can do this because $g_1 \geq g_2$ gives that $g_1 + d_{12} - g_2 \geq d_{12} > 0$ and $d_{12} \leq g_2$ gives that $g_1 + d_{12} - g_2 \leq g_1$, so there are enough vertices of color 1.  This keeps all flag f-numbers of $\Gamma$ the same as of $\Delta_2$ except possibly the number of facets.  If we plug the new values of $d_{12}$ and $d_{21}$ into the above formula, we get $$f_{123}(\Gamma) = g_1g_2g_3 + g_2g_3 + d_{23}(g_1 + d_{12} - g_2) + d_{31}g_2.$$
Now we can compute
\begin{eqnarray*}
f_{123}(\Gamma) - f_{123}(\Delta_2) & = & g_1g_2g_3 + g_2g_3 + d_{23}(g_1 + d_{12} - g_2) + d_{31}g_2 \\ & & \qquad - (g_1g_2g_3 + d_{12}g_3 + d_{23}g_1 + d_{31}g_2) \\ & = & (g_2 - d_{12})g_3 + d_{23}(d_{12} - g_2) \\ & = & (g_2 - d_{12})(g_3 - d_{23}) \geq 0
\end{eqnarray*}
The last line follows because both factors are nonnegative by construction.

We must check that $\Gamma$ satisfies all of the needed conditions.  First, for all $S \subset [3]$ with $S \not = [3]$, we have $f_S(\Gamma) = f_S(\Delta_2) \leq f_S(\Delta_1) = f_S(\Delta)$.  For condition 2, we have just shown $f_{123}(\Gamma) \geq f_{123}(\Delta_2) = f_{123}(\Delta_1) \geq f_{123}(\Delta)$.  For condition 4, the vertex $v_{g_1+1}^1$ is now adjacent to the first $g_2$ vertices of color 2 and the first $g_3$ vertices of color 3.  As such, we it is no longer an extra vertex, and we can increase $g_1$ by 1.  This leaves only the other two vertices as extra vertices.  Condition 3 is clear.  Finally, $\Gamma$ comes from the specified construction with the new value of $g_1$, $p = 3$, and $q = 2$.

If $g_2 \geq g_3$, we can do the same procedure as before, this time rearranging edges of color 23 to make $v_{g_2+1}^2$ no longer an extra vertex.  Likewise, if $g_3 \geq g_1$, we rearrange edges of color 13 to make $v_{g_3+1}^3$ no longer an extra vertex.  This leaves only the case where $g_1 < g_2 < g_3 < g_1$, which is impossible.  \endproof

With this last lemma, finding the greatest number of facets possible given the rest of the flag f-vector merely requires finding the optimal values of $g_1$, $g_2$, $g_3$, $p$ and $q$.  What remains is to restrict how many values of these constants are necessary to check in order to ensure that we have found the maximum number of facets.

In subsequent lemmas, we sometimes have more than one complex constructed as $\Gamma$ was in the previous lemma, except using different constants.  To avoid confusion, we refer to the constants associated to a particular complex as $g_1(\Gamma), p(\Gamma)$, and so forth.

\begin{definition}
\textup{Let $\Delta$ be a simplicial complex.  We define}
\begin{eqnarray*}
\mathcal{A}(\Delta) & = & \{\Gamma\ |\ \Gamma\ \textup{satisfies conditions 1, 3, and 4 of Lemma~\ref{twovert} and is} \\ & & \textup{constructed as described in the statement of that lemma} \}, \\ m(\Delta) & = & \max\{f_{123}(\Gamma)\ |\ \Gamma \in \mathcal{A}(\Delta)\}, \\ \mathcal{B}(\Delta) & = & \{\Gamma\ \in \mathcal{A}(\Delta)\ |\ f_{123}(\Gamma) = m(\Delta)\}, \\ n(\Delta) & = & \max\{f_{12}(\Gamma) + f_{13}(\Gamma) + f_{23}(\Gamma)\ |\ \Gamma \in \mathcal{B}(\Delta)\} \\ \mathcal{C}(\Delta) & = & \{\Gamma \in \mathcal{B}(\Delta)\ |\ f_{12}(\Gamma) + f_{13}(\Gamma) + f_{23}(\Gamma) = n(\Delta)\}, \qquad \textup{and} \\ r(\Gamma) & = & 6 - p(\Gamma) - q(\Gamma) \qquad \textup{for any } \Gamma \in \mathcal{A}(\Delta).
\end{eqnarray*}
\end{definition}

$\mathcal{A}(\Delta)$ is the set of complexes constructed from $\Delta$ as described in Lemma~\ref{twovert}.  The lemma guarantees that there will be some $\Gamma \in \mathcal{A}(\Delta)$ such that $f_{123}(\Gamma) \geq f_{123}(\Delta)$, which means that $m(\Delta) \geq f_{123}(\Delta)$.  The definition allows that there could also be some $\Gamma \in \mathcal{A}(\Delta)$ such that $f_{123}(\Gamma) < f_{123}(\Delta)$.

The problem of characterizing the flag f-vectors of 3-colored complexes is equivalent to computing $m(\Delta)$.  It is immediate from the definitions that $\mathcal{B}(\Delta) \not = \emptyset$ and $\mathcal{C}(\Delta) \not = \emptyset$.

Intuitively, $p(\Gamma)$ and $q(\Gamma)$ are two of the numbers in $\{1, 2, 3\}$, and $r(\Gamma)$ is the third.

Next we have an easy lemma that will allow us to conveniently discard some complexes later.

\begin{lemma} \label{twovertmax}
Let $\Theta$ and $\Delta$ be simplicial complexes such that $f_S(\Theta) \leq f_S(\Delta)$ for all $S \subset [3]$ with $S \not = [3]$.  Then $f_{123}(\Theta) \leq m(\Delta)$.  Furthermore, if $\Sigma$ is a complex such that $f_{123}(\Sigma) < f_{123}(\Theta)$, then $\Sigma \not \in \mathcal{B}(\Delta)$.
\end{lemma}

\proof  By Lemma~\ref{twovert}, we can construct a simplicial complex $\Gamma$ from $\Theta$ satisfying all of the conditions of the lemma.  Since $f_S(\Gamma) \leq f_S(\Theta) \leq f_S(\Delta)$ for all $S \subset [3]$ with $S \not = [3]$, if $f_{123}(\Gamma) \geq f_{123}(\Delta)$, then $\Gamma$ satisfies all the conditions of Lemma~\ref{twovert} for $\Delta$ as well as for $\Theta$.  Thus, $\Gamma \in \mathcal{A}(\Delta)$, so $m(\Delta) \geq f_{123}(\Gamma) \geq f_{123}(\Theta)$.  If $f_{123}(\Gamma) < f_{123}(\Delta)$, then $m(\Delta) \geq f_{123}(\Delta) > f_{123}(\Gamma) \geq f_{123}(\Theta)$.  The second claim follows immediately from $f_{123}(\Sigma) < f_{123}(\Theta) \leq m(\Delta)$. \endproof

The next two lemmas assert that either computing $m(\Delta)$ is rather trivial or else complexes in $\mathcal{C}(\Delta)$ use as many edges as $\Delta$.

\begin{lemma} \label{alledges}
Let $\Delta$ be a 3-colored simplicial complex.  Then at least one of the following holds for every $\Gamma \in \mathcal{C}(\Delta)$.
\begin{enumerate}
\item $f_{123}(\Gamma) = f_1(\Delta)f_{23}(\Delta)$;
\item $f_{123}(\Gamma) = f_2(\Delta)f_{13}(\Delta)$;
\item $f_{123}(\Gamma) = f_3(\Delta)f_{12}(\Delta)$; or
\item $f_{12}(\Gamma) = f_{12}(\Delta)$, $f_{13}(\Gamma) = f_{13}(\Delta)$, and $f_{23}(\Gamma) = f_{23}(\Delta)$.
\end{enumerate}
\end{lemma}

\proof  Note first that if one of the options of the lemma holds for some $\Gamma \in \mathcal{C}(\Delta)$, then it must hold for all complexes in $\mathcal{C}(\Delta)$, as all have exactly the same number of facets (a requirement of $\mathcal{B}(\Delta)$) and total edges (a requirement of $\mathcal{C}(\Delta)$), and the only way for $\Gamma$ to have as many edges as $\Delta$ in total is for it to have exactly the same number of edges of each color set.  Thus, it suffices to prove the lemma for just one $\Gamma \in \mathcal{C}(\Delta)$.

\textbf{Case I:}  All three of the equalities in condition 4 of the lemma fail.

This means we can add an extra edge of each color set and still have $f_S(\Gamma) \leq f_S(\Delta)$ for all $S \subset [3]$ with $S \not = [3]$.

\textbf{Case I A:}  Some pair of vertices of distinct colors is not adjacent.

We can add an edge to connect this pair of vertices, and make both vertices adjacent to some vertex of the third color by adding an edge if necessary.  This adds another facet, so by Lemma~\ref{twovertmax}, $\Gamma \not \in \mathcal{B}(\Delta) \supseteq \mathcal{C}(\Delta)$, a contradiction.

\textbf{Case I B:}  Every pair of vertices of distinct colors is adjacent.

Because $f_1(\Delta)f_2(\Delta) \geq f_{12}(\Delta) > f_{12}(\Gamma) = f_1(\Gamma)f_2(\Gamma)$, either $f_1(\Delta) > f_1(\Gamma)$ or $f_2(\Delta) > f_2(\Gamma)$.  Assume without loss of generality that $f_1(\Delta) > f_1(\Gamma)$.  Add another vertex of color 1 and make it adjacent to a vertex of each other color to add a facet.  Hence, Lemma~\ref{twovertmax} gives $\Gamma \not \in \mathcal{B}(\Delta) \supseteq \mathcal{C}(\Delta)$, a contradiction.

\textbf{Case II:}  Exactly two of the equalities in condition 4 of the lemma fail.

Assume without loss of generality that $f_{12}(\Gamma) < f_{12}(\Delta)$ and $f_{13}(\Gamma) < f_{13}(\Delta)$.

\textbf{Case II A:}  $f_1(\Gamma) < f_1(\Delta)$

We can construct a new complex $\Gamma_1$ from $\Gamma$ by adding another vertex of color 1 to $\Gamma$ and making it adjacent to at least one vertex of each of the other colors.  This increases the number of facets, so $f_{123}(\Gamma_1) > f_{123}(\Gamma)$.  By Lemma~\ref{twovertmax}, $m(\Delta) \geq f_{123}(\Gamma_1) > f_{123}(\Gamma)$, so $\Gamma \not \in \mathcal{B}(\Delta) \supseteq \mathcal{C}(\Delta)$, a contradiction.

\textbf{Case II B:}  $f_1(\Gamma) = f_1(\Delta)$

An edge of color 23 can be contained in at most $f_1(\Delta)$ facets of $\Gamma$, as a facet is uniquely determined by the choice of an edge of color 23 and a vertex of color 1.

\textbf{Case II B 1:}  Every edge of color 23 is contained in exactly $f_1(\Delta)$ facets.

That only two of the equalities of condition 4 fail means that $f_{23}(\Gamma) = f_{23}(\Delta)$.  Since $f_{123}(\Gamma) = f_1(\Delta)f_{23}(\Delta)$, option 1 in the lemma holds.

\textbf{Case II B 2:}  There is some edge of color 23 contained in fewer than $f_1(\Delta)$ facets of $\Gamma$.

Let the edge in question be $\{v_i^2, v_j^3\}$.  Since $\Gamma$ is color-shifted, these two vertices together with $v_{f_1(\Delta)}^1$ do not form a facet of $\Gamma$.  Add edges as necessary to make $v_i^2$ and $v_j^3$ adjacent to $v_{f_1(\Delta)}^1$.  This adds an extra facet, and we had a spare edge available of both of the relevant color sets.  Thus, by Lemma~\ref{twovertmax}, $\Gamma \not \in \mathcal{B}(\Delta) \supseteq \mathcal{C}(\Delta)$, a contradiction.

\textbf{Case III:}  Exactly one of the equalities in condition 4 of the lemma fails.

Assume without loss of generality that $f_{12}(\Gamma) < f_{12}(\Delta)$.

\textbf{Case III A:}  $q(\Gamma) = 3$

Assume without loss of generality that $p(\Gamma) = 1$.

\textbf{Case III A 1:}  $g_1(\Gamma) = f_1(\Delta)$

This means that applying the $p(\Gamma) = 1$ step doesn't add any faces, as a vertex can't be added.  Since $f_1(\Delta)f_2(\Delta) \geq f_{12}(\Delta) > f_{12}(\Gamma) = f_1(\Gamma)f_2(\Gamma)$, we have $f_2(\Gamma) < f_2(\Delta)$.  Define $\Gamma_1$ by $g_i(\Gamma_1) = g_i(\Gamma)$ for all $i \in [3]$, $p(\Gamma_1) = 3$, and $q(\Gamma_1) = 2$.  After adding the first extra vertex, we have $\Gamma$ exactly.  The second extra vertex uses at least one additional edge of color 12, so $f_{12}(\Gamma) + f_{13}(\Gamma) + f_{23}(\Gamma) < f_{12}(\Gamma_1) + f_{13}(\Gamma_1) + f_{23}(\Gamma_1) \leq n(\Delta)$.  Therefore, $\Gamma \not \in \mathcal{C}(\Delta)$, a contradiction.

\textbf{Case III A 2:}  $g_1(\Gamma) < f_1(\Delta)$

\textbf{Case III A 2 a:}  $f_{13}(\Gamma) \geq (g_1(\Gamma)+1)g_3(\Gamma)$

The extra vertex of color 1 is adjacent to all $g_2(\Gamma)$ vertices of color 2 (since we must have leftover edges of color 12) and at least $g_3(\Gamma)$ vertices of color 3, as the construction requires making $v_{g_1(\Gamma)+1}^1$ adjacent to as many other vertices as the restrictions on edges allow.

\textbf{Case III A 2 a i:}  $g_1(\Gamma) + 1 < f_1(\Delta)$

We can create a new complex $\Gamma_1$ as in the construction of Lemma~\ref{twovert} using $g_1(\Gamma_1) = g_1(\Gamma)+1$, $g_2(\Gamma_1) = g_2(\Gamma)$, $g_3(\Gamma_1) = g_3(\Gamma)$, $p(\Gamma_1) = 3$, and $q(\Gamma_1) = 1$.  After adding the first extra vertex, we have the complex $\Gamma$ exactly.  Adding the second extra vertex uses at least one additional edge of color 12, while $f_{123}(\Gamma_1) = f_{123}(\Gamma)$.  Thus, $f_{12}(\Gamma) + f_{13}(\Gamma) + f_{23}(\Gamma) < f_{12}(\Gamma_1) + f_{13}(\Gamma_1) + f_{23}(\Gamma_1) \leq n(\Delta)$, so $\Gamma \not \in \mathcal{C}(\Delta)$, a contradiction.

\textbf{Case III A 2 a ii:}  $g_1(\Gamma) + 1 = f_1(\Delta)$

We can create a new complex $\Gamma_1$ as in the construction of Lemma~\ref{twovert} using $g_1(\Gamma_1) = g_1(\Gamma)+1$, $g_2(\Gamma_1) = g_2(\Gamma)$, $g_3(\Gamma_1) = g_3(\Gamma)$, $p(\Gamma_1) = 3$, and $q(\Gamma_1) = 2$.  After adding the first extra vertex, we have the complex $\Gamma$ exactly.  Since $f_1(\Delta)f_2(\Delta) \geq f_{12}(\Delta) > f_{12}(\Gamma) = f_1(\Gamma)f_2(\Gamma)$ and $f_1(\Gamma) = g_1(\Gamma) + 1 = f_1(\Delta)$, we have $f_2(\Gamma) < f_2(\Delta)$.  Adding the second extra vertex uses at least one additional edge of color 12, while $f_{123}(\Gamma_1) = f_{123}(\Gamma)$.  Thus, $f_{12}(\Gamma) + f_{13}(\Gamma) + f_{23}(\Gamma) < f_{12}(\Gamma_1) + f_{13}(\Gamma_1) + f_{23}(\Gamma_1) \leq n(\Delta)$, so $\Gamma \not \in \mathcal{C}(\Delta)$, a contradiction.

\textbf{Case III A 2 b:}  $f_{13}(\Gamma) < (g_1(\Gamma)+1)g_3(\Gamma)$

Adding the vertex $v_{g_1(\Gamma)+1}^1$ because $p(\Gamma) = 1$ uses up all of the remaining edges of color 13.  Hence, $v_{g_3(\Gamma)+1}^3$, which is added because $q(\Gamma) = 3$, cannot be adjacent to any vertices of color 1, and therefore is not in any facets.

\textbf{Case III A 2 b i:}  $f_2(\Delta) = g_2(\Gamma)$

All edges of color 13 are of the form $\{v_i^1, v_j^3\}$ with $i \leq g_1(\Gamma)+1$ and $j \leq g_3(\Gamma)$.  By construction, every vertex of color 2 is adjacent to all of the vertices in $\{v_1^1, \dots, v_{g_1(\Gamma)+1}^1, v_1^3, \dots, v_{g_3(\Gamma)}^3\}$.  Thus, every choice of an edge of color 13 and a vertex of color 2 forms a facet.  Therefore, $f_{123}(\Gamma) = f_2(\Gamma)f_{13}(\Gamma) = f_2(\Delta)f_{13}(\Delta)$, which is option 2 in the lemma.

\textbf{Case III A 2 b ii:}  $f_2(\Delta) > g_2(\Gamma)$

Define a complex $\Gamma_1$ by $g_i(\Gamma_1) = g_i(\Gamma)$ for all $i \in [3]$, $p(\Gamma_1) = 1$, and $q(\Gamma_1) = 2$.

\textbf{Case III A 2 b ii (1):}  $f_{23}(\Delta) > g_2(\Gamma)g_3(\Gamma)$

The difference in facets between $\Gamma$ and $\Gamma_1$ is the number added by the second extra vertex.  This vertex is not contained in any facets of $\Gamma$.  The new vertex $v_{g_2(\Gamma)+1}^2$ adjacent in $\Gamma_1$ to vertices of both color 1 and color 3, and is thus contained in at least one facet.  Therefore, by Lemma~\ref{twovertmax}, $\Gamma \not \in \mathcal{B}(\Delta) \supseteq \mathcal{C}(\Delta)$.

\textbf{Case III A 2 b ii (2):}  $f_{23}(\Delta) = g_2(\Gamma)g_3(\Gamma)$

All edges of color 23 were in $\Gamma$ before adding either extra vertex.  All edges of color 13 were in $\Gamma$ after adding the first extra vertex (for $p(\Gamma) = 1$) and before adding the second.  As such, adding the extra vertex of color 3 in $\Gamma$ does not use any more edges, nor does it add any facets.  The second extra vertex of $\Gamma_1$ does at least use some additional edges of color 12.  As such, $f_{123}(\Gamma) = f_{123}(\Gamma_1)$ and $f_{12}(\Gamma) + f_{13}(\Gamma) + f_{23}(\Gamma) < f_{12}(\Gamma_1) + f_{13}(\Gamma_1) + f_{23}(\Gamma_1) \leq n(\Delta)$, so $\Gamma \not \in \mathcal{C}(\Delta)$.

\textbf{Case III B:}  $p(\Gamma) = 3$

Assume without loss of generality that $q(\Gamma) = 1$.

\textbf{Case III B 1:}  $g_1(\Gamma) = f_1(\Delta)$

Making the second extra vertex of color 2 rather than color 1 increases the number of edges used by the same argument as in Case III A 1, so $\Gamma \not \in \mathcal{C}(\Delta)$.

\textbf{Case III B 2:}  $g_1(\Gamma) < f_1(\Delta)$

Define $\Gamma_1$ by $g_i(\Gamma_1) = g_i(\Gamma)$ for all $i \in [3]$, $p(\Gamma_1) = 1$, and $q(\Gamma_1) = 3$.  The only edges that can differ between $\Gamma$ and $\Gamma_1$ are that $\Gamma$ has some extra edges containing $v_{g_3(\Gamma)+1}^3$ but not $v_{g_1(\Gamma)+1}^1$, while $\Gamma_1$ has the same number of extra edges containing $v_{g_1(\Gamma)+1}^1$ but not $v_{g_3(\Gamma)+1}^3$.  The only vertex that may not be adjacent to all $g_2(\Gamma)$ vertices of color 2 is $v_{g_3(\Gamma)+1}^3$.  Thus, any edge in $\Gamma_1$ but not $\Gamma$ is in at least as many facets as each edge in $\Gamma$ but not $\Gamma_1$.  Therefore, $f_{123}(\Gamma_1) \geq f_{123}(\Gamma)$.  Since the rest of their flag f-vectors are the same, if $\Gamma \in \mathcal{C}(\Delta)$, then $\Gamma_1 \in \mathcal{C}(\Delta)$.  We now apply Case~III~A to $\Gamma_1$.

\textbf{Case III C:}  $r(\Gamma) = 3$

Because $f_{12}(\Delta) > f_{12}(\Gamma)$, all vertices of color 1 are adjacent to all vertices of color 2, including the extra vertex of each color.  Thus, $\Gamma$ is exactly the same complex regardless of whether $p(\Gamma) = 1$ and $q(\Gamma) = 2$ or vice versa.

Note that using all edges of colors 13 and 23 means that $f_{13}(\Gamma) \leq (g_1(\Gamma) + 1)g_3(\Gamma)$ and $f_{23}(\Gamma) \leq (g_2(\Gamma) + 1)g_3(\Gamma)$.

\textbf{Case III C 1:}  $f_1(\Delta) = g_1(\Gamma)$ or $f_2(\Delta) = g_2(\Gamma)$

Assume without loss of generality that $f_1(\Delta) = g_1(\Gamma)$.  Since $f_{13}(\Delta) = f_{13}(\Gamma) \geq g_1(\Gamma)g_3(\Gamma)$ and there are leftover edges of color 12, every vertex of color 1 is adjacent to every vertex not of color 1.  Therefore, $f_{123}(\Gamma) = f_1(\Gamma)f_{23}(\Gamma) = f_1(\Delta)f_{23}(\Delta)$, which is option 1 of the lemma.

\textbf{Case III C 2:}  $f_1(\Delta) > g_1(\Gamma)$ and $f_2(\Delta) > g_2(\Gamma)$

\textbf{Case III C 2 a:}  $f_{13}(\Gamma) = (g_1(\Gamma) + 1)g_3(\Gamma)$ or $f_{23}(\Gamma) = (g_2(\Gamma) + 1)g_3(\Gamma)$

Assume without loss of generality that $f_{13}(\Gamma) = (g_1(\Gamma) + 1)g_3(\Gamma)$.  As noted above, we can assume that $p(\Gamma) = 1$.

\textbf{Case III C 2 a i:}  $f_1(\Gamma) = f_1(\Delta)$

Every vertex of color 1 is adjacent to every vertex of color 2 or 3, so it forms a facet together with every edge of color 23.  Therefore, $f_{123}(\Gamma) = f_1(\Gamma)f_{23}(\Gamma) = f_1(\Delta)f_{23}(\Delta)$, which is option 1 of the lemma.

\textbf{Case III C 2 a ii:}  $f_1(\Gamma) < f_1(\Delta)$

Define $\Gamma_1$ by $g_1(\Gamma_1) = g_1(\Gamma) + 1$, $g_2(\Gamma_1) = g_2(\Gamma)$, $g_3(\Gamma_1) = g_3(\Gamma)$, $p(\Gamma_1) = 2$, and $q(\Gamma_1) = 1$.  In this case, $\Gamma_1$ consists of $\Gamma$ plus an extra vertex of color 1 that is contained in some edges of color 12.  Thus, $\Gamma_1$ has all the facets of $\Gamma$ and more edges, so $\Gamma \not \in \mathcal{C}(\Delta)$.

\textbf{Case III C 2 b:}  $f_{13}(\Gamma) < (g_1(\Gamma) + 1)g_3(\Gamma)$, $f_{23}(\Gamma) < (g_2(\Gamma) + 1)g_3(\Gamma)$, and either $f_{13}(\Gamma) \geq (g_1(\Gamma) + 1)(g_3(\Gamma) - 1)$ or $f_{23}(\Gamma) \geq (g_2(\Gamma) + 1)(g_3(\Gamma) - 1)$

Assume without loss of generality that $f_{13}(\Gamma) \geq (g_1(\Gamma) + 1)(g_3(\Gamma) - 1)$.  Again, we can assume that $p(\Gamma) = 1$.  Define $\Gamma_1$ by $g_1(\Gamma_1) = g_1(\Gamma) + 1$, $g_2(\Gamma_1) = g_2(\Gamma)$, $g_3(\Gamma_1) = g_3(\Gamma) - 1$, $p(\Gamma_1) = 3$, and $q(\Gamma_1) = 2$.  It is easy to check that $\Gamma_1$ has exactly the same edges of colors 12 and 23 as $\Gamma$.  The only possible difference in edges is that $\Gamma_1$ may have some extra edges containing $v_{g_1(\Gamma)+1}^1$ but not $v_{g_3(\Gamma)}^3$ while $\Gamma$ may have some extra edges containing $v_{g_3(\Gamma)}^3$ but not $v_{g_1(\Gamma)+1}^1$.  An edge of the former type is contained in at least $g_2(\Gamma)$ facets, as the first $g_2(\Gamma)$ vertices of color 2 are adjacent to all vertices of other colors.  An edge of the latter type is contained in at most $g_2(\Gamma)$ facets because $v_{g_3(\Gamma)}^3$ is adjacent to only $g_2(\Gamma)$ vertices of color 2, as $f_{23}(\Gamma) < (g_2(\Gamma) + 1)g_3(\Gamma)$.  Therefore $f_{123}(\Gamma_1) \geq f_{123}(\Gamma)$.  Since $\Gamma_1$ has the same number of edges as $\Gamma$, we get $\Gamma_1 \in \mathcal{C}(\Delta)$.  We now apply Case~III~B to $\Gamma_1$.

\textbf{Case III C 2 c:}  $f_{13}(\Gamma) < (g_1(\Gamma) + 1)(g_3(\Gamma) - 1)$ and $f_{23}(\Gamma) < (g_2(\Gamma) + 1)(g_3(\Gamma) - 1)$

Define $\Gamma_1$ by $g_1(\Gamma_1) = g_1(\Gamma)$, $g_2(\Gamma_1) = g_2(\Gamma)$, $g_3(\Gamma_1) = g_3(\Gamma) - 1$, $p(\Gamma_1) = 1$, and $q(\Gamma_1) = 2$.  The facets of $\Gamma$ not in $\Gamma_1$ are those containing $v_{g_3(\Gamma)}^3$.  There are $g_1(\Gamma)g_2(\Gamma)$ such facets.  The facets of $\Gamma_1$ not in $\Gamma$ are those containing the edges of $\Gamma_1$ but not $\Gamma$.  Both of the vertices of each such edge of color 13 are adjacent to at least the first $g_2(\Gamma)$ vertices of color 2.  Each of the $g_1(\Gamma)$ such edges adds at least $g_2(\Gamma)$ facets.  Likewise, the vertices of a new edge of color 23 are adjacent to at least the first $g_1(\Gamma)$ vertices of color 1, so the edge is contained in at least $g_1(\Gamma)$ facets.  Thus, $\Gamma_1$ contains at least $2g_1(\Gamma)g_2(\Gamma)$ facets that $\Gamma_2$ does not.  Therefore, $f_{123}(\Gamma_1) > f_{123}(\Gamma)$, so $\Gamma \not \in \mathcal{B}(\Delta) \supseteq \mathcal{C}(\Delta)$.

\textbf{Case IV:}  All of the equalities of condition 4 of the lemma hold.

This makes condition 4 of the lemma true.  \endproof

If any of the first three options in the preceding lemma apply, then we are done, as these are also upper bounds.  For example, each facet corresponds uniquely to a vertex of color 1 together with an edge of color 23, so $f_{123}(\Delta) \leq f_1(\Delta)f_{23}(\Delta)$.  The next lemma says that it is easy to tell when one of the first three options hold.

\begin{lemma} \label{vertedge}
Let $\Delta$ be a simplicial complex and let $\Gamma \in \mathcal{C}(\Delta)$.
\begin{enumerate}
\item $f_{123}(\Gamma) = f_1(\Delta)f_{23}(\Delta)$ if and only if $\Big\lfloor {f_{12}(\Delta) \over f_1(\Delta)} \Big\rfloor \Big\lfloor {f_{13}(\Delta) \over f_1(\Delta)} \Big\rfloor \geq f_{23}(\Delta)$.
\item $f_{123}(\Gamma) = f_2(\Delta)f_{13}(\Delta)$ if and only if $\Big\lfloor {f_{12}(\Delta) \over f_2(\Delta)} \Big\rfloor \Big\lfloor {f_{23}(\Delta) \over f_2(\Delta)} \Big\rfloor \geq f_{13}(\Delta)$.
\item $f_{123}(\Gamma) = f_3(\Delta)f_{12}(\Delta)$ if and only if $\Big\lfloor {f_{13}(\Delta) \over f_3(\Delta)} \Big\rfloor \Big\lfloor {f_{23}(\Delta) \over f_3(\Delta)} \Big\rfloor \geq f_{12}(\Delta)$.
\end{enumerate}
\end{lemma}

\proof  It suffices to prove one, as the others follow by relabeling the colors.  Suppose first that $\Big\lfloor {f_{12}(\Delta) \over f_1(\Delta)} \Big\rfloor \Big\lfloor {f_{13}(\Delta) \over f_1(\Delta)} \Big\rfloor \geq f_{23}(\Delta)$.  Let $\Gamma_1$ be defined by $g_1(\Gamma_1) = f_1(\Delta)$, $g_2(\Gamma_1) = \Big\lfloor {f_{12}(\Delta) \over f_1(\Delta)} \Big\rfloor$, $g_3(\Gamma_1) = \Big\lfloor {f_{23}(\Delta) \over g_2(\Gamma_1)} \Big\rfloor$, $p(\Gamma_1) = 3$, and $q(\Gamma_1) = 2$.  It follows that $g_3(\Gamma_1) > {f_{23}(\Delta) \over g_2(\Gamma_1)} - 1$, from which $(g_3(\Gamma_1) + 1)g_2(\Gamma_1) > f_{23}(\Delta)$.  Thus, all edges of color 23 have their vertices in the set $\{v_1^2, \dots, v_{g_2(\Gamma_1)}^2, v_1^3, \dots, v_{g_3(\Gamma_1)+1}^3\}$.

Since $\Big\lfloor {f_{12}(\Delta) \over f_1(\Delta)} \Big\rfloor \Big\lfloor {f_{13}(\Delta) \over f_1(\Delta)} \Big\rfloor \geq f_{23}(\Delta)$, we have that $g_2(\Gamma_1) \Big\lfloor {f_{13}(\Delta) \over f_1(\Delta)} \Big\rfloor \geq f_{23}(\Delta)$, and so $\Big\lfloor {f_{13}(\Delta) \over f_1(\Delta)} \Big\rfloor \geq {f_{23}(\Delta) \over g_2(\Gamma_1)}$.  If the right side is not an integer, then $\Big\lfloor {f_{13}(\Delta) \over f_1(\Delta)} \Big\rfloor \geq \Big\lfloor {f_{23}(\Delta) \over g_2(\Gamma_1)} \Big\rfloor + 1 = g_3(\Gamma_1) + 1$, so all vertices contained in edges of color 23 are adjacent to all vertices of color 1.  If ${f_{23}(\Delta) \over g_2(\Gamma_1)}$ is an integer, then $v_{g_3(\Gamma_1)+1}^3$ is not contained in any edges of color 23, and again, all vertices contained in edges of color 23 are adjacent to all vertices of color 1.  Either way, every edge of color 23 forms a facet together with each vertex of color 1, so $f_{123}(\Gamma_1) = f_1(\Gamma_1)f_{23}(\Gamma_1) = f_1(\Delta)f_{23}(\Delta)$.  By Lemma~\ref{twovertmax}, for any $\Gamma \in \mathcal{C}(\Delta)$, $f_{123}(\Gamma) \geq f_{123}(\Gamma_1) = f_1(\Delta)f_{23}(\Delta)$.  Since we also have $f_{123}(\Gamma) \leq f_1(\Gamma)f_{23}(\Gamma) \leq f_1(\Delta)f_{23}(\Delta)$, the result follows.

Conversely, suppose that $f_{123}(\Gamma) = f_1(\Delta)f_{23}(\Delta)$.  We must have $f_{23}(\Gamma) = f_{23}(\Delta)$, and every edge of $\Gamma$ of color 23 must form a facet with each of the $f_1(\Delta)$ vertices of color 1.  Thus, every vertex of an edge of color 23 must be adjacent to every vertex of color 1.  If there are $d_2$ such vertices of color 2 and $d_3$ such vertices of color 3, then the number of required edges is $f_1(\Delta)d_2$ of color 12 and $f_1(\Delta)d_3$ of color 13.  Since we are only allowed so many edges of each color set, we have $f_1(\Delta)d_2 \leq f_{12}(\Delta)$ and $f_1(\Delta)d_3 \leq f_{13}(\Delta)$.  These yield $d_2 \leq {f_{12}(\Delta) \over f_1(\Delta)}$ and $d_3 \leq {f_{12}(\Delta) \over f_1(\Delta)}$, respectively.  Since $d_2$ and $d_3$ must be integers, the inequality still holds if we take the integer parts of the right sides.  This yields $d_2 \leq \Big\lfloor {f_{12}(\Delta) \over f_1(\Delta)} \Big\rfloor$ and $d_3 \leq \Big\lfloor {f_{12}(\Delta) \over f_1(\Delta)} \Big\rfloor$, respectively.  The number of edges of color 23 on $d_2$ vertices of color 2 and $d_3$ vertices of color 3 is at most $d_2d_3$, so we have $f_{23}(\Delta) \leq d_2d_3 \leq \Big\lfloor {f_{12}(\Delta) \over f_1(\Delta)} \Big\rfloor \Big\lfloor {f_{12}(\Delta) \over f_1(\Delta)} \Big\rfloor$, as desired.  \endproof

If the preceding lemma does not give the value of $m(\Delta)$, then option 4 of Lemma~\ref{alledges} must hold.  In this case, we have some reasonably strong restrictions on what some $\Gamma \in \mathcal{C}(\Delta)$ must look like from Lemmas~\ref{twovert} and~\ref{alledges}.  The next task is to derive enough restrictions on $g_1(\Gamma)$, $g_2(\Gamma)$, $g_3(\Gamma)$, $p(\Gamma)$, and $q(\Gamma)$ for there to be only a small number of constructions to try in order to compute $m(\Delta)$.

\begin{definition}
Let $\Delta$ be a 3-colored simplicial complex.  Define $$\mathcal{D}(\Delta) = \{\Gamma \in \mathcal{C}(\Delta)\ |\ f_{123}(\Gamma) < \min\{f_1(\Delta)f_{23}(\Delta), f_2(\Delta)f_{13}(\Delta), f_3(\Delta)f_{12}(\Delta)\}\}.$$
If $\Gamma \in \mathcal{D}(\Delta)$, then define
\begin{eqnarray*}
j_1(\Gamma) & = & f_{23}(\Delta) - g_2(\Gamma)g_3(\Gamma) \\ j_2(\Gamma) & = & f_{13}(\Delta) - g_1(\Gamma)g_3(\Gamma) \\ j_3(\Gamma) & = & f_{12}(\Delta) - g_1(\Gamma)g_2(\Gamma).
\end{eqnarray*}
\end{definition}

If Lemma~\ref{vertedge} gives $m(\Delta)$, then $\mathcal{D}(\Delta) = \emptyset$.  Otherwise, $\mathcal{D}(\Delta) = \mathcal{C}(\Delta)$.  The $j_i's$ are the number of edges left over of a given color set before adding the additional vertices for $p(\Gamma)$ and $q(\Gamma)$.

This next lemma tells us when swapping $p(\Gamma)$ with $q(\Gamma)$ is beneficial.  Intuitively, this means switching the order in which the two extra vertices are added.

\begin{lemma} \label{choosepq}
Let $\Delta$ be a 3-colored simplicial complex and let $\Gamma \in \mathcal{A}(\Delta)$ with $j_{p(\Gamma)}(\Gamma) \geq j_{q(\Gamma)}(\Gamma)$.  Define $\Gamma_1$ by $g_i(\Gamma_1) = g_i(\Gamma)$ for all $i \in [3]$, $p(\Gamma_1) = q(\Gamma)$, and $q(\Gamma_1) = p(\Gamma)$.  Then $\Gamma_1 \in \mathcal{A}(\Delta)$ and $f_{123}(\Gamma_1) \geq f_{123}(\Gamma)$.  Furthermore, if $\Gamma \in \mathcal{D}(\Delta)$, then $\Gamma_1 \in \mathcal{D}(\Delta)$.
\end{lemma}

\proof  All of the edges that $g_1(\Gamma_1)$, $g_2(\Gamma_1)$, and $g_3(\Gamma_1)$ require $\Gamma_1$ to have are in $\Gamma$, so $\Delta$ has enough edges available for $\Gamma_1$ to be well-defined.  That $\Gamma_1 \in \mathcal{A}(\Delta)$ is immediate from the construction.  If $f_{123}(\Gamma_1) \geq f_{123}(\Gamma)$, then $f_{123}(\Gamma_1) \geq f_{123}(\Gamma) = m(\Delta)$, so $\Gamma_1 \in \mathcal{B}(\Delta)$.  Because $\Gamma_1$ has the same number of vertices of each color as $\Gamma$, it uses just as many edges of each color set as $\Gamma$.  If $\Gamma \in \mathcal{D}(\Delta)$, then $\Gamma \in \mathcal{C}(\Delta)$, so we have $\Gamma_1 \in \mathcal{C}(\Delta)$.  If $\Gamma \in \mathcal{D}(\Delta)$, then $\mathcal{D}(\Delta) \not = \emptyset$, and so $\mathcal{D}(\Delta) = \mathcal{C}(\Delta)$.  Thus, $\Gamma_1 \in \mathcal{D}(\Delta)$.  Therefore, it suffices to show that $f_{123}(\Gamma_1) \geq f_{123}(\Gamma)$.

It is immediate from the construction that both complexes have exactly the same edges of color sets $\{p(\Gamma), r(\Gamma)\}$ and $\{q(\Gamma), r(\Gamma)\}$.  All that can differ is the edges of color set $\{p(\Gamma), q(\Gamma)\}$.  Among these, all that can differ is that $\Gamma$ may have some extra edges containing $v_{g_{p(\Gamma)}+1}^{p(\Gamma)}$ but not $v_{g_{q(\Gamma)}+1}^{q(\Gamma)}$ while $\Gamma_1$ may have some extra edges containing $v_{g_{q(\Gamma)}+1}^{q(\Gamma)}$ but not $v_{g_{p(\Gamma)}+1}^{p(\Gamma)}$.  Any vertex of color $p(\Gamma)$ or $q(\Gamma)$ other than the two extra vertices is adjacent to exactly $g_{r(\Gamma)}(\Gamma)$ vertices of color $r(\Gamma)$ in both $\Gamma$ and $\Gamma_1$.  The vertex $v_{g_{p(\Gamma)}+1}^{p(\Gamma)}$ is adjacent to $j_{q(\Gamma)}(\Gamma) \leq g_{r(\Gamma)}(\Gamma)$ vertices of color $r(\Gamma)$.  The vertex $v_{g_{q(\Gamma)}+1}^{q(\Gamma)}$ is adjacent to $j_{p(\Gamma)}(\Gamma) \leq g_{r(\Gamma)}(\Gamma)$ vertices of color $r(\Gamma)$.  Hence, an extra edge of $\Gamma$ is contained in $j_{q(\Gamma)}(\Gamma)$ facets, while an extra edge of $\Gamma_1$ is contained in $j_{p(\Gamma)}(\Gamma)$ facets.  Since $j_{p(\Gamma)}(\Gamma) \geq j_{q(\Gamma)}(\Gamma)$, we have $f_{123}(\Gamma_1) \geq f_{123}(\Gamma)$.  \endproof

The next lemma says that, for some $\Gamma \in \mathcal{D}(\Delta)$, if we know $r(\Gamma)$ and $g_{r(\Gamma)}(\Gamma)$, then a particular construction is guaranteed to give us a $\Gamma_1 \in \mathcal{D}(\Delta)$. $\Gamma_1$ may or may not be the same as $\Gamma$.  This reduces the problem to trying to find these two parameters.

\begin{lemma}  \label{knowr}
Let $\Delta$ be a 3-colored simplicial complex.  Suppose that $\Gamma \in \mathcal{D}(\Delta)$, $p(\Gamma) = 1$, and $q(\Gamma) = 2$.  Define $\Gamma_1$ by $g_3(\Gamma_1) = g_3(\Gamma)$,
\begin{eqnarray*}
g_1(\Gamma_1) & = & \left\{ \begin{array}{ccl} \Big\lfloor {f_{13}(\Delta) \over g_3(\Gamma)} \Big\rfloor & \textup{if} & {f_{13}(\Delta) \over g_3(\Gamma)} \not \in \mathbb{Z} \\ {f_{13}(\Delta) \over g_3(\Gamma)} - 1 & \textup{if} & {f_{13}(\Delta) \over g_3(\Gamma)} \in \mathbb{Z} \textup{ and } {f_{13}(\Delta) \over g_3(\Gamma)} \Big\lceil {f_{23}(\Delta) \over g_3(\Gamma)} - 1 \Big\rceil > f_{12}(\Delta) \\ {f_{13}(\Delta) \over g_3(\Gamma)} & \textup{if} & {f_{13}(\Delta) \over g_3(\Gamma)} \in \mathbb{Z} \textup{ and } {f_{13}(\Delta) \over g_3(\Gamma)} \Big\lceil {f_{23}(\Delta) \over g_3(\Gamma)} - 1 \Big\rceil \leq f_{12}(\Delta) \end{array} \right., \\ g_2(\Gamma_1) & = & \left\{ \begin{array}{ccl} \Big\lfloor {f_{23}(\Delta) \over g_3(\Gamma)} \Big\rfloor & \textup{if} & {f_{23}(\Delta) \over g_3(\Gamma)} \not \in \mathbb{Z} \\ {f_{23}(\Delta) \over g_3(\Gamma)} - 1 & \textup{if} & {f_{23}(\Delta) \over g_3(\Gamma)} \in \mathbb{Z} \textup{ and } {f_{23}(\Delta) \over g_3(\Gamma)} g_1(\Gamma_1) > f_{12}(\Delta) \\ {f_{23}(\Delta) \over g_3(\Gamma)} & \textup{if} & {f_{23}(\Delta) \over g_3(\Gamma)} \in \mathbb{Z} \textup{ and } {f_{23}(\Delta) \over g_3(\Gamma)} g_1(\Gamma_1) \leq f_{12}(\Delta) \end{array} \right., \\ p(\Gamma_1) & = & \left\{ \begin{array}{ccl} 1 & \textup{if} & j_2(\Gamma_1) \geq j_1(\Gamma_1) \\ 2 & \textup{if} & j_2(\Gamma_1) < j_1(\Gamma_1) \end{array} \right., \qquad \textup{ and} \\ q(\Gamma_1) & = & 3 - p(\Gamma_1).
\end{eqnarray*}
Then $\Gamma_1 \in \mathcal{D}(\Delta)$.
\end{lemma}

\proof  We have that $g_3(\Gamma_1) = g_3(\Gamma)$.  If $g_1(\Gamma_1) = g_1(\Gamma)$ and $g_2(\Gamma_1) = g_2(\Gamma)$, then Lemma~\ref{choosepq} promises that $\Gamma_1 \in \mathcal{D}(\Delta)$.  Thus, for the lemma to be false, we must have either $g_1(\Gamma_1) \not = g_1(\Gamma)$ or $g_2(\Gamma_1) \not = g_2(\Gamma)$.

Since $\Gamma$ has no extra vertex of color 3, we must have $f_1(\Gamma) \geq {f_{13}(\Gamma) \over f_3(\Gamma)} = {f_{13}(\Delta) \over g_3(\Gamma)}$ and $f_2(\Gamma) \geq {f_{23}(\Gamma) \over f_3(\Gamma)} = {f_{23}(\Delta) \over g_3(\Gamma)}$.  One can check the various cases in the definition to easily see that $g_1(\Gamma_1) \leq f_1(\Gamma) \leq f_1(\Delta)$ and $g_2(\Gamma_1) \leq f_2(\Gamma) \leq f_2(\Delta)$.  We also have that $g_3(\Gamma_1) = g_3(\Gamma) \leq f_3(\Delta)$, so there are enough vertices for $\Gamma_1$ to be well-defined.

Because $\Gamma$ uses at least ${f_{13}(\Delta) \over g_3(\Gamma)} - 1$ vertices of color $1$ and at least ${f_{23}(\Delta) \over g_3(\Gamma)} - 1$ vertices of color $2$ before adding extra vertices, there are enough edges of color 12 to do this.  From the definition, $\Gamma_1$ does not require more vertices than this of either color unless there are enough edges.  Thus, $\Gamma_1$ is well-defined.  It is immediate from the definition that $\Gamma_1 \in \mathcal{A}(\Delta)$.

As the only extra vertex of color 13 that can contain edges is of color $1$, we have $j_2(\Gamma) \leq g_3(\Gamma)$.  Since $\Gamma \in \mathcal{D}(\Delta)$, it must use all edges of this color set, so we have $g_3(\Gamma)g_1(\Gamma) \leq f_{13}(\Gamma) \leq g_3(\Gamma)g_1(\Gamma) + g_3(\Gamma)= g_3(\Gamma)(g_1(\Gamma) + 1)$.  Divide by $g_3(\Gamma)$ and we have $g_1(\Gamma) \leq {f_{13}(\Delta) \over g_3(\Gamma)} \leq g_1(\Gamma) + 1$.  This can be rearranged as ${f_{13}(\Delta) \over g_3(\Gamma)} - 1 \leq g_1(\Gamma) \leq {f_{13}(\Delta) \over g_3(\Gamma)}$.  If ${f_{13}(\Delta) \over g_3(\Gamma)} \not \in \mathbb{Z}$, then this forces $g_1(\Gamma) = \Big\lfloor {f_{13}(\Delta) \over g_3(\Gamma)} \Big\rfloor = g_1(\Gamma_1)$.  Similarly, if ${f_{23}(\Delta) \over g_3(\Gamma)} \not \in \mathbb{Z}$, we get $g_2(\Gamma) = g_2(\Gamma_1)$.

\textbf{Case I:}  ${f_{13}(\Delta) \over g_3(\Gamma)} \not \in \mathbb{Z}$

As seen above, we have $g_1(\Gamma) = g_1(\Gamma_1)$.

\textbf{Case I A:}  ${f_{23}(\Delta) \over g_3(\Gamma)} \not \in \mathbb{Z}$

This case has $g_1(\Gamma_1) = g_1(\Gamma)$ and $g_2(\Gamma_1) = g_2(\Gamma)$.

\textbf{Case I B:}  ${f_{23}(\Delta) \over g_3(\Gamma)} \in \mathbb{Z}$

That $\Gamma$ is well-defined and uses all edges of color 13 corresponds to the inequalities $g_3(\Gamma)g_2(\Gamma) \leq f_{23}(\Gamma) \leq g_3(\Gamma)(g_2(\Gamma) + 1)$, which force either $g_2(\Gamma) = {f_{23}(\Delta) \over g_3(\Gamma)}$ or $g_2(\Gamma) = {f_{23}(\Delta) \over g_3(\Gamma)} - 1$.

\textbf{Case I B 1:}  $f_{12}(\Delta) < {f_{23}(\Delta) \over g_3(\Gamma)} g_1(\Gamma_1)$

There are not enough edges of color 12 to have $g_2(\Gamma) = {f_{23}(\Delta) \over g_3(\Gamma)}$ or $g_2(\Gamma_1) = {f_{23}(\Delta) \over g_3(\Gamma)}$.  This means that $g_2(\Gamma) = {f_{23}(\Delta) \over g_3(\Gamma)} - 1 = g_2(\Gamma_1)$, and we are done.

\textbf{Case I B 2:}  $f_{12}(\Delta) \geq {f_{23}(\Delta) \over g_3(\Gamma)} g_1(\Gamma_1)$

The definition of $\Gamma_1$ gives $g_2(\Gamma_1) = {f_{23}(\Delta) \over g_3(\Gamma)}$.  If we also have $g_2(\Gamma) = {f_{23}(\Delta) \over g_3(\Gamma)}$, then we are done.  Otherwise, $g_2(\Gamma) = {f_{23}(\Delta) \over g_3(\Gamma)} - 1$.  If this happens, we can compute $j_1(\Gamma) = g_3(\Gamma) > j_2(\Gamma) = j_2(\Gamma_1) > 0 = j_1(\Gamma_1)$ (with the strict inequalities because ${f_{13}(\Delta) \over g_3(\Gamma)} \not \in \mathbb{Z}$).  This means that $p(\Gamma_1) = 1$.  From this, Lemma~\ref{choosepq} asserts that if we define $\Gamma_2$ and by $g_i(\Gamma_2) = g_i(\Gamma)$ for all $i \in [3]$, $p(\Gamma_2) = 2$, and $q(\Gamma_2) = 1$, then $f_{123}(\Gamma_2) \geq f_{123}(\Gamma)$ and $\Gamma_2 \in \mathcal{D}(\Delta)$.

Note that $\Gamma_1$ is merely $\Gamma_2$ with possibly an extra isolated vertex added.  Thus, $f_{123}(\Gamma_1) \geq f_{123}(\Gamma_2)$, so $\Gamma_1 \in \mathcal{B}(\Delta)$.  Furthermore, since $\Gamma_2$ uses all available edges, so does $\Gamma_1$, and so $\Gamma_1 \in \mathcal{C}(\Delta)$.  Since $\Gamma_2 \in \mathcal{D}(\Delta)$, we get $\mathcal{D}(\Delta) = \mathcal{C}(\Delta)$, and so $\Gamma_1 \in \mathcal{D}(\Delta)$.

\textbf{Case II:}  ${f_{13}(\Delta) \over g_3(\Gamma)} \in \mathbb{Z}$

That $\Gamma$ is well-defined and uses all edges of color 13 corresponds to the inequalities $g_3(\Gamma)g_1(\Gamma) \leq f_{13}(\Gamma) \leq g_3(\Gamma)(g_1(\Gamma) + 1)$, which force either $g_1(\Gamma) = {f_{13}(\Delta) \over g_3(\Gamma)}$ or $g_1(\Gamma) = {f_{13}(\Delta) \over g_3(\Gamma)} - 1$.

\textbf{Case II A:}  ${f_{23}(\Delta) \over g_3(\Gamma)} \not \in \mathbb{Z}$

As we have seen, this gives $g_2(\Gamma) = g_2(\Gamma_1) = \Big\lfloor {f_{23}(\Delta) \over g_3(\Gamma)} \Big\rfloor$.

\textbf{Case II A 1:}  $f_{12}(\Delta) < {f_{13}(\Delta) \over g_3(\Gamma)} g_2(\Gamma)$

This gives $g_1(\Gamma_1) = {f_{13}(\Delta) \over g_3(\Gamma)} - 1$.  Since there are not enough edges for $\Gamma_1$ to have $g_1(\Gamma) = {f_{13}(\Delta) \over g_3(\Gamma)}$, we must have $g_1(\Gamma) = {f_{13}(\Delta) \over g_3(\Gamma)} - 1 = g_1(\Gamma_1)$.

\textbf{Case II A 2:}  $f_{12}(\Delta) \geq {f_{13}(\Delta) \over g_3(\Gamma)} g_2(\Gamma)$

This time, the definition gives $g_1(\Gamma_1) = {f_{13}(\Delta) \over g_3(\Gamma)}$.  In order to have $g_1(\Gamma_1) \not = g_1(\Gamma)$, we must have $g_1(\Gamma) = {f_{13}(\Delta) \over g_3(\Gamma)} - 1$.  This only swaps colors 1 and 2 from Case~I~B~2, so $\Gamma_1 \in \mathcal{D}(\Delta)$ by the same argument as there.

\textbf{Case II B:}  ${f_{23}(\Delta) \over g_3(\Gamma)} \in \mathbb{Z}$

\textbf{Case II B 1:}  $f_{12}(\Delta) \leq {f_{13}(\Delta) \over g_3(\Gamma)} {f_{23}(\Delta) \over g_3(\Gamma)}$

From the definition, it is clear that once the extra vertices are added, $\Gamma_1$ has at least ${f_{13}(\Delta) \over g_3(\Gamma)}$ vertices of color $1$ and ${f_{23}(\Delta) \over g_3(\Gamma)}$ vertices of color $2$.  Furthermore, the order of the extra vertices dictates that $\Gamma_1$ must have at least this many of each color before any more vertices of either color are added.  Thus, the $f_{12}(\Delta)$ edges of color 12 in $\Gamma_1$ all have both vertices among the first ${f_{13}(\Delta) \over g_3(\Gamma)}$ vertices of color $1$ and the first ${f_{23}(\Delta) \over g_3(\Gamma)}$ of color $2$.  All of these vertices are adjacent to all vertices of color $3$, so every edge of color 12 in $\Gamma_1$ is contained in $g_3(\Gamma)$ facets.  Therefore,
$$f_{123}(\Gamma_1) = g_3(\Gamma)f_{12}(\Delta) = f_3(\Gamma)f_{12}(\Gamma) \geq f_{123}(\Gamma).$$
As such, since $\Gamma \in \mathcal{B}(\Delta)$, we have $\Gamma_1 \in \mathcal{B}(\Delta)$.  Since $\Gamma_1$ has at least ${f_{13}(\Delta) \over g_3(\Gamma)}$ vertices of color $1$ and ${f_{23}(\Delta) \over g_3(\Gamma)}$ vertices of color $2$, it uses all available edges, and so $\Gamma_1 \in \mathcal{C}(\Delta) = \mathcal{D}(\Delta)$.

\textbf{Case II B 2:}  $f_{12}(\Delta) > {f_{13}(\Delta) \over g_3(\Gamma)} {f_{23}(\Delta) \over g_3(\Gamma)}$

We get $g_1(\Gamma_1) = {f_{13}(\Delta) \over g_3(\Gamma)}$ and $g_2(\Gamma_1) = {f_{23}(\Delta) \over g_3(\Gamma)}$.  As $\Gamma$ has only ${f_{13}(\Delta) \over g_3(\Gamma)}$ vertices of color 1 and ${f_{23}(\Delta) \over g_3(\Gamma)}$ vertices of color $2$ adjacent to any vertices of color $3$, any edge of color 12 contained in any facets must have its vertices among the first ${f_{13}(\Delta) \over g_3(\Gamma)}$ of color 1 and the first ${f_{23}(\Delta) \over g_3(\Gamma)}$ of color 2.  There are ${f_{13}(\Delta) \over g_3(\Gamma)}{f_{23}(\Delta) \over g_3(\Gamma)}$ such edges possible, each of which is contained in $g_3(\Gamma)$ facets, so we have $f_{123}(\Gamma) \leq {f_{13}(\Delta) \over g_3(\Gamma)}{f_{23}(\Delta) \over g_3(\Gamma)}g_3(\Gamma) = g_1(\Gamma_1)g_2(\Gamma_1)g_3(\Gamma_1) \leq f_{123}(\Gamma_1)$.  Since $\Gamma_1$ has at least as many edges of each color set as $\Gamma$ and $\Gamma \in \mathcal{D}(\Delta)$, we get $\Gamma_1 \in \mathcal{D}(\Delta)$.  \endproof

While the previous lemma assumes particular values of $p(\Gamma)$ and $q(\Gamma)$ for notational simplicity, it also applies to other values of $p(\Gamma)$ and $q(\Gamma)$ by relabeling colors.

The next lemma says that we know at least one $g_i(\Gamma)$ immediately.  If $i = r(\Gamma)$, then the previous lemma settles the problem.  If not, this at least puts considerable restrictions on what $g_{r(\Gamma)}(\Gamma)$ can be.

\begin{definition}
\textup{Define $\mathcal{E}(\Delta)$ by $\Gamma \in \mathcal{E}(\Delta)$ exactly if $\Gamma \in \mathcal{D}(\Delta)$ and $g_i(\Gamma) = b_i(\Delta)$ for some $i \in [3]$.}
\end{definition}

\begin{lemma} \label{oneb}
Let $\Delta$ be a 3-colored simplicial complex.  Either $\mathcal{D}(\Delta) = \mathcal{E}(\Delta)$ or else there are two complexes $\Gamma_1, \Gamma_2 \in \mathcal{E}(\Delta)$ with $r(\Gamma_1) \not = r(\Gamma_2)$.
\end{lemma}

\proof  In order for the first option of the lemma to not hold, there must be some $\Gamma \in \mathcal{D}(\Delta)$ with $\Gamma \not \in \mathcal{E}(\Delta)$.  We must either have $g_i(\Gamma) > b_i(\Delta)$ for at least two values of $i \in [3]$ or else $g_i(\Gamma) < b_i(\Delta)$ for at least two values of $i \in [3]$.  Suppose that it is the former.  Assume without loss of generality that $g_1(\Gamma) > b_1(\Delta)$ and $g_2(\Gamma) > b_2(\Delta)$.  Since these are all integers, $g_1(\Gamma) \geq b_1(\Delta) + 1$ and $g_2(\Gamma) \geq b_2(\Delta) + 1$.  We can compute
\begin{eqnarray*}
f_{12}(\Gamma) & \geq & g_1(\Gamma)g_2(\Gamma) \\ & \geq & (b_1(\Delta) + 1)(b_2(\Delta) + 1) \\ & = & \Bigg( \Bigg\lfloor \sqrt{{f_{12}(\Delta)f_{13}(\Delta) \over f_{23}(\Delta)}} \Bigg\rfloor + 1 \Bigg) \Bigg( \Bigg\lfloor \sqrt{{f_{12}(\Delta)f_{23}(\Delta) \over f_{13}(\Delta)}} \Bigg\rfloor + 1 \Bigg) \\ & > & \sqrt{{f_{12}(\Delta)f_{13}(\Delta) \over f_{23}(\Delta)}} \sqrt{{f_{12}(\Delta)f_{23}(\Delta) \over f_{13}(\Delta)}} \\ & = & f_{12}(\Delta) \\ & = & f_{12}(\Gamma).
\end{eqnarray*}
This is obviously impossible.

Otherwise, we must have $g_i(\Gamma) < b_i(\Delta)$ for at least two values of $i \in [3]$.  Assume without loss of generality that $g_1(\Gamma) < b_1(\Delta)$ and $g_2(\Gamma) < b_2(\Delta)$.  Since these are all integers, $g_1(\Gamma) \leq b_1(\Delta) - 1$ and $g_2(\Gamma) \leq b_2(\Delta) - 1$.  We compute
\begin{eqnarray*}
f_{12}(\Gamma) & \leq & (g_1(\Gamma)+1)(g_2(\Gamma)+1) \\ & \leq & b_1(\Delta)b_2(\Delta) \\ & = & \Bigg\lfloor \sqrt{{f_{12}(\Delta)f_{13}(\Delta) \over f_{23}(\Delta)}} \Bigg\rfloor \Bigg\lfloor \sqrt{{f_{12}(\Delta)f_{23}(\Delta) \over f_{13}(\Delta)}} \Bigg\rfloor \\ & \leq & \sqrt{{f_{12}(\Delta)f_{13}(\Delta) \over f_{23}(\Delta)}} \sqrt{{f_{12}(\Delta)f_{23}(\Delta) \over f_{13}(\Delta)}} \\ & = & f_{12}(\Delta) \\ & = & f_{12}(\Gamma).
\end{eqnarray*}
Because the opposite ends of the chain of inequalities are equal, equality must hold throughout.  For the first inequality to be an equality, we must have $\{p(\Gamma), q(\Gamma)\} = \{1, 2\}$.  We can assume without loss of generality that $p(\Gamma) = 1$ and $q(\Gamma) = 2$.  The second inequality means that $b_1(\Delta) = g_1(\Gamma) + 1$ and $b_2(\Delta) = g_2(\Gamma) + 1$.  The third gives that $\sqrt{{f_{12}(\Delta)f_{13}(\Delta) \over f_{23}(\Delta)}}$ and $\sqrt{{f_{12}(\Delta)f_{23}(\Delta) \over f_{13}(\Delta)}}$ are integers, so taking their floors does not change them.

If $g_3(\Gamma) = b_3(\Delta)$, then $\Gamma \in \mathcal{E}(\Delta)$, which contradicts the choice of $\Gamma$.  If $g_3(\Gamma) < b_3(\Delta)$, then we can apply the previous paragraph using 1 and 3 to get $\{p(\Gamma), q(\Gamma)\} = \{1, 3\}$, which contradicts $\{p(\Gamma), q(\Gamma)\} = \{1, 2\}$.

The only other possibility is if $g_3(\Gamma) > b_3(\Delta)$. In this case, we compute
\begin{eqnarray*}
f_{13}(\Delta) & = & \sqrt{{f_{12}(\Delta)f_{13}(\Delta) \over f_{23}(\Delta)}} \sqrt{{f_{13}(\Delta)f_{23}(\Delta) \over f_{12}(\Delta)}} \\ & < & \sqrt{{f_{12}(\Delta)f_{13}(\Delta) \over f_{23}(\Delta)}} \bigg( \bigg\lfloor \sqrt{{f_{13}(\Delta)f_{23}(\Delta) \over f_{12}(\Delta)}} \bigg\rfloor + 1 \bigg) \\ & = & b_1(\Delta)(b_3(\Delta)+1).
\end{eqnarray*}
Likewise, we can compute that $f_{23}(\Delta) < b_2(\Delta)(b_3(\Delta)+1)$.

Suppose that $g_3(\Gamma) \geq b_3(\Delta) + 2$.  We get that \begin{eqnarray*}
j_2(\Gamma) & = & f_{13}(\Delta) - g_1(\Gamma)g_3(\Gamma) \\
& < & b_1(\Delta)(b_3(\Delta)+1) - (b_1(\Delta)-1)(b_3(\Delta) + 2) \\ & = & b_1(\Delta)b_3(\Delta) + b_1(\Delta) - b_1(\Delta)b_3(\Delta) + b_3(\Delta) - 2b_1(\Delta) + 2 \\ & = & b_3(\Delta) - b_1(\Delta) + 2.
\end{eqnarray*}
This gives $g_3(\Gamma) - 1 \geq b_3(\Delta) + 1 \geq b_1(\Delta) + j_2(\Gamma)$.  By the same argument, $g_3(\Gamma) - 1 \geq b_2(\Delta) + j_1(\Gamma)$.  Thus, we can define $\Gamma_1$ by $p(\Gamma_1) = 1$, $q(\Gamma_1) = 2$, $g_3(\Gamma_1) = g_3(\Gamma) - 1$, $g_1(\Gamma_1) = g_1(\Gamma)$, and $g_2(\Gamma_1) = g_2(\Gamma)$ and have $\Gamma_1$ use all available edges.

We can compute that $\Gamma$ has $g_1(\Gamma)g_2(\Gamma)g_3(\Gamma)$ facets containing neither extra vertex, $g_2(\Gamma)j_2(\Gamma)$ facets containing $v_{b_1(\Delta)}^1$ but not $v_{b_2(\Delta)}^2$, $g_1(\Gamma)j_1(\Gamma)$ facets containing $v_{b_2(\Delta)}^2$ but not $v_{b_1(\Delta)}^1$, and $\min\{j_1(\Gamma), j_2(\Gamma)\}$ facets containing both extra vertices.  Similarly, we compute that $\Gamma_1$ has $g_1(\Gamma)g_2(\Gamma)(g_3(\Gamma)-1)$ facets containing neither extra vertex, $g_2(\Gamma)(j_2(\Gamma)+g_1(\Gamma))$ facets containing $v_{b_1(\Delta)}^1$ but not $v_{b_2(\Delta)}^2$, $g_1(\Gamma)(j_1(\Gamma)+g_2(\Gamma))$ facets containing $v_{b_2(\Delta)}^2$ but not $v_{b_1(\Delta)}^1$, and $\min\{j_1(\Gamma)+g_2(\Gamma), j_2(\Gamma)+g_1(\Gamma)\}$ facets containing both extra vertices.  Thus,
\begin{eqnarray*}
f_{123}(\Gamma_1) & = & g_1(\Gamma)g_2(\Gamma)(g_3(\Gamma)-1) + g_2(\Gamma)(j_2(\Gamma)+g_1(\Gamma)) \\ & & + g_1(\Gamma)(j_1(\Gamma)+g_2(\Gamma)) + \min\{j_1(\Gamma)+g_2(\Gamma), j_2(\Gamma)+g_1(\Gamma)\} \\ & = & g_1(\Gamma)g_2(\Gamma)g_3(\Gamma) + g_1(\Gamma)j_1(\Gamma) + g_2(\Gamma)j_2(\Gamma) \\ & & + g_1(\Gamma)g_2(\Gamma) + \min\{j_1(\Gamma)+g_2(\Gamma), j_2(\Gamma)+g_1(\Gamma)\} \\ & > & g_1(\Gamma)g_2(\Gamma)g_3(\Gamma) + g_1(\Gamma)j_1(\Gamma) + g_2(\Gamma)j_2(\Gamma) + \min\{j_1(\Gamma), j_2(\Gamma)\} \\ & = & f_{123}(\Gamma).
\end{eqnarray*}
Therefore, by Lemma~\ref{twovertmax}, $\Gamma \not \in \mathcal{B}(\Delta) \supset \mathcal{D}(\Delta)$, a contradiction.

Otherwise, $g_3(\Gamma) = b_3(\Delta) + 1$.  In this case, we define $\Gamma_1$ by $p(\Gamma_1) = 3$, $q(\Gamma_1) = 1$, and $g_i(\Gamma_1) = b_i(\Gamma)$ for all $i \in [3]$.  Since $f_{13}(\Delta) < b_1(\Delta)(b_3(\Delta)+1)$ and $f_{23}(\Delta) < b_1(\Delta)(b_3(\Delta)+1)$, the first extra vertex of $\Gamma_1$ uses up all remaining edges, and the second extra vertex is contained in no edges at all, so it doesn't matter if there is another vertex of color 1 available. We also find that $f_{13}(\Delta) = (b_1(\Delta) - 1)(b_3(\Delta) + 1) + j_2(\Gamma) < b_1(\Delta)(b_3(\Delta)+1)$, from which $j_2(\Gamma) < b_3(\Delta) + 1$, and so $j_2(\Gamma) \leq b_3(\Delta)$.  Similarly, $j_1(\Gamma) \leq b_3(\Delta)$.

We compute $j_2(\Gamma_1) = j_2(\Gamma) + b_1(\Delta) - 1 - b_3(\Delta)$ and $j_1(\Gamma_1) = j_1(\Gamma) + b_2(\Delta) - 1 - b_3(\Delta)$.  The first extra vertex of $\Gamma_1$ contains $(j_2(\Gamma) + b_1(\Delta) - 1 - b_3(\Delta))(j_1(\Gamma) + b_2(\Delta) - 1 - b_3(\Delta))$ facets.  We get that
\begin{eqnarray*}
f_{123}(\Gamma_1) & = & b_1(\Delta)b_2(\Delta)b_3(\Delta) \\ & & + (j_2(\Gamma) + b_1(\Delta) - 1 - b_3(\Delta))(j_1(\Gamma) + b_2(\Delta) - 1 - b_3(\Delta)) \\ & = & b_1(\Delta)b_2(\Delta)b_3(\Delta) + (b_1(\Delta) - 1)(b_2(\Delta) - 1) - (b_1(\Delta) - 1)b_3(\Delta) \\ & & - (b_2(\Delta) - 1)b_3(\Delta) + (b_1(\Delta) - 1)j_1(\Gamma) + (b_2(\Delta) - 1)j_2(\Gamma) \\ & & + (j_1(\Gamma) - b_3(\Delta))(j_2(\Gamma) - b_3(\Delta)) \\ & = & (b_1(\Delta) - 1)(b_2(\Delta) - 1)(b_3(\Delta) + 1) + (b_1(\Delta) - 1)j_1(\Gamma) \\ & & + (b_2(\Delta) - 1)j_2(\Gamma) + (b_3(\Delta) - j_1(\Gamma))(b_3(\Delta) - j_2(\Gamma)) + b_3(\Delta) \\ & \geq & (b_1(\Delta) - 1)(b_2(\Delta) - 1)(b_3(\Delta) + 1) + (b_1(\Delta) - 1)j_1(\Gamma) \\ & & + (b_2(\Delta) - 1)j_2(\Gamma) + \min\{j_1(\Gamma), j_2(\Gamma)\} \\ & = & f_{123}(\Gamma).
\end{eqnarray*}
The inequality comes because $j_1(\Gamma) \leq b_3(\Delta)$ and $j_2(\Gamma) \leq b_3(\Delta)$.

Since $\Gamma \in \mathcal{B}(\Delta)$, we get $\Gamma_1 \in \mathcal{B}(\Delta)$.  As we have already seen that $\Gamma_1$ uses all available edges, $\Gamma_1 \in \mathcal{C}(\Delta)$.  Since $\Gamma \in \mathcal{D}(\Delta)$, we know that $\mathcal{D}(\Delta) \not = \emptyset$.  The alternative is $\mathcal{D}(\Delta) = \mathcal{C}(\Delta)$, and so $\Gamma_1 \in \mathcal{D}(\Delta)$.  Since $b_3(\Delta) = g_3(\Gamma_1)$, we get $\Gamma_1 \in \mathcal{E}(\Delta)$.

There was nothing special about choosing $q(\Gamma_1) = 1$, as the second extra vertex was not used at all.  If we defined $\Gamma_2$ in exactly the same way as $\Gamma_1$ except that $q(\Gamma_2) = 2$, then $\Gamma_2 \in \mathcal{E}(\Delta)$ by the same argument as $\Gamma_1$.  This completes the proof because $r(\Gamma_2) = 1 \not = 2 = r(\Gamma_1)$, which is the second option of the lemma.  \endproof

The basic approach to find a complex in $\mathcal{E}(\Delta)$ is to drop the requirement that a complex maximizes the number of facets and check all of the complexes that satisfy the rest of the conditions of $\mathcal{E}(\Delta)$ and could plausibly maximize the number of facets.  Knowing that $\mathcal{E}(\Delta) \not = \emptyset$ tells us that at least one such complex must maximize the number of facets, and hence be in $\mathcal{E}(\Delta)$.  Whichever complex has the most facets from among the ones we check must be such a complex.  The next definition makes this more explicit.

\begin{definition}
\textup{Let $\Delta$ be a 3-colored simplicial complex.  Define $\mathcal{F}(\Delta)$ by $\Gamma \in \mathcal{F}(\Delta)$ exactly if}
\begin{enumerate}
\item \textup{$\Gamma \in \mathcal{A}(\Delta)$,}
\item \textup{$\Gamma$ has exactly as many edges as $\Delta$,}
\item \textup{$\mathcal{D}(\Delta) \not = \emptyset$, and}
\item \textup{$g_i(\Gamma) = b_i(\Delta)$ for some $i \in [3]$.}
\end{enumerate}
\end{definition}

For any $\Gamma \in \mathcal{F}(\Delta)$, we must have $g_i(\Gamma) = b_i(\Delta)$ for some $i \in [3]$, and there are only three ways to pick $i$.  It is clear from the definition that there are only three ways to pick $r(\Gamma)$.  The basic plan to compute $m(\Delta)$ is to try all possible combinations of a choice of $i$ and of $r(\Gamma)$, for nine cases in all.  The three cases where $g_{r(\Gamma)}(\Gamma) = b_{r(\Gamma)}(\Delta)$ are quickly handled by Lemma~\ref{knowr}.

It is sometimes convenient to assume without loss of generality that $f_{12}(\Delta) \leq f_{13}(\Delta) \leq f_{23}(\Delta)$.  We can do this because if it is not true for a given complex $\Delta$, we can fix that by relabeling the colors.  In some lemmas, we relax this assumption a bit for the sake of generality.

The next lemma says that we can handle all of the cases where $r(\Gamma) = 3$ by checking only the case where $g_3(\Gamma) = f_3(\Delta)$.

\begin{lemma}  \label{nor3}
Let $\Delta$ be a 3-colored simplicial complex such that $f_{12}(\Delta) \leq f_{23}(\Delta)$ and $f_{13}(\Delta) \leq f_{23}(\Delta)$.  If $\Gamma_0 \in \mathcal{F}(\Delta)$, $r(\Gamma_0) = 3$, and $f_3(\Delta) > g_3(\Gamma_0)$, then at least one of the following holds:
\begin{enumerate}
\item there is some $\Gamma \in \mathcal{F}(\Delta)$ with $r(\Gamma) \not = 3$ and $f_{123}(\Gamma) \geq f_{123}(\Gamma_0)$;
\item there is some $\Gamma \in \mathcal{A}(\Delta)$ with $f_{123}(\Gamma) > f_{123}(\Gamma_0)$; or
\item there is some $\Gamma \in \mathcal{F}(\Delta)$ with $r(\Gamma) = 3$, $f_{123}(\Gamma) \geq f_{123}(\Gamma_0)$, and $g_3(\Gamma) = f_3(\Delta)$.
\end{enumerate}
\end{lemma}

\proof  We have that $p(\Gamma_0)$ and $q(\Gamma_0)$ are 1 and 2 in some order.  Assume without loss of generality that $p(\Gamma_0) = 1$.

We break the proof into several cases.  Each time that we construct a complex $\Gamma$, we need to check that it is well-defined, in $\mathcal{A}(\Delta)$, and if we might have $f_{123}(\Gamma) = f_{123}(\Gamma_0)$, also that $\Gamma \in \mathcal{F}(\Delta)$.  To show that $\Gamma$ is well-defined, it suffices to show that there are enough edges and vertices available to construct the complex.  That $\Gamma \in \mathcal{A}(\Delta)$ follows from the definition.  To check that $\Gamma \in \mathcal{F}(\Delta)$, the first and fourth conditions are true by construction and the third part holds because it is a property of $\Delta$ and must hold to get $\Gamma_0 \in \mathcal{F}(\Delta)$.  It thus suffices to check the second condition.

In constructions where $g_i(\Gamma) \leq g_i(\Gamma_0)$ for all $i \in [3]$, there are enough vertices because $\Gamma$ uses at most as many vertices of each color as $\Gamma_0$, except that $\Gamma$ could use one additional vertex of color 3, which is available because $f_3(\Delta) > g_3(\Gamma_0)$.  Without new non-extra vertices, no additional edges are forced to be in the complex by the non-extra vertices, so $\Gamma$ is well-defined because $\Gamma_0$ is.

We can assume that $f_1(\Delta) > g_1(\Gamma_0)$ and $f_2(\Delta) > g_2(\Gamma_0)$, as if not, then one of the extra vertices is completely missing, so we can drop it, set $p(\Gamma)$ as the color of the remaining extra vertex, make $q(\Gamma) = 3$, and get $\Gamma_0 \subset \Gamma$. We can also assume that $g_3(\Gamma_0) < f_3(\Delta)$, as otherwise, we can take $\Gamma = \Gamma_0$ and meet the third option of the lemma.

We note that $j_1(\Gamma_0) \leq g_3(\Gamma_0)$ and $j_2(\Gamma_0) \leq g_3(\Gamma_0)$, as this is necessary for $\Gamma_0$ to use all edges of colors 23 and 13, respectively, as it has no extra vertex of color 3.

\textbf{Case I:}  $j_1(\Gamma_0) \leq j_2(\Gamma_0)$

\textbf{Case I A:}  $j_3(\Gamma_0) \leq g_2(\Gamma_0)$

This case means that the first extra vertex of $\Gamma_0$ uses all available edges of color 12.  In particular, this means that the second extra vertex adds no additional facets.  We can compute that the first extra vertex of $\Gamma_0$ adds $j_3(\Gamma_0)j_2(\Gamma_0)$ facets, so
$$f_{123}(\Gamma_0) = g_1(\Gamma_0)g_2(\Gamma_0)g_3(\Gamma_0) + j_3(\Gamma_0)j_2(\Gamma_0).$$

\textbf{Case I A 1:}  $j_1(\Gamma_0) \leq g_2(\Gamma_0)$

Define $\Gamma$ by $p(\Gamma) = 1$, $q(\Gamma) = 3$, and $g_i(\Gamma) = g_i(\Gamma_0)$ for all $i \in [3]$.  Because the second extra vertex of $\Gamma_0$ adds no additional facets, every facet of $\Gamma_0$ is also in $\Gamma$, so $f_{123}(\Gamma) \geq f_{123}(\Gamma_0)$.  The first extra vertex of $\Gamma$ uses all edges of colors 12 and 13 because it also does so in $\Gamma_0$.  Finally, because $j_1(\Gamma_0) \leq g_2(\Gamma_0)$, the second extra vertex of $\Gamma$ uses all remaining edges of color 23.  Hence, $\Gamma \in \mathcal{F}(\Delta)$, and we have the first option of the lemma.

\textbf{Case I A 2:}  $j_1(\Gamma_0) > g_2(\Gamma_0)$

\textbf{Case I A 2 a:}  $j_2(\Gamma_0) < g_1(\Gamma_0)$

Define $\Gamma$ by $q(\Gamma) = 2$, $p(\Gamma) = 3$, and $g_i(\Gamma) = g_i(\Gamma_0)$ for all $i \in [3]$.  We can chain the inequalities of this case to get
$$g_1(\Gamma_0) > j_2(\Gamma_0) \geq j_1(\Gamma_0) > g_2(\Gamma_0) \geq j_3(\Gamma_0).$$
In particular, $g_1(\Gamma_0) > j_3(\Gamma_0)$, so the second extra vertex of $\Gamma$ uses all available edges of color 12.  That $j_2(\Gamma_0) < g_1(\Gamma_0)$ means that the first extra vertex of $\Gamma$ uses all available edges of color 13.  The second extra vertex of $\Gamma$ uses any leftover edges of color 23 because it does in $\Gamma_0$.  Thus, $\Gamma \in \mathcal{F}(\Delta)$.

We can compute that the first extra vertex of $\Gamma$ adds $j_2(\Gamma_0)g_2(\Gamma_0)$ facets, and the second one adds $j_3(\Gamma_0)(j_1(\Gamma_0) - g_2(\Gamma_0))$ facets.  This allows us to compute
\begin{eqnarray*}
f_{123}(\Gamma) & = & g_1(\Gamma_0)g_2(\Gamma_0)g_3(\Gamma_0) + j_2(\Gamma_0)g_2(\Gamma_0) + j_3(\Gamma_0)(j_1(\Gamma_0) - g_2(\Gamma_0)) \\ & \geq & g_1(\Gamma_0)g_2(\Gamma_0)g_3(\Gamma_0) + j_2(\Gamma_0)g_2(\Gamma_0) \\ & \geq & g_1(\Gamma_0)g_2(\Gamma_0)g_3(\Gamma_0) + j_2(\Gamma_0)j_3(\Gamma_0) \\ & = & f_{123}(\Gamma_0),
\end{eqnarray*}
which yields the first option in the lemma.

\textbf{Case I A 2 b:}  $j_2(\Gamma_0) \geq g_1(\Gamma_0)$

Let $w = \min\big\{ \big\lfloor {j_2(\Gamma_0) \over g_1(\Gamma_0)} \big\rfloor, \big\lfloor {j_1(\Gamma_0) \over g_2(\Gamma_0)} \big\rfloor, f_3(\Delta) - g_3(\Gamma_0) \big\}$.  Note that $w \geq 1$ because $j_2(\Gamma_0) \geq g_1(\Gamma_0)$, $j_1(\Gamma_0) > g_2(\Gamma_0)$, and $f_3(\Delta) > g_3(\Gamma_0)$.  Define $\Gamma_1$ by $p(\Gamma_1) = 1$, $q(\Gamma_1) = 2$, $g_1(\Gamma_1) = g_1(\Gamma_0)$, $g_2(\Gamma_1) = g_2(\Gamma_0)$, and $g_3(\Gamma_1) = g_3(\Gamma_0) + w$.  There are enough edges to do this because the first two terms assert that $w$ is small enough not to use more edges than allowed of color sets 13 or 23, respectively.  The third term of $w$ that there are enough vertices of color 3 available.  Thus, $\Gamma_1$ is well-defined.  Because the extra vertices in $\Gamma_1$ are able to use at least as many edges of each color set as those of $\Gamma_0$ and have at most as many such edges available to use, $\Gamma_1 \in \mathcal{F}(\Delta)$.

The first extra vertex of $\Gamma_1$ uses all edges of color 12, so the second extra vertex adds no additional facets.  Meanwhile, the first extra vertex of $\Gamma_1$ adds $j_2(\Gamma_1)j_3(\Gamma_1) = j_3(\Gamma_0)(j_2(\Gamma_0) - wg_1(\Gamma_0))$ facets.  This yields
\begin{eqnarray*}
f_{123}(\Gamma_1) & = & g_1(\Gamma_0)g_2(\Gamma_0)(g_3(\Gamma_0)+w) + j_3(\Gamma_0)(j_2(\Gamma_0) - wg_1(\Gamma_0)) \\ & = & g_1(\Gamma_0)g_2(\Gamma_0)g_3(\Gamma_0) + j_3(\Gamma_0)j_2(\Gamma_0) + wg_1(\Gamma_0)(g_2(\Gamma_0) - j_3(\Gamma_0)) \\ & \geq & g_1(\Gamma_0)g_2(\Gamma_0)g_3(\Gamma_0) + j_3(\Gamma_0)j_2(\Gamma_0) \\ & = & f_{123}(\Gamma_0).
\end{eqnarray*}

\textbf{Case I A 2 b i:}  $j_3(\Gamma_0) < g_2(\Gamma_0)$

This ensures that the inequality above is strict, so we can take $\Gamma = \Gamma_1$ and have the second option of the lemma.

\textbf{Case I A 2 b ii:}  $w = f_3(\Delta) - g_3(\Gamma_0)$

This ensures that $g_3(\Gamma_1) = f_3(\Delta)$, so we can take $\Gamma = \Gamma_1$ and have the third option of the lemma.

\textbf{Case I A 2 b iii:}  $j_3(\Gamma_0) = g_2(\Gamma_0)$ and $w < f_3(\Delta) - g_3(\Gamma_0)$

\textbf{Case I A 2 b iii (a):}  $j_1(\Gamma_1) \leq g_2(\Gamma_1)$

Define $\Gamma$ by $g_i(\Gamma) = g_i(\Gamma_1)$ for all $i \in [3]$, $p(\Gamma) = 1$ and $q(\Gamma) = 3$.  We have that $\Gamma$ is well-defined because $\Gamma_1$ is.  The first extra vertex of $\Gamma$ uses all available edges of colors 12 and 13 because the first extra vertex of $\Gamma_1$ does also.  Because $j_1(\Gamma_1) \leq g_2(\Gamma_1)$, the second extra vertex of $\Gamma$ uses all remaining edges of color 23.  Hence, $\Gamma \in \mathcal{F}(\Delta)$.  Furthermore, since the only differing edges between $\Gamma$ and $\Gamma_1$ are the ones in the second extra vertex, none of which are in any facets, we have $f_{123}(\Gamma) = f_{123}(\Gamma_1) \geq f_{123}(\Gamma_0)$, so we satisfy the first option of the lemma.

\textbf{Case I A 2 b iii (b):}  $j_1(\Gamma_1) > g_2(\Gamma_1)$

This means that $w \not = \big\lfloor {j_1(\Gamma_0) \over g_2(\Gamma_0)} \big\rfloor$.  Since $w < f_3(\Delta) - g_3(\Gamma_0)$, we must have $w = \big\lfloor {j_2(\Gamma_0) \over g_1(\Gamma_0)} \big\rfloor$.  Hence, $j_2(\Gamma_1) = j_2(\Gamma_0) - wg_1(\Gamma_0) < g_1(\Gamma_0)$.  That $w \not = \big\lfloor {j_1(\Gamma_0) \over g_2(\Gamma_0)} \big\rfloor$ means that $w+1 \leq \big\lfloor {j_1(\Gamma_0) \over g_2(\Gamma_0)} \big\rfloor \leq {j_1(\Gamma_0) \over g_2(\Gamma_0)}$.  This yields ${j_2(\Gamma_0) \over g_1(\Gamma_0)} < w+1 \leq {j_1(\Gamma_0) \over g_2(\Gamma_0)}$, so $j_2(\Gamma_0)g_2(\Gamma_0) < j_1(\Gamma_0)g_1(\Gamma_0)$.  Since $j_2(\Gamma_0) \geq j_1(\Gamma_0)$, we must have $g_2(\Gamma_0) < g_1(\Gamma_0)$.

\textbf{Case I A 2 b iii (b) (i):}  $j_1(\Gamma_1) > j_2(\Gamma_1)$

Define $\Gamma$ by $g_i(\Gamma) = g_i(\Gamma_1)$ for all $i \in [3]$, $p(\Gamma) = 2$ and $q(\Gamma) = 1$.  Then $\Gamma$ is well-defined because $\Gamma_1$ is.  We can compute that $\Gamma$ has $g_1(\Gamma)g_2(\Gamma)g_3(\Gamma)$ facets before adding extra vertices.  The first extra vertex of $\Gamma$ uses up all remaining edges of color 12 because $j_3(\Gamma) = j_3(\Gamma_0) = g_2(\Gamma_0) < g_1(\Gamma_0)$.  Furthermore, $\Gamma$ uses all remaining edges of color 23 because $j_1(\Gamma) < j_1(\Gamma_0) \leq g_3(\Gamma_0) < g_3(\Gamma)$.  Thus, we can compute
\begin{eqnarray*}
f_{123}(\Gamma) & = & g_1(\Gamma)g_2(\Gamma)g_3(\Gamma) + j_3(\Gamma)j_1(\Gamma) \\ & = & g_1(\Gamma)g_2(\Gamma)g_3(\Gamma) + j_3(\Gamma_1)j_1(\Gamma_1) \\ & > & g_1(\Gamma)g_2(\Gamma)g_3(\Gamma) + j_3(\Gamma_1)j_2(\Gamma_1) \\ & = & f_{123}(\Gamma_1) \\ & \geq & f_{123}(\Gamma_0),
\end{eqnarray*}
which yields the second option of the lemma because $\Gamma \in \mathcal{A}(\Delta)$ by construction.

\textbf{Case I A 2 b iii (b) (ii):}  $j_1(\Gamma_1) \leq j_2(\Gamma_1)$

Define $\Gamma_2$ by $g_1(\Gamma_2) = g_1(\Gamma_0)$, $g_2(\Gamma_2) = g_2(\Gamma_0)$, $g_3(\Gamma_2) = \big\lceil {f_{13}(\Delta) \over g_1(\Gamma_0) + 1} \big\rceil$, $p(\Gamma_2) = 1$, and $q(\Gamma_2) = 2$.  Let $y = g_3(\Gamma_0) - g_3(\Gamma_2)$.  Since
$$f_{13}(\Delta) = g_1(\Gamma_0)g_3(\Gamma_0) + j_2(\Gamma_0) \leq g_1(\Gamma_0)g_3(\Gamma_0) + g_3(\Gamma_0) = (g_1(\Gamma_0)+1)g_3(\Gamma_0),$$
we have $g_3(\Gamma_2) = \big\lceil {f_{13}(\Delta) \over g_1(\Gamma_0) + 1} \big\rceil \leq \lceil g_3(\Gamma_0) \rceil = g_3(\Gamma_0)$, and so $y \geq 0$.  Hence, $\Gamma_2$ is well-defined because $\Gamma_0$ is.  It uses all edges of color 12 because $\Gamma_0$ does.  $\Gamma_2$ uses all edges of color 13 because $g_1(\Gamma_0) + 1$ vertices of color 1 and $g_3(\Gamma_2)$ of color 3 can use up to $(g_1(\Gamma_0) + 1)g_3(\Gamma_2) = (g_1(\Gamma_0) + 1)\big\lceil {f_{13}(\Delta) \over g_1(\Gamma_0) + 1} \big\rceil \geq f_{13}(\Delta)$ edges of color 13.  Finally, $\Gamma_2$ uses all edges of color 23 because
$$j_1(\Gamma_2) = j_1(\Gamma_0) + yg_2(\Gamma_0) \leq j_2(\Gamma_0) + yg_1(\Gamma_0) = j_2(\Gamma_2) \leq g_3(\Gamma_2),$$
as $\Gamma_2$ also uses all edges of color 13. Therefore, $\Gamma_2 \in \mathcal{F}(\Delta)$.

Before adding any extra vertices, $\Gamma_2$ has $g_1(\Gamma_0)g_2(\Gamma_0)(g_3(\Gamma_0)-y)$ facets.  The first extra vertex adds $g_2(\Gamma_0)(j_2(\Gamma_0) + yg_1(\Gamma_0))$ facets.  This uses all edges of color 12, so the second extra vertex adds no additional facets.  Thus, we compute
\begin{eqnarray*}
f_{123}(\Gamma_2) & = & g_1(\Gamma_0)g_2(\Gamma_0)(g_3(\Gamma_0)-y) + g_2(\Gamma_0)(j_2(\Gamma_0) + yg_1(\Gamma_0)) \\ & = & g_1(\Gamma_0)g_2(\Gamma_0)g_3(\Gamma_0) + g_2(\Gamma_0)j_2(\Gamma_0) \\ & = & g_1(\Gamma_0)g_2(\Gamma_0)g_3(\Gamma_0) + j_3(\Gamma_0)j_2(\Gamma_0) \\ & = & f_{123}(\Gamma_0).
\end{eqnarray*}

Define $\Gamma$ by $g_1(\Gamma) = g_1(\Gamma_2) + 1$, $g_2(\Gamma) = g_2(\Gamma_2)$, $g_3(\Gamma) = g_3(\Gamma_2) - 1$, $p(\Gamma) = 3$, and $q(\Gamma) = 2$.  There are enough edges of color 23 to do this because $\Gamma_2$ is well-defined.  There are enough edges of color 12 because
$$f_{12}(\Delta) = g_1(\Gamma_0)g_2(\Gamma_0) + j_3(\Gamma_0) = g_1(\Gamma_0)g_2(\Gamma_0) + g_2(\Gamma_0) = g_1(\Gamma)g_2(\Gamma).$$
There are enough edges of color 13 because
\begin{eqnarray*}
g_1(\Gamma)g_3(\Gamma) & = & (g_1(\Gamma_0)+1) \bigg( \bigg\lceil {f_{13}(\Delta) \over g_1(\Gamma_0) + 1} \bigg\rceil - 1 \bigg) \\ & < & (g_1(\Gamma_0)+1) \bigg({f_{13}(\Delta) \over g_1(\Gamma_0) + 1} \bigg) \\ & = & f_{13}(\Delta).
\end{eqnarray*}
Therefore, $\Gamma$ is well-defined.

It is easy to check that $\Gamma$ and $\Gamma_2$ use exactly the same edges of colors 12 and 23.  They have exactly the same vertices of all colors, so $\Gamma$ uses all available edges because $\Gamma_2$ does.  Hence, $\Gamma \in \mathcal{F}(\Delta)$.  The only edges that can differ are that $\Gamma$ could include some edges containing $v_{g_1(\Gamma_2)+1}^1$ but not $v_{g_3(\Gamma_2)}^3$, while $\Gamma_2$ contains exactly the same number of edges containing $v_{g_3(\Gamma_2)}^3$ but not $v_{g_1(\Gamma_2)+1}^1$.  Both of these vertices are adjacent to exactly $g_2(\Gamma_0)$ vertices of color 2, so every edge that differs between $\Gamma$ and $\Gamma_2$ is contained in exactly $g_2(\Gamma_0)$ facets.  Therefore, $f_{123}(\Gamma) = f_{123}(\Gamma_2) \geq f_{123}(\Gamma_0)$, which gives us the first option of the lemma.

\textbf{Case I B:}  $j_3(\Gamma_0) > g_2(\Gamma_0)$

The first extra vertex of $\Gamma_0$ adds $g_2(\Gamma_0)j_2(\Gamma_0)$ facets, while the second extra vertex adds $j_1(\Gamma_0)(j_3(\Gamma_0) - g_2(\Gamma_0))$.  Thus,
$$f_{123}(\Gamma_0) = g_1(\Gamma_0)g_2(\Gamma_0)g_3(\Gamma_0) + g_2(\Gamma_0)j_2(\Gamma_0) + j_1(\Gamma_0)(j_3(\Gamma_0) - g_2(\Gamma_0)).$$

\textbf{Case I B 1:}  $j_2(\Gamma_0) + g_1(\Gamma_0) < g_3(\Gamma_0)$

Define $\Gamma$ by $p(\Gamma) = 1$, $q(\Gamma) = 2$, $g_1(\Gamma) = g_1(\Gamma_0)$, $g_2(\Gamma) = g_2(\Gamma_0)$, and $g_3(\Gamma) = g_3(\Gamma_0) - 1$.  We have already seen that this is well-defined, and it is immediate from the definition that $\Gamma \in \mathcal{A}(\Delta)$.  The first extra vertex of $\Gamma$ adds $g_2(\Gamma_0)(j_2(\Gamma_0) + g_1(\Gamma_0))$ facets.

\textbf{Case I B 1 a:}  $j_1(\Gamma_0) + g_2(\Gamma_0) < g_3(\Gamma_0)$

This allows the second extra vertex of $\Gamma$ to use all remaining edges of color 23, so it adds $(j_1(\Gamma_0) + g_2(\Gamma_0))(j_3(\Gamma_0) - g_2(\Gamma_0))$ facets.  Thus, we have
\begin{eqnarray*}
f_{123}(\Gamma) & = & g_1(\Gamma_0)g_2(\Gamma_0)(g_3(\Gamma_0)-1) + g_2(\Gamma_0)(j_2(\Gamma_0) + g_1(\Gamma_0)) \\ & & \qquad + (j_1(\Gamma_0) + g_2(\Gamma_0))(j_3(\Gamma_0) - g_2(\Gamma_0)) \\ & > & g_1(\Gamma_0)g_2(\Gamma_0)g_3(\Gamma_0) + g_2(\Gamma_0)j_2(\Gamma_0) + j_1(\Gamma_0)(j_3(\Gamma_0) - g_2(\Gamma_0)) \\ & = & f_{123}(\Gamma_0),
\end{eqnarray*}
giving us the second option of the lemma.

\textbf{Case I B 1 b:}  $j_1(\Gamma_0) + g_2(\Gamma_0) \geq g_3(\Gamma_0)$

The second extra vertex of $\Gamma$ can only use $g_3(\Gamma) = g_3(\Gamma_0) - 1$ edges of color 23, so it only adds $(g_3(\Gamma_0) - 1)(j_3(\Gamma_0) - g_2(\Gamma_0))$ facets.  Still, we have
$$g_3(\Gamma_0) - 1 \geq j_2(\Gamma_0) + g_1(\Gamma_0) > j_2(\Gamma_0) \geq j_1(\Gamma_0),$$
from which we can compute
\begin{eqnarray*}
f_{123}(\Gamma) & = & g_1(\Gamma_0)g_2(\Gamma_0)(g_3(\Gamma_0)-1) + g_2(\Gamma_0)(j_2(\Gamma_0) + g_1(\Gamma_0)) \\ & & \qquad + (g_3(\Gamma_0) - 1)(j_3(\Gamma_0) - g_2(\Gamma_0)) \\ & > & g_1(\Gamma_0)g_2(\Gamma_0)g_3(\Gamma_0) + g_2(\Gamma_0)j_2(\Gamma_0) + j_1(\Gamma_0)(j_3(\Gamma_0) - g_2(\Gamma_0)) \\ & = & f_{123}(\Gamma_0),
\end{eqnarray*}
which is the second option in the lemma.

\textbf{Case I B 2:}  $j_2(\Gamma_0) + g_1(\Gamma_0) \geq g_3(\Gamma_0)$

\textbf{Case I B 2 a:}  $j_2(\Gamma_0) = g_3(\Gamma_0)$

Define $\Gamma$ by $g_1(\Gamma) = g_1(\Gamma_0)+1$, $g_2(\Gamma) = g_2(\Gamma_0)$, $g_3(\Gamma) = g_3(\Gamma_0)$, $p(\Gamma) = 2$, and $q(\Gamma) = 3$.  The first extra vertex of $\Gamma_0$ is adjacent to all vertices of color 3 because $j_2(\Gamma_0) = g_3(\Gamma_0)$ and the first $g_2(\Gamma_0)$ vertices of color 2 because $j_3(\Gamma_0) > g_2(\Gamma_0)$.  All that $\Gamma$ does is to make this vertex no longer an extra vertex, and then tack on an extra isolated vertex of color 3 at the end.  Hence, $\Gamma_0 \subset \Gamma$, so $\Gamma$ must use at least as many edges and facets as $\Gamma_0$.  This means that $\Gamma \in \mathcal{F}(\Delta)$, and satisfies the first option of the lemma.

\textbf{Case I B 2 b:}  $j_2(\Gamma_0) < g_3(\Gamma_0)$

Define $\Gamma$ by $g_1(\Gamma) = g_1(\Gamma_0) + 1$, $g_2(\Gamma) = g_2(\Gamma_0)$, $g_3(\Gamma) = g_3(\Gamma_0) - 1$, $p(\Gamma) = 3$, and $q(\Gamma) = 2$.  We know that there are enough edges of color 23 for $\Gamma$ because there are enough for $\Gamma_0$, which needs more.  There are enough edges of color 12 because $\Gamma$ only needs $g_2(\Gamma_0)$ more than $\Gamma_0$ and $j_3(\Gamma_0) > g_2(\Gamma_0)$.  There are enough edges of color 13 because
\begin{eqnarray*}
f_{13}(\Delta) & = & g_1(\Gamma_0)g_3(\Gamma_0) + j_2(\Gamma_0) \\ & \geq & g_1(\Gamma_0)g_3(\Gamma_0) + g_3(\Gamma_0) - g_1(\Gamma_0) \\ & > & (g_1(\Gamma_0) + 1)(g_3(\Gamma_0) - 1) \\ & = & g_1(\Gamma)g_3(\Gamma).
\end{eqnarray*}
The vertices of $\Gamma$ are precisely the vertices of $\Gamma_0$ (with the last vertex of color 1 coming because $g_1(\Gamma_0) < f_1(\Delta)$), so $\Gamma$ is well-defined.  Furthermore, because $\Gamma$ and $\Gamma_0$ have the same vertices, they can each use the same number of edges of each color set.  Therefore, $\Gamma \in \mathcal{F}(\Delta)$ because $\Gamma_0$ is also.

One can easily check that $\Gamma$ and $\Gamma_0$ have exactly the same edges of colors 12 and 23.  The edges of color 13 can only differ in that $\Gamma$ may have some edges that $\Gamma_0$ lacks containing $v_{g_1(\Gamma_0)+1}^1$ but not $v_{g_3(\Gamma_0)}^3$, while $\Gamma_0$ could have some edges that are missing from $\Gamma$ and contain $v_{g_3(\Gamma_0)}^3$ but not $v_{g_1(\Gamma_0)+1}^1$.  Any edge that only $\Gamma$ has is contained in at least $g_2(\Gamma_0)$ facets, as both of its vertices are present before adding any extra vertices.  Any edge that only $\Gamma_0$ has is contained in at most $g_2(\Gamma_0)$ facets, as $v_{g_3(\Gamma_0)}^3$ is not adjacent to $v_{g_2(\Gamma_0) + 1}^2$ because $j_2(\Gamma_0) < g_3(\Gamma)$.  Therefore, each differing edge of $\Gamma$ has at least as many facets as each one of $\Gamma_0$, so $f_{123}(\Gamma) \geq f_{123}(\Gamma_0)$, which gives us the first option in the lemma.

\textbf{Case II:}  $j_1(\Gamma_0) > j_2(\Gamma_0)$

Define $\Gamma_1$ by $p(\Gamma_1) = 2$, $q(\Gamma_1) = 1$, and $g_i(\Gamma_1) = g_i(\Gamma_0)$ for all $i \in [3]$.  We have that $f_{123}(\Gamma_1) \geq f_{123}(\Gamma_0)$ by Lemma~\ref{choosepq}.  Because $\Gamma_0$ and $\Gamma_1$ have the same vertices, including the same extra vertices, $\Gamma_1 \in \mathcal{F}(\Delta)$.  Applying Case I to $\Gamma_1$ gives that the lemma holds for $\Gamma_0$.  \endproof

The next lemma settles the cases where $g_3(\Gamma) = b_3(\Delta)$ and $r(\Gamma) \not = 3$, as well as the case $g_2(\Gamma) = b_2(\Delta)$ and $r(\Gamma) = 1$.

\begin{lemma}  \label{bibr}
Let $\Delta$ be a 3-colored simplicial complex and let $\Gamma \in \mathcal{E}(\Delta)$.  If $g_i(\Gamma) = b_i(\Delta) \geq b_{r(\Gamma)}(\Delta)$ and $r(\Gamma) \not = i$, then $\Big\lceil {f_{ir(\Gamma)}(\Delta) \over b_i(\Delta) + 1} \Big\rceil \leq \Big\lfloor {f_{ir(\Gamma)}(\Delta) \over b_i(\Delta)} \Big\rfloor$.  Furthermore, either $g_{r(\Gamma)}(\Gamma) = \Big\lceil {f_{ir(\Gamma)}(\Delta) \over b_i(\Delta) + 1} \Big\rceil$ or else $g_{r(\Gamma)}(\Gamma) = \Big\lfloor {f_{ir(\Gamma)}(\Delta) \over b_i(\Delta)} \Big\rfloor$.
\end{lemma}

\proof  Because $\Gamma \in \mathcal{E}(\Delta) \subset \mathcal{D}(\Delta)$, it must use all available edges.  Let $k = 6 - r(\Gamma) - i$.  Since there is no extra vertex of color $r(\Gamma)$, we have $j_k(\Gamma) \leq g_{r(\Gamma)}(\Gamma)$.  We trivially must have $j_k(\Gamma) \geq 0$, so we have $g_{r(\Gamma)}(\Gamma)b_i(\Delta) \leq g_{r(\Gamma)}(\Gamma)b_i(\Delta) + j_k(\Gamma) \leq g_{r(\Gamma)}(\Gamma)b_i(\Delta) + g_{r(\Gamma)}(\Gamma)$.  The middle term is $f_{ir(\Gamma)}(\Delta)$, so we have $g_{r(\Gamma)}(\Gamma)b_i(\Delta) \leq f_{ir(\Gamma)}(\Delta) \leq g_{r(\Gamma)}(\Gamma)(b_i(\Delta) + 1)$.  Dividing the two inequalities as appropriate, we get $g_{r(\Gamma)}(\Gamma) \leq {f_{ir(\Gamma)}(\Delta) \over b_i(\Delta)}$ and $g_{r(\Gamma)}(\Gamma) \geq {f_{ir(\Gamma)}(\Delta) \over b_i(\Delta) + 1}$, respectively.  Chain these together to get ${f_{ir(\Gamma)}(\Delta) \over b_i(\Delta) + 1} \leq g_{r(\Gamma)}(\Gamma) \leq {f_{ir(\Gamma)}(\Delta) \over b_i(\Delta)}$.  Since $g_{r(\Gamma)}(\Gamma)$ is an integer, we have $\Big\lceil {f_{ir(\Gamma)}(\Delta) \over b_i(\Delta) + 1} \Big\rceil \leq g_{r(\Gamma)}(\Gamma) \leq \Big\lfloor {f_{ir(\Gamma)}(\Delta) \over b_i(\Delta)} \Big\rfloor$, which gives one inequality of the lemma.

Next, we compute
\begin{eqnarray*}
{f_{ir(\Gamma)}(\Delta) \over b_i(\Delta)} - {f_{ir(\Gamma)}(\Delta) \over b_i(\Delta) + 1} & = & f_{ir(\Gamma)}(\Delta) \Big( {1 \over b_i(\Delta)} - {1 \over b_i(\Delta) + 1} \Big) \\ & = & {f_{ir(\Gamma)}(\Delta) \over b_i(\Delta)(b_i(\Delta) + 1)} \\ & < & {(b_i(\Delta) + 1)(b_{r(\Gamma)}(\Delta) + 1) \over b_i(\Delta)(b_i(\Delta) + 1)} \\ & = & {b_{r(\Gamma)}(\Delta) + 1 \over b_i(\Delta)} \\ & \leq & {b_i(\Delta) + 1 \over b_i(\Delta)} \\ & = & 1 + {1 \over b_i(\Delta)} \\ & \leq & 2.
\end{eqnarray*}
Thus, $g_{r(\Gamma)}(\Gamma)$ is an integer contained in an interval of length less than two.  There can be at most two such integers.  We have seen that $\Big\lceil {f_{ir(\Gamma)}(\Delta) \over b_i(\Delta) + 1} \Big\rceil$ is the smallest possible such integer and $\Big\lfloor {f_{ir(\Gamma)}(\Delta) \over b_i(\Delta)} \Big\rfloor$ is the largest, so if $g_{r(\Gamma)}(\Gamma)$ has two possible values, these must be both of them.  If $g_{r(\Gamma)}(\Gamma)$ has only one possible value, then these expressions both give that one value.  \endproof

This leaves only the case where $g_1(\Gamma) = b_1(\Delta)$ and $r(\Gamma) = 2$.  If $f_{23}(\Delta)$ is close to $f_{12}(\Delta)$, something analogous to the above lemma could deal with this case.  The ratio between them could be arbitrarily large, however, which could leave arbitrarily many possible values of $g_2(\Gamma)$ that correspond to a known $g_1(\Gamma)$.

\begin{definition}
\textup{Let $\Delta$ be a 3-colored simplicial complex.  Define}
\begin{eqnarray*}
v(\Delta, t) & = & b_1(\Delta)f_{23}(\Delta) + (f_{12}(\Delta) - b_1(\Delta)t)\bigg(f_{13}(\Delta) - {b_1(\Delta)f_{23}(\Delta) \over t} \bigg) \textup{ and} \\ s(\Delta) & = & \sqrt{{f_{12}(\Delta)f_{23}(\Delta) \over f_{13}(\Delta)}}.
\end{eqnarray*}
\end{definition}

Note that it is immediate from the definition that $b_2(\Delta) = \lfloor s(\Delta) \rfloor$.

The next lemma gives an upper bound on $f_{123}(\Delta)$ that depends on $g_2(\Gamma)$.  If we construct a complex and count its facets, the number of facets of the complex is a lower bound for $m(\Delta)$.  If the next lemma says that $f_{123}(\Gamma)$ is less than this lower bound, then we know immediately that $\Gamma \not \in \mathcal{B}(\Delta)$, and there is no need to actually construct the complex.  This greatly restricts how many possible values of $g_2(\Gamma)$ we need to check.

\begin{lemma}  \label{b1quadr2}
Let $\Delta$ be a 3-colored simplicial complex and let $\Gamma \in \mathcal{F}(\Delta)$ with $g_1(\Gamma) = b_1(\Delta)$.  If $r(\Gamma) = 2$, then $f_{123}(\Gamma) \leq v(\Delta, g_2(\Gamma))$.
\end{lemma}

\proof  The lemma is stated as it is to make it clear that the upper bound depends only on the choice of $g_2(\Gamma)$, but it is easier to prove an alternate form.  We can compute
\begin{eqnarray*}
v(\Delta, g_2(\Gamma)) & = & b_1(\Delta)f_{23}(\Delta) + (f_{12}(\Delta) - b_1(\Delta)g_2(\Gamma))\bigg(f_{13}(\Delta) - {b_1(\Delta)f_{23}(\Delta) \over g_2(\Gamma)} \bigg) \\ & = & b_1(\Delta)f_{23}(\Delta) + j_3(\Gamma)\bigg(f_{13}(\Delta) - {b_1(\Delta)(g_2(\Gamma)g_3(\Gamma) + j_1(\Gamma)) \over g_2(\Gamma)} \bigg) \\ & = & b_1(\Delta)f_{23}(\Delta) + j_3(\Gamma)\bigg(f_{13}(\Delta) - b_1(\Delta)g_3(\Gamma) - {b_1(\Delta)j_1(\Gamma) \over g_2(\Gamma)} \bigg) \\ & = & b_1(\Delta)f_{23}(\Delta) + j_2(\Gamma)j_3(\Gamma) - {b_1(\Delta)j_1(\Gamma)j_3(\Gamma) \over g_2(\Gamma)}.
\end{eqnarray*}

\textbf{Case I:}  $p(\Gamma) = 3$

\textbf{Case I A:}  $j_2(\Gamma) \leq b_1(\Delta)$

The first extra vertex uses all remaining edges of both of its color sets.  This leaves no remaining edges of color 13 for use by the second extra vertex, so the second extra vertex adds no additional facets.  Therefore,
$$f_{123}(\Gamma) = b_1(\Delta)g_2(\Gamma)g_3(\Gamma) + j_1(\Gamma)j_2(\Gamma).$$

\textbf{Case I A 1:}  $j_1(\Gamma) \leq j_3(\Gamma)$
\begin{eqnarray*}
f_{123}(\Gamma) & = & b_1(\Delta)g_2(\Gamma)g_3(\Gamma) + j_1(\Gamma)j_2(\Gamma) \\ & = & b_1(\Delta)(f_{23}(\Delta) - j_1(\Gamma)) + j_2(\Gamma)j_3(\Gamma) + j_2(\Gamma)(j_1(\Gamma) - j_3(\Gamma)) \\ & \leq & b_1(\Delta)f_{23}(\Delta) + j_2(\Gamma)j_3(\Gamma) - {b_1(\Delta)j_1(\Gamma)j_3(\Gamma) \over g_2(\Gamma)}.
\end{eqnarray*}
The last line comes because $j_3(\Gamma) \geq j_1(\Gamma)$ and $j_3(\Gamma) \leq g_2(\Gamma)$, as there is no extra vertex of color 2.

\textbf{Case I A 2:}  $j_3(\Gamma) < j_1(\Gamma)$
\begin{eqnarray*}
f_{123}(\Gamma) & = & b_1(\Delta)g_2(\Gamma)g_3(\Gamma) + j_1(\Gamma)j_2(\Gamma) \\ & = & b_1(\Delta)(f_{23}(\Delta) - j_1(\Gamma)) + j_2(\Gamma)j_3(\Gamma) + j_2(\Gamma)(j_1(\Gamma) - j_3(\Gamma)) \\ & \leq & b_1(\Delta)f_{23}(\Delta) + j_2(\Gamma)j_3(\Gamma) - b_1(\Delta)j_1(\Gamma) + b_1(\Delta)(j_1(\Gamma) - j_3(\Gamma)) \\ & = & b_1(\Delta)f_{23}(\Delta) + j_2(\Gamma)j_3(\Gamma) - b_1(\Delta)j_3(\Gamma) \\ & \leq & b_1(\Delta)f_{23}(\Delta) + j_2(\Gamma)j_3(\Gamma) - {b_1(\Delta)j_1(\Gamma)j_3(\Gamma) \over g_2(\Gamma)}.
\end{eqnarray*}
As in the previous case, the last line comes because $j_1(\Gamma) \leq g_2(\Gamma)$.

\textbf{Case I B:}  $j_2(\Gamma) \geq b_1(\Delta)$

The first extra vertex is adjacent to all previous vertices of color 1, so it adds $j_1(\Gamma)b_1(\Delta)$ facets.  The second extra vertex adds $j_3(\Gamma)(j_2(\Gamma) - b_1(\Delta))$ facets.  Thus, we have
\begin{eqnarray*}
f_{123}(\Gamma) & = & b_1(\Delta)g_2(\Gamma)g_3(\Gamma) + j_1(\Gamma)b_1(\Delta) + j_3(\Gamma)(j_2(\Gamma) - b_1(\Delta)) \\ & = & b_1(\Delta)(f_{23}(\Delta) - j_1(\Gamma)) + j_1(\Gamma)b_1(\Delta) + j_3(\Gamma)(j_2(\Gamma) - b_1(\Delta)) \\ & = & b_1(\Delta)f_{23}(\Delta) + j_2(\Gamma)j_3(\Gamma) - b_1(\Delta)j_3(\Gamma) \\ & \leq & b_1(\Delta)f_{23}(\Delta) + j_2(\Gamma)j_3(\Gamma) - {b_1(\Delta)j_1(\Gamma)j_3(\Gamma) \over g_2(\Gamma)}.
\end{eqnarray*}

\textbf{Case II:}  $p(\Gamma) = 1$

\textbf{Case II A:}  $j_2(\Gamma) \leq g_3(\Gamma)$

The first extra vertex uses all available edges of color 13 and adds $j_2(\Gamma)j_3(\Gamma)$ facets.  This leaves no edges of this color set to be used by the second available vertex, so the other vertex adds no more facets.  This gives us
\begin{eqnarray*}
f_{123}(\Gamma) & = & b_1(\Delta)g_2(\Gamma)g_3(\Gamma) + j_2(\Gamma)j_3(\Gamma) \\ & = & b_1(\Delta)(f_{23}(\Delta) - j_1(\Gamma)) + j_2(\Gamma)j_3(\Gamma) \\ & = & b_1(\Delta)f_{23}(\Delta) + j_2(\Gamma)j_3(\Gamma) - b_1(\Delta)j_1(\Gamma) \\ & \leq & b_1(\Delta)f_{23}(\Delta) + j_2(\Gamma)j_3(\Gamma) - {b_1(\Delta)j_1(\Gamma)j_3(\Gamma) \over g_2(\Gamma)}.
\end{eqnarray*}

\textbf{Case II B:}  $j_2(\Gamma) > g_3(\Gamma)$

\textbf{Case II B 1:}  $j_3(\Gamma) \geq j_1(\Gamma)$

There are enough spare edges of color 13 to make the first extra vertex adjacent to all previous vertices of color 3.  Thus, the first extra vertex adds $g_3(\Gamma)j_3(\Gamma)$ facets.  The second extra vertex brings an additional $j_1(\Gamma)(j_2(\Gamma) - g_3(\Gamma))$ facets.  We can use these to compute
\begin{eqnarray*}
f_{123}(\Gamma) & = & b_1(\Delta)g_2(\Gamma)g_3(\Gamma) + g_3(\Gamma)j_3(\Gamma) + j_1(\Gamma)(j_2(\Gamma) - g_3(\Gamma)) \\ & \leq & b_1(\Delta)(f_{23}(\Delta) - j_1(\Gamma)) + g_3(\Gamma)j_3(\Gamma) + j_3(\Gamma)(j_2(\Gamma) - g_3(\Gamma)) \\ & = & b_1(\Delta)f_{23}(\Delta) - b_1(\Delta)j_1(\Gamma) + j_2(\Gamma)j_3(\Gamma) \\ & \leq & b_1(\Delta)f_{23}(\Delta) + j_2(\Gamma)j_3(\Gamma) - {b_1(\Delta)j_1(\Gamma)j_3(\Gamma) \over g_2(\Gamma)}.
\end{eqnarray*}

\textbf{Case II B 2:}  $j_3(\Gamma) < j_1(\Gamma)$

Define $\Gamma_1$ by $g_i(\Gamma_1) = g_i(\Gamma)$ for all $i \in [3]$, $p(\Gamma_1) = 3$, and $q(\Gamma_1) = 1$.  By Lemma~\ref{choosepq}, $f_{123}(\Gamma) \leq f_{123}(\Gamma_1)$. Furthermore, $\Gamma_1$ satisfies the bound of this lemma by Case~I.  \endproof

The bound of the last lemma depends only on $g_2(\Gamma)$ and the flag f-vector of $\Delta$.  If we multiply out the bound, the coefficients on the $g_2(\Gamma)$ and ${1 \over g_2(\Gamma)}$ terms are both negative.  Thus, for sufficiently large or small $g_2(\Gamma)$, the bound is small.  Once the bound is smaller than the number of facets of a known complex, we can discard the choice of $g_2(\Gamma)$ as obviously not giving a complex that maximizes the number of facets.  We would like to say that this happens very quickly.  The next lemma says that we can get close to the bound of the previous lemma.

\begin{lemma}  \label{b1nearquad}
Let $\Delta$ be a 3-colored simplicial complex with $f_{12}(\Delta) \leq f_{13}(\Delta)$ and $f_{12}(\Delta) \leq f_{23}(\Delta)$.  Suppose that there is a $\Gamma_0 \in \mathcal{F}(\Delta)$ with $g_1(\Gamma_0) = b_1(\Delta)$ and $r(\Gamma_0) = 2$.  Then there is a $\Gamma \in \mathcal{F}(\Delta)$ such that $f_{123}(\Gamma) \geq v(\Delta, g_2(\Gamma_0)) - f_{12}(\Delta)$.
\end{lemma}

\proof  Define $\Gamma$ by $g_i(\Gamma) = g_i(\Gamma_0)$ for all $i \in [3]$.  If $j_2(\Gamma_0) \geq b_1(\Delta)$, then let $p(\Gamma) = 3$ and $q(\Gamma) = 1$.  Otherwise, let $p(\Gamma) = 1$ and $q(\Gamma) = 3$.  Either $\Gamma$ is the same complex as $\Gamma_0$ or else it swaps $p(\Gamma_0)$ with $q(\Gamma_0)$.  This does not affect any of the criteria for $\mathcal{F}(\Delta)$, so $\Gamma \in \mathcal{F}(\Delta)$.

We can compute $b_1(\Delta) = \Big\lfloor {f_{12}(\Delta) f_{13}(\Delta) \over f_{23}(\Delta)} \Big\rfloor \leq \Big\lfloor {f_{13}(\Delta) f_{23}(\Delta) \over f_{12}(\Delta)} \Big\rfloor = b_3(\Delta)$ because $f_{12}(\Delta) \leq f_{23}(\Delta)$.  Suppose first that $g_3(\Gamma) \leq b_1(\Delta) - 2$.  In this case, we have
\begin{eqnarray*}
f_{13}(\Delta) & \leq & (g_1(\Gamma) + 1)(g_3(\Gamma) + 1) \\ & \leq & (b_1(\Delta) + 1)(b_1(\Delta) - 1) \\ & = & b_1(\Delta)^2 - 1 \\ & < & b_1(\Delta)^2 \\ & \leq & b_1(\Delta)b_3(\Delta) \\ & \leq & f_{13}(\Delta),
\end{eqnarray*}
a contradiction.  Therefore, $g_3(\Gamma) \geq b_1(\Delta) - 1$.  Thus, if $j_2(\Gamma) < b_1(\Delta)$, then $j_2(\Gamma) \leq g_3(\Gamma)$.  As such, in the proof of Lemma~\ref{b1quadr2}, we are either in case I B or case II A.

As with the previous lemma, it is more convenient to prove the alternative form
$$f_{123}(\Gamma) \geq b_1(\Delta)f_{23}(\Delta) + j_2(\Gamma)j_3(\Gamma) - {b_1(\Delta)j_1(\Gamma)j_3(\Gamma) \over g_2(\Gamma)} - f_{12}(\Delta).$$
For simplicity, let
$$z(\Gamma) = b_1(\Delta)f_{23}(\Delta) + j_2(\Gamma)j_3(\Gamma) - {b_1(\Delta)j_1(\Gamma)j_3(\Gamma) \over g_2(\Gamma)} - f_{12}(\Delta),$$
so that we are trying to prove that $f_{123}(\Gamma) \geq z(\Gamma)$.  We break this into the same cases as before and do not repeat the computations, but only check how far from inequality we are.

If $j_2(\Gamma) \geq b_1(\Delta)$, then from the arithmetic of Case~I~B of Lemma~\ref{b1quadr2}, we have
\begin{eqnarray*}
f_{123}(\Gamma) - z(\Gamma) & = & b_1(\Delta)j_3(\Gamma) \bigg( {j_1(\Gamma) \over g_2(\Gamma)} - 1 \bigg) + f_{12}(\Delta) \\ & \geq & -b_1(\Delta)j_3(\Gamma) + f_{12}(\Delta) \\ & \geq & f_{12}(\Delta) - g_1(\Gamma)g_2(\Gamma) \geq 0.
\end{eqnarray*}
Similarly, if $j_2(\Gamma) < b_1(\Delta)$, then the arithmetic of Case~II~A yields
\begin{eqnarray*}
f_{123}(\Gamma) - z(\Gamma) & = & b_1(\Delta)j_1(\Gamma) \bigg( {j_3(\Gamma) \over g_2(\Gamma)} - 1 \bigg) + f_{12}(\Delta) \\ & \geq & -b_1(\Delta)j_1(\Gamma) + f_{12}(\Delta) \\ & \geq & f_{12}(\Delta) - g_1(\Gamma)g_2(\Gamma) \geq 0. \qquad \endproof
\end{eqnarray*}

It is useful in certain computations to be able to divide by $b_1(\Delta)$.  This would be problematic if $b_1(\Delta) = 0$. The next lemma says that the case where $b_1(\Delta) = 0$ is easily handled, so it can otherwise be discarded.

\begin{lemma}  \label{b10}
If $\Delta$ is a 3-colored simplicial complex such that $b_1(\Delta) = 0$ and $\mathcal{D}(\Delta) \not = \emptyset$, then $m(\Delta) = f_{12}(\Delta)f_{13}(\Delta)$.
\end{lemma}

\proof  As $b_1(\Delta) = 0$, we have $\sqrt{{f_{12}(\Delta)f_{13}(\Delta) \over f_{23}(\Delta)}} < 1$, or equivalently, $f_{12}(\Delta)f_{13}(\Delta) < f_{23}(\Delta)$.  If $f_2(\Delta) < f_{12}(\Delta)$, then we have ${f_{12}(\Delta) \over f_2(\Delta)} > 1$, so $\Big\lfloor {f_{12}(\Delta) \over f_2(\Delta)} \Big\rfloor \geq 1$.  Furthermore, ${f_{23}(\Delta) \over f_2(\Delta)} > {f_{23}(\Delta) \over f_{12}(\Delta)} > f_{13}(\Delta)$, so $\Big\lfloor {f_{23}(\Delta) \over f_2(\Delta)} \Big\rfloor \geq f_{13}(\Delta)$.  Thus, $\Big\lfloor {f_{12}(\Delta) \over f_2(\Delta)} \Big\rfloor \Big\lfloor {f_{23}(\Delta) \over f_2(\Delta)} \Big\rfloor \geq f_{13}(\Delta)$, so by Lemma~\ref{vertedge}, $\mathcal{D}(\Delta) = \emptyset$, a contradiction.  By the same argument, if $f_3(\Delta) < f_{13}(\Delta)$, then $\mathcal{D}(\Delta) = \emptyset$.

Otherwise, define $\Gamma$ by $g_1(\Gamma) = 1$, $g_2(\Gamma) = f_{12}(\Delta)$, $g_3(\Gamma) = f_{13}(\Delta)$, and no extra vertices.  We have seen that there are enough vertices of each color to do this.  There are clearly enough edges of colors 12 and 13.  Since $f_{23}(\Gamma) = f_{12}(\Delta)f_{13}(\Delta) < f_{23}(\Delta)$, there are also enough edges of color 23.  Hence, $\Gamma$ is well-defined.  Since $f_{123}(\Gamma) = f_{12}(\Delta)f_{13}(\Delta)$, we get $m(\Delta) \geq f_{12}(\Delta)f_{13}(\Delta)$.

Conversely, each choice of an edge of color 12 and edge of color 13 specifies at least one vertex of each color, so there can be at most one facet containing these two edges.  Each facet must use an edge of each color set, so any $\Gamma_1 \in \mathcal{A}(\Delta)$ can have at most $f_{12}(\Delta)f_{13}(\Delta)$ facets.  Therefore, $m(\Delta) \leq f_{12}(\Delta)f_{13}(\Delta)$, and so the statement of the lemma follows.  \endproof

\begin{lemma}  \label{g2interval}
Let $\Delta$ be a 3-colored simplicial complex with $f_{12}(\Delta) \leq f_{13}(\Delta) \leq f_{23}(\Delta)$ and let $\Gamma_1, \Gamma_2 \in \mathcal{F}(\Delta)$ such that $g_1(\Gamma_1) = g_1(\Gamma_2) = b_1(\Delta)$ and $r(\Gamma_1) = r(\Gamma_2) = 2$.  If $x \in \mathbb{Z}$ such that $g_2(\Gamma_1) < x < g_2(\Gamma_2)$, then there is a complex $\Gamma_3 \in \mathcal{F}(\Delta)$ such that $g_2(\Gamma_3) = x$, $g_1(\Gamma_3) = b_1(\Delta)$ and $r(\Gamma_3) = 2$.
\end{lemma}

\proof  Try to define $\Gamma_3$ by $g_2(\Gamma_3) = x$, $g_1(\Gamma_3) = \big\lfloor {f_{12}(\Delta) \over x} \big\rfloor$, $g_3(\Gamma_3) = \big\lfloor {f_{23}(\Delta) \over x} \big\rfloor$, $p(\Gamma_3) = 1$, and $q(\Gamma_3) = 3$.  It follows from the definition that there are enough edges of colors 12 and 23 for $\Gamma_3$ to be well-defined, and that $\Gamma_3$ uses all edges of these two color sets if $f_1(\Delta) > g_1(\Gamma_3)$ and $f_3(\Delta) > g_3(\Gamma_3)$.  Since $\Gamma_1 \in \mathcal{F}(\Delta)$, we must have $\mathcal{D}(\Delta) \not = \emptyset$, so $\Gamma_3$ satisfies this condition for $\mathcal{F}(\Delta)$.

Because $\Gamma_1, \Gamma_2 \in \mathcal{F}(\Delta)$, we must have $b_1(\Delta)g_2(\Gamma_1) < b_1(\Delta)g_2(\Gamma_2) \leq f_{12}(\Delta)$ and $f_{12}(\Delta) \leq (b_1(\Delta) + 1)g_2(\Gamma_1)$.  If $f_1(\Delta) \leq b_1(\Delta)$, we would have $f_{12}(\Gamma_1) \leq b_1(\Delta)g_2(\Gamma_1) < f_{12}(\Delta)$, a contradiction.  Hence, $f_1(\Delta) > b_1(\Delta)$.  Furthermore,
$$b_1(\Delta)x < b_1(\Delta)g_2(\Gamma_2) \leq f_{12}(\Delta) \leq (b_1(\Delta) + 1)g_2(\Gamma_1) < (b_1(\Delta) + 1)x.$$
That $b_1(\Delta) < {f_{12}(\Delta) \over x} < b_1(\Delta) + 1$ guarantees that $g_1(\Gamma_3) = b_1(\Delta) < f_1(\Delta)$.  Similarly, we have
$$g_3(\Gamma_2)x < g_3(\Gamma_2)g_2(\Gamma_2) \leq f_{23}(\Delta) \leq (g_3(\Gamma_1) + 1)g_2(\Gamma_1) < (g_3(\Gamma_1) + 1)x.$$
This means $g_3(\Gamma_2) < {f_{23}(\Delta) \over x} < g_3(\Gamma_1) + 1$, from which $g_3(\Gamma_2) \leq g_3(\Gamma_3) < g_3(\Gamma_1)$.  Since $\Gamma_1$ is well-defined, $f_3(\Delta) \geq g_3(\Gamma_1) > g_3(\Gamma_3)$.

Next, we can compute
$$g_1(\Gamma_3)g_3(\Gamma_3) < g_1(\Gamma_1)g_3(\Gamma_1) \leq f_{13}(\Delta) \leq (g_1(\Gamma_2) + 1)(g_3(\Gamma_2) + 1)$$
$$\leq (g_1(\Gamma_3) + 1)(g_3(\Gamma_3) + 1).$$
This ensures that $\Gamma_3$ has enough edges of color 13 to be well-defined because $\Gamma_1$ is and enough vertices to use all of the edges because $\Gamma_2$ does.

Because $\Gamma_3$ uses at most as many vertices of color 2 as $\Gamma_2$ and at most as many of colors 1 and 3 as $\Gamma_1$, $\Gamma_3$ is well-defined.  We immediately have $\Gamma_3 \in \mathcal{A}(\Delta)$ by construction.  We have seen that $\Gamma_3$ uses all available edges of each color set, so $\Gamma_3 \in \mathcal{F}(\Delta)$.  \endproof

Finally we are in a position to show that the remaining case of $b_1(\Delta) = g_1(\Gamma)$ and $r(\Gamma) = 2$ can be checked by brute force in a reasonable number of steps.

\begin{lemma}  \label{g2bruteforce}
Let $\Delta$ be a 3-colored simplicial complex with $f_{12}(\Delta) \leq f_{13}(\Delta) \leq f_{23}(\Delta)$ and $b_1(\Delta) \geq 1$.  Suppose that there is some $\Gamma_0 \in \mathcal{D}(\Delta)$ with $g_1(\Gamma_0) = b_1(\Delta)$ and $r(\Gamma_0) = 2$.  Then we can find some $\Gamma \in \mathcal{D}(\Delta)$ by checking fewer than $6 + 2\sqrt{2}{\sqrt{f_{12}(\Delta)f_{23}(\Delta)} \over f_{13}(\Delta)}$ potential values of $g_2(\Gamma)$ and applying Lemma~\ref{knowr} to each potential value of $g_2(\Gamma)$ and $r(\Gamma) = 2$.
\end{lemma}

\proof  We start with some preliminary computations.  We wish to find the value of $t > 0$ that maximizes $v(\Delta, t)$.  One can readily compute ${\partial \over \partial t} v(\Delta, t) = -b_1(\Delta)f_{13}(\Delta) + {1 \over t^2}b_1(\Delta)f_{12}(\Delta)f_{23}(\Delta)$.  Setting the derivative equal to zero and solving for $t$ gives $t = s(\Delta)$.  Furthermore, ${\partial^2 \over \partial t^2} v(\Delta, t) = -{1 \over t^3}b_1(\Delta)f_{12}(\Delta)f_{23}(\Delta)$, which is negative for all $t > 0$, so this is a maximum.

Next, we compute how far from maximizing $v(\Delta, t)$ a given value of $t$ is.  For the former, if we define $z$ by $t = s(\Delta) + z$,  we compute
\begin{eqnarray*}
& & v(\Delta, s(\Delta)) - v(\Delta, s(\Delta) + z) \\ & = & b_1(\Delta)f_{23}(\Delta) + f_{12}(\Delta)f_{13}(\Delta) + b_1(\Delta)^2f_{23}(\Delta) - b_1(\Delta)f_{13}(\Delta)s(\Delta) \\ & & - {b_1(\Delta)f_{12}(\Delta)f_{23}(\Delta) \over s(\Delta)} - b_1(\Delta)f_{23}(\Delta) - f_{12}(\Delta)f_{13}(\Delta) - b_1(\Delta)^2f_{23}(\Delta) \\ & & + b_1(\Delta)f_{13}(\Delta)(s(\Delta) + z) + {b_1(\Delta)f_{12}(\Delta)f_{23}(\Delta) \over s(\Delta) + z} \\ & = & b_1(\Delta) \bigg( f_{13}(\Delta)z - {f_{12}(\Delta)f_{23}(\Delta)z \over s(\Delta)(s(\Delta) + z)} \bigg) \\ & = & b_1(\Delta) \bigg( f_{13}(\Delta)z - {f_{13}(\Delta)s(\Delta)z \over (s(\Delta) + z)} \bigg) \\ & = & b_1(\Delta)f_{13}(\Delta)z \bigg( 1 - {s(\Delta) \over (s(\Delta) + z)} \bigg) \\ & = & {b_1(\Delta)f_{13}(\Delta)z^2 \over (s(\Delta) + z)}.
\end{eqnarray*}

Suppose that there is a complex $\Gamma_0 \in \mathcal{F}(\Delta)$ with $g_1(\Gamma_0) = b_1(\Delta)$, $g_2(\Gamma_0) = b_2(\Delta)$, and $r(\Gamma_0) = 2$. For this complex, we get $s(\Delta) + z = b_2(\Delta)$.  Since $b_2(\Delta) = \lfloor s(\Delta) \rfloor$, we get $|z| < 1$.  As such, we get $v(\Delta, s(\Delta) + z) = v(\Delta, s(\Delta)) - {b_1(\Delta)f_{13}(\Delta)z^2 \over b_2(\Delta)}$.  Applying Lemma~\ref{b1nearquad} yields
\begin{eqnarray*}
f_{123}(\Gamma_0) & \geq & v(\Delta, s(\Delta)) - {b_1(\Delta)f_{13}(\Delta)z^2 \over b_2(\Delta)} - f_{12}(\Delta) \\ & > & v(\Delta, s(\Delta)) - {b_1(\Delta)f_{13}(\Delta) \over b_2(\Delta)} - f_{12}(\Delta).
\end{eqnarray*}

For simplicity, let $x(\Delta) = {b_1(\Delta)f_{12}(\Delta) \over b_2(\Delta)} + f_{12}(\Delta)$.  It thus suffices to check the values of $g_2(\Gamma)$ where $v(\Delta, g_2(\Gamma)) > f_{123}(\Gamma_0)$.  We set $g_2(\Gamma) = s(\Delta) + z$ and use Lemma~\ref{b1quadr2} compute
\begin{eqnarray*}
0 & < & v(\Delta, g_2(\Gamma)) - f_{123}(\Gamma_0) \\ & \leq & v(\Delta, s(\Delta) + z) - v(\Delta, s(\Delta)) + x(\Delta) \\ & = & v(\Delta, s(\Delta)) - {b_1(\Delta)f_{13}(\Delta)z^2 \over (s(\Delta) + z)} - v(\Delta, s(\Delta)) + x(\Delta) \\ & = & x(\Delta) - {b_1(\Delta)f_{13}(\Delta)z^2 \over (s(\Delta) + z)}.
\end{eqnarray*}

The inequality $x(\Delta) - {b_1(\Delta)f_{13}(\Delta)z^2 \over (s(\Delta) + z)} > 0$ clearly holds if $z = 0$ and fails if $z$ gets far enough away from zero.  Thus, to find the values of $z$ that make it true, it suffices to find the values that give equality and take the interval between them.  We compute $b_1(\Delta)f_{13}(\Delta)z^2 - zx(\Delta) - s(\Delta)x(\Delta) = 0$.  The quadratic formula gives
$$z = {x(\Delta) \pm \sqrt{x(\Delta)^2 + 4b_1(\Delta)f_{13}(\Delta)s(\Delta)x(\Delta)} \over 2b_1(\Delta)f_{13}(\Delta)}.$$
The difference between the two roots is
\begin{eqnarray*}
& & {\sqrt{x(\Delta)^2 + 4b_1(\Delta)f_{13}(\Delta)s(\Delta)x(\Delta)} \over b_1(\Delta)f_{13}(\Delta)} \\ & = & \sqrt{\bigg({x(\Delta) \over b_1(\Delta)f_{13}(\Delta)}\bigg)^2 + {4b_1(\Delta)f_{13}(\Delta)s(\Delta)x(\Delta) \over b_1(\Delta)^2f_{13}(\Delta)^2}} \\ & = & \sqrt{\bigg({x(\Delta) \over b_1(\Delta)f_{13}(\Delta)}\bigg)^2 + {4s(\Delta)x(\Delta) \over b_1(\Delta)f_{13}(\Delta)}} \\ & = & \sqrt{\bigg({x(\Delta) \over b_1(\Delta)f_{13}(\Delta)}\bigg)\bigg(4s(\Delta) + {x(\Delta) \over b_1(\Delta)f_{13}(\Delta)}\bigg)} \\ & = & \sqrt{\bigg({1 \over b_2(\Delta)} + {f_{12}(\Delta) \over b_1(\Delta)f_{13}(\Delta)}\bigg)\bigg(4s(\Delta) + {1 \over b_2(\Delta)} + {f_{12}(\Delta) \over b_1(\Delta)f_{13}(\Delta)}\bigg)} \\ & < & \sqrt{{2\sqrt{f_{13}(\Delta)} \over \sqrt{f_{12}(\Delta)f_{23}(\Delta)}} + {2f_{12}(\Delta)\sqrt{f_{23}(\Delta)} \over f_{13}(\Delta)\sqrt{f_{12}(\Delta)f_{13}(\Delta)}}} \\ & & * \sqrt{{4\sqrt{f_{12}(\Delta)f_{23}(\Delta)} \over \sqrt{f_{13}(\Delta)}} + {2\sqrt{f_{13}(\Delta)} \over \sqrt{f_{12}(\Delta)f_{23}(\Delta)}} + {2f_{12}(\Delta)\sqrt{f_{23}(\Delta)} \over f_{13}(\Delta)\sqrt{f_{12}(\Delta)f_{13}(\Delta)}}} \\ & = & \sqrt{8 + 8{f_{12}(\Delta)f_{23}(\Delta) \over f_{13}(\Delta)^2} + {4f_{13}(\Delta) \over f_{12}(\Delta)f_{23}(\Delta)} + {8 \over f_{13}(\Delta)} + {4f_{12}(\Delta)f_{23}(\Delta) \over f_{13}(\Delta)^3}}
\\ & < & \sqrt{24 + 8{f_{12}(\Delta)f_{23}(\Delta) \over f_{13}(\Delta)^2}} \\ & < & 5 + 2\sqrt{2}{\sqrt{f_{12}(\Delta)f_{23}(\Delta)} \over f_{13}(\Delta)}
\end{eqnarray*}
Above, we used that $f_{23}(\Delta) \leq f_{12}(\Delta)f_{13}(\Delta)$ (because $b_1(\Delta) \geq 1$) and $b_2(\Delta) \geq {1 \over 2}s_2(\Delta)$ (because $b_2(\Delta) \geq b_1(\Delta) \geq 1$ and $s_2(\Delta) - b_2(\Delta) < 1$).  Hence, there are fewer than $6 + 2\sqrt{2}{\sqrt{f_{12}(\Delta)f_{23}(\Delta)} \over f_{13}(\Delta)}$ integers in the interval, and so fewer than $6 + 2\sqrt{2}{\sqrt{f_{12}(\Delta)f_{23}(\Delta)} \over f_{13}(\Delta)}$ possible values of $g_2(\Gamma)$ to check.

Above, we assumed that one could take $g_2(\Gamma_0) = b_2(\Delta)$ and $g_1(\Gamma_0) = b_1(\Delta)$.  If this is not the case, then by Lemma~\ref{g2interval}, either all $\Gamma \in \mathcal{F}(\Delta)$ with $g_1(\Gamma) = b_1(\Delta)$ have $g_2(\Gamma) > b_2(\Delta)$ or else all have $g_2(\Gamma) < b_2(\Delta)$.  Since $g_2(\Gamma)$ is an integer, all such $\Gamma$ have $g_2(\Gamma)$ on the same side of $s(\Delta)$.

We have seen that $v(\Delta, t)$ attains its maximum at $t = s(\Delta)$ and that ${\partial^2 \over \partial t^2} v(\Delta, t) < 0$ for all $t > 0$.  Let
\begin{eqnarray*}
m_2(\Delta) & = & \max\{v(\Delta, g_2(\Gamma)) - f_{12}(\Delta)\ |\ \Gamma \in \mathcal{F}(\Delta), r(\Gamma) = 2, g_1(\Gamma) = b_1(\Delta)\}, \\ c_1(\Delta) & = & \min\{g_2(\Gamma)\ |\ \Gamma \in \mathcal{F}(\Delta), r(\Gamma) = 2, g_1(\Gamma) = b_1(\Delta), \\ & & \qquad v(\Delta, g_2(\Gamma)) > m_2(\Delta)\}, \qquad \textup{and} \\ c_2(\Delta) & = & \max\{g_2(\Gamma)\ |\ \Gamma \in \mathcal{F}(\Delta), r(\Gamma) = 2, g_1(\Gamma) = b_1(\Delta), \\ & & \qquad v(\Delta, g_2(\Gamma)) > m_2(\Delta)\}.
\end{eqnarray*}
Suppose that $c_1(\Delta) > b_2(\Delta)$.  If $c_2(\Delta) - c_1(\Delta) \geq 5 + 2\sqrt{2}{\sqrt{f_{12}(\Delta)f_{23}(\Delta)} \over f_{13}(\Delta)}$, then
\begin{eqnarray*}
& & v(\Delta, c_1(\Delta)) - v(\Delta, c_2(\Delta)) \\ & > & v(\Delta, s(\Delta)) - v(\Delta, s(\Delta) + c_2(\Delta) - c_1(\Delta)) \qquad \textup{(because ${\partial^2 \over \partial t^2} v(\Delta, t) < 0$)} \\ & > & v(\Delta, s(\Delta)) - v \Bigg(\Delta, s(\Delta) + 5 + 2\sqrt{2}{\sqrt{f_{12}(\Delta)f_{23}(\Delta)} \over f_{13}(\Delta)} \Bigg) \\ & > & v (\Delta, s(\Delta)) - v \bigg(\Delta, s(\Delta) + {x(\Delta) + \sqrt{x(\Delta)^2 + 4b_1(\Delta)f_{13}(\Delta)s(\Delta)x(\Delta)} \over 2b_1(\Delta)f_{13}(\Delta)}\bigg) \\ & = & f_{12}(\Delta).
\end{eqnarray*}

This gives that $m_2(\Delta) \geq v(\Delta, c_1(\Delta)) - f_{12}(\Delta) > v(\Delta, c_2(\Delta))$, a contradiction.
Similarly, if $c_2(\Delta) < b_2(\Delta)$ and $c_2(\Delta) - c_1(\Delta) \geq 5 + 2\sqrt{2}{\sqrt{f_{12}(\Delta)f_{23}(\Delta)} \over f_{13}(\Delta)}$, we get $v(\Delta, c_2(\Delta)) - v(\Delta, c_1(\Delta)) > f_{12}(\Delta)$, which gives that $m_2(\Delta) \geq v(\Delta, c_2(\Delta)) - f_{12}(\Delta) > v(\Delta, c_1(\Delta))$, a contradiction.  Therefore, if either $c_1(\Delta) > b_2(\Delta)$ or $c_2(\Delta) < b_2(\Delta)$, the lemma holds.  It is clear from the definitions that $c_1(\Delta) \leq c_2(\Delta)$, so the only other possibility is that $c_1(\Delta) \leq b_2(\Delta) \leq c_2(\Delta)$.  In this case, by Lemma~\ref{g2interval}, there is such a complex $\Gamma_0$, and so we have already seen that the lemma holds. \endproof

In the last lemma, the question arose of what values of $g_2(\Gamma)$ can occur if $r(\Gamma) = 2$ and $g_1(\Gamma) = b_1(\Delta)$.  The next lemma settles this issue.

\begin{lemma}  \label{b1g2}
Let $\Delta$ be a 3-colored simplicial complex with $f_{12}(\Delta) \leq f_{13}(\Delta) \leq f_{23}(\Delta)$.  Suppose further that $\mathcal{D}(\Delta) \not = \emptyset$ and that there is a complex $\Gamma \in \mathcal{F}(\Delta)$ such that $g_1(\Gamma) = b_1(\Delta)$ and $r(\Gamma) = 2$.  This guarantees that $f_1(\Delta) \geq b_1(\Delta)$.  Furthermore, if $f_1(\Delta) = b_1(\Delta)$, then $g_2(\Gamma) = {f_{12}(\Delta) \over f_1(\Delta)}$.  If $f_1(\Delta) > b_1(\Delta)$, then
\begin{enumerate}
\item $f_{13}(\Delta) \leq f_3(\Delta)(b_1(\Delta)+1)$,
\item $g_2(\Gamma) \geq {f_{12}(\Delta) \over f_1(\Delta)}$,
\item $g_2(\Gamma) \geq {f_{12}(\Delta) \over b_1(\Delta) + 1}$,
\item $g_2(\Gamma) \geq {f_{23}(\Delta) \over \big\lfloor {f_{13}(\Delta) \over b_1(\Delta)} \big\rfloor + 1}$,
\item $g_2(\Gamma) \geq {f_{23}(\Delta) \over f_3(\Delta)}$,
\item $g_2(\Gamma) \leq {f_{12}(\Delta) \over b_1(\Delta)}$,
\item $g_2(\Gamma) \leq {f_{23}(\Delta) \over \big\lceil {f_{13}(\Delta) \over b_1(\Delta) + 1} \big\rceil - 1}$, and
\item $g_2(\Gamma) \leq f_2(\Delta)$.
\end{enumerate}
\end{lemma}

\proof  We must have $f_1(\Delta) \geq g_1(\Gamma)$, and so if there is a $\Gamma$ with $b_1(\Delta) = g_1(\Gamma)$, then we must have $f_1(\Delta) \geq b_1(\Delta)$.

Suppose that $f_1(\Delta) = b_1(\Delta)$.  Note that this means that $\Gamma$ can have no extra vertex of color 1, in addition to having no extra vertex of color 2 because $r(\Gamma) = 2$.  There are few enough edges of color 12 for $\Gamma$ to use them all if and only if $f_1(\Delta)g_2(\Gamma) \geq f_{12}(\Delta)$.  There are enough edges of color 12 for $\Gamma$ to be well-defined if and only if $b_1(\Delta)g_2(\Gamma) \leq f_{12}(\Delta)$.  Hence, equality must hold, and so we get $g_2(\Gamma) = {f_{12}(\Delta) \over f_1(\Delta)}$.

Otherwise, we have $f_1(\Delta) > b_1(\Delta)$.  There are enough vertices of color 2 to define the complex if and only if $g_2(\Gamma) \leq f_2(\Delta)$, which is condition 8.  There are enough vertices of color 1 to handle the edges of color 12 if and only if $f_1(\Delta)g_2(\Gamma) \geq f_{12}(\Delta)$, or equivalently, $g_2(\Gamma) \geq {f_{12}(\Delta) \over f_1(\Delta)}$, which is condition 2.  There are enough vertices of color 3 to deal with the edges of color 23 if and only if $f_3(\Delta)g_2(\Gamma) \geq f_{23}(\Delta)$, or equivalently, $g_2(\Gamma) \geq {f_{23}(\Delta) \over f_3(\Delta)}$, which is condition 5.

There are enough edges of color 12 for the complex to be well-defined if and only if $b_1(\Delta)g_2(\Gamma) \leq f_{12}(\Delta)$, or equivalently, $g_2(\Gamma) \leq {f_{12}(\Delta) \over b_1(\Delta)}$, which is condition 6.  There are enough edges of color 23 for the complex to be well-defined if and only if $g_2(\Gamma)g_3(\Gamma) \leq f_{23}(\Delta)$, or equivalently, $g_3(\Gamma) \leq {f_{23}(\Delta) \over g_2(\Gamma)}$.  There are enough edges of color 13 for the complex to be well-defined if and only if $b_1(\Delta)g_3(\Gamma) \leq f_{13}(\Delta)$, or equivalently, $g_3(\Gamma) \leq {f_{13}(\Delta) \over b_1(\Delta)}$.  Since $g_3(\Gamma)$ is an integer, we can take floors of its upper bounds to get $g_3(\Gamma) \leq \Big\lfloor {f_{23}(\Delta) \over g_2(\Gamma)} \Big\rfloor$ and $g_3(\Gamma) \leq \Big\lfloor {f_{13}(\Delta) \over b_1(\Delta)} \Big\rfloor$.

There are few enough edges of color 12 for $\Gamma$ to use them all if and only if $f_{12}(\Delta) \leq (b_1(\Delta) + 1)g_2(\Gamma)$, or equivalently, $g_2(\Gamma) \geq {f_{12}(\Delta) \over b_1(\Delta) + 1}$, which is condition 3.  There are few enough edges of color 23 for the complex to use them all if and only if $g_2(\Delta)(g_3(\Gamma)+1) \geq f_{23}(\Delta)$, or equivalently, $g_2(\Gamma) \geq {f_{23}(\Delta) \over g_3(\Delta)+1}$.

There are few enough edges of color 13 for $\Gamma$ to use them all if and only if both $f_{13}(\Delta) \leq f_3(\Delta)(b_1(\Delta)+1)$ and $(b_1(\Delta) + 1)(g_3(\Gamma) + 1) \geq f_{13}(\Delta)$.  The former is condition 1, and the latter is equivalent to $g_3(\Gamma) \geq {f_{13}(\Delta) \over b_1(\Delta) + 1} - 1$.  Since $g_3(\Gamma)$ is an integer, we can take the ceiling and get $g_3(\Gamma) \geq \Big\lceil {f_{13}(\Delta) \over b_1(\Delta) + 1} \Big\rceil - 1$.

Thus, we have that in order to make $g_3(\Gamma)$ compatible with the choice of $g_1(\Gamma) = b_1(\Delta)$ and not use more vertices of color 3 than are available, our bounds are $g_3(\Gamma) \geq \Big\lceil {f_{13}(\Delta) \over b_1(\Delta) + 1} \Big\rceil - 1$, $g_3(\Gamma) \leq \Big\lfloor {f_{13}(\Delta) \over b_1(\Delta)} \Big\rfloor$, and $g_3(\Gamma) \leq f_3(\Delta)$.  In order to then make $g_2(\Gamma)$ compatible with $g_3(\Gamma)$, our bounds are $g_2(\Gamma) \geq {f_{23}(\Delta) \over g_3(\Gamma)+1}$ and $g_2(\Gamma)g_3(\Gamma) \leq f_{23}(\Delta)$, from which $g_2(\Gamma) \leq {f_{23}(\Delta) \over g_3(\Gamma)}$.  We plug in our bounds on $g_3(\Gamma)$ to get $g_2(\Gamma) \leq {f_{23}(\Delta) \over \big\lceil {f_{13}(\Delta) \over b_1(\Delta) + 1} \big\rceil - 1 }$ and $g_2(\Gamma) \geq {f_{23}(\Delta) \over \big\lfloor {f_{13}(\Delta) \over b_1(\Delta)} \big\rfloor + 1 }$, which are conditions 7 and 4, respectively.  \endproof

It is possible to prove more in the above lemma, but it isn't necessary for our purposes, as we are mainly interested in restricting how many cases there are to check.  More precisely, if $f_1(\Delta) = b_1(\Delta)$, then there is a complex $\Gamma \in \mathcal{F}(\Delta)$ such that $g_1(\Gamma) = b_1(\Delta)$ and $r(\Gamma) = 2$ if and only if $g_2(\Gamma) = {f_{12}(\Delta) \over f_1(\Delta)}$, $g_2(\Gamma) \leq f_2(\Delta)$, $\max\Big\{  {f_{23}(\Delta) \over z}, {f_{13}(\Delta) \over f_1(\Delta)} \Big\} \leq f_3(\Delta),$ and
$$\max\Big\{ \Big\lceil {f_{23}(\Delta) \over z} \Big\rceil - 1, \Big\lceil {f_{13}(\Delta) \over f_1(\Delta)} \Big\rceil - 1 \Big\} \leq \min\Big\{ \Big\lfloor {f_{23}(\Delta) \over z} \Big\rfloor, \Big\lfloor {f_{13}(\Delta) \over f_1(\Delta)} \Big\rfloor\Big\}.$$
In addition, if $f_1(\Delta) > b_1(\Delta)$, then the converse of the above lemma holds as well.  That is, if $z$ is an integer that satisfies all eight of the listed conditions for $g_2(\Gamma)$, then there is a complex $\Gamma \in \mathcal{F}(\Delta)$ such that $g_1(\Gamma) = b_1(\Delta)$, $r(\Gamma) = 2$, and $g_2(\Gamma) = z$.

Finally, we reach the main theorem.  This basically summarizes the lemmas of this section, and gives a method guaranteed to produce a complex with the maximal number of facets subject to the restrictions on the number of vertices and edges of various color sets.

\begin{theorem}  \label{maintheorem}
Given positive integers $f_1(\Delta)$, $f_2(\Delta)$, $f_3(\Delta)$, $f_{12}(\Delta)$, $f_{13}(\Delta)$, and $f_{23}(\Delta)$, the following procedure will suffice to compute $m(\Delta)$.
\begin{enumerate}
\item  Check whether the inequalities $f_1(\Delta)f_2(\Delta) \geq f_{12}(\Delta)$, $f_1(\Delta)f_3(\Delta) \geq f_{13}(\Delta)$, and $f_2(\Delta)f_3(\Delta) \geq f_{23}(\Delta)$ all hold.  If not, then there is no $\Delta$ with the desired flag f-numbers, so stop.
\item  Check the inequalities of Lemma~\ref{vertedge}.  If any of them hold, then the lemma gives $m(\Delta)$, so stop.
\item  Relabel the colors if necessary to ensure that $f_{12}(\Delta) \leq f_{13}(\Delta) \leq f_{23}(\Delta)$.
\item  Compute $b_1(\Delta)$, $b_2(\Delta)$, and $b_3(\Delta)$.  If $b_1(\Delta) = 0$, then Lemma~\ref{b10} gives $m(\Delta)$, so stop.
\item  Attempt to construct complexes where $g_{r(\Gamma)}(\Gamma) = b_{r(\Gamma)}(\Delta)$ for each of $r(\Gamma) = 1$ and $r(\Gamma) = 2$ as described in Lemma~\ref{knowr}.  Compute $f_{123}(\Gamma)$ for each such complex that is well-defined.
\item  Attempt to construct complexes where $g_3(\Gamma) = b_3(\Delta)$ and $r(\Gamma) = 2$ as explained in Lemma~\ref{bibr}.  Compute $f_{123}(\Gamma)$ for each such complex that is well-defined.
\item  Repeat the previous step using $r(\Gamma) = 1$.
\item  Repeat the previous step using $g_2(\Gamma) = b_2(\Delta)$ (and $r(\Gamma) = 1$).
\item  Attempt to construct a complex $\Gamma$ with $r(\Gamma) = 3$ and $g_3(\Gamma) = f_3(\Delta)$ as explained in Lemma~\ref{knowr}.  Compute $f_{123}(\Gamma)$ if the complex is well-defined.
\item  Use Lemma~\ref{b1g2} to compute the maximum and minimum possible values of $g_2(\Gamma)$ if $g_1(\Gamma) = b_1(\Delta)$ and $r(\Gamma) \not = 1$.
\item If it is possible to have $g_2(\Gamma) = b_2(\Delta)$, then construct such a complex as explained in Lemma~\ref{g2bruteforce}.  Decrease $g_2(\Gamma)$ by 1 and construct the complexes again repeatedly until either it is not possible to construct complexes or Lemma~\ref{b1quadr2} says that it decreasing $g_2(\Gamma)$ further will necessarily give no more facets than an already known complex.  Likewise, try $g_2(\Gamma) = b_2(\Delta) + 1$ and increase $g_2(\Gamma)$ by 1 and construct complexes repeatedly until they are not defined or the lemma says that increasing $g_2(\Gamma)$ further will necessarily give no more facets than an already known complex.
\item If Lemma~\ref{b1g2} gives a lower bound on $g_2(\Gamma)$ that is greater than $b_2(\Delta)$, then try setting $g_2(\Gamma)$ to this lower bound and construct a complex as explained in Lemmas~\ref{g2bruteforce}.  Increase $g_2(\Gamma)$ by 1 and construct complexes again repeatedly until we stop as in the previous step.
\item If Lemma~\ref{b1g2} gives a upper bound on $g_2(\Gamma)$ that is less than $b_2(\Delta)$, then try setting $g_2(\Gamma)$ to this upper bound and construct a complex as explained in Lemma~\ref{g2bruteforce}.  Decrease $g_2(\Gamma)$ by 1 and construct complexes again repeatedly until we stop as in the previous step.
\item Compare the values of $f_{123}(\Gamma)$ for the various complexes constructed.  The largest such value is $m(\Delta)$.
\end{enumerate}
Furthermore, this process requires computing the number of facets of fewer than $15 + 2\sqrt{2}{\sqrt{f_{12}(\Delta)f_{23}(\Delta)} \over f_{13}(\Delta)}$ complexes.
\end{theorem}

\proof  If the inequalities in point 1 hold, then we can easily construct $\Delta$ by picking arbitrary subsets of the appropriate sizes of the possible edges of each color set.  In this case, it is clear from the definitions that $\mathcal{C}(\Delta) \not = \emptyset$.  By Lemma~\ref{alledges}, either $\mathcal{D}(\Delta) \not = \emptyset$ or else Lemma~\ref{vertedge} completes the problem in step 2.  In the former case, Lemma~\ref{oneb} gives that $\mathcal{E}(\Delta) \not = \emptyset$.

There are three ways to pick a value of $i$ such that $g_i(\Gamma) = b_i(\Delta)$ and three ways to pick a value of $r(\Gamma)$, for nine possibilities in all.  Part five handles two of these nine cases, and parts six through eight each handle one.  If $b_1(\Delta) = 0$, then Lemma~\ref{b10} solves the problem.  Otherwise, parts 10 through 13 handle a sixth case.

If $\mathcal{D}(\Delta) = \mathcal{E}(\Delta)$, then Lemma~\ref{nor3} says that either step 9 finds a complex in $\mathcal{E}(\Delta)$ or else one of the other six cases has such a complex.  If $\mathcal{E}(\Delta)$ is a proper subset of $\mathcal{D}(\Delta)$, then Lemma~\ref{oneb} ensures that one of the other six cases produces a complex in $\mathcal{E}(\Delta)$.  Therefore, we are guaranteed to find a complex in $\mathcal{E}(\Delta)$ by this procedure if there is one.

The following table summarizes the nine cases and says which lemmas give the upper bounds on how many complexes it could be necessary to construct for that particular case.
$$\begin{array}{cccc}
g_1(\Gamma) = b_1(\Delta) & r(\Gamma) = 1 & \textup{trivial} & 1 \\ g_1(\Gamma) = b_1(\Delta) & r(\Gamma) = 2 & \textup{Lemma~\ref{g2bruteforce}} & <6 + 2\sqrt{2}{\sqrt{f_{12}(\Delta)f_{23}(\Delta)} \over f_{13}(\Delta)} \\g_2(\Gamma) = b_2(\Delta) & r(\Gamma) = 1 & \textup{Lemma~\ref{bibr}} & 2 \\ g_2(\Gamma) = b_2(\Delta) & r(\Gamma) = 2 & \textup{trivial} & 1 \\ g_3(\Gamma) = b_3(\Delta) & r(\Gamma) = 1 & \textup{Lemma~\ref{bibr}} & 2 \\ g_3(\Gamma) = b_3(\Delta) & r(\Gamma) = 2 & \textup{Lemma~\ref{bibr}} & 2 \\  & r(\Gamma) = 3 & \textup{Lemma~\ref{nor3}} & 1
\end{array}$$
Add up all of the cases to get fewer than $15 + 2\sqrt{2}{\sqrt{f_{12}(\Delta)f_{23}(\Delta)} \over f_{13}(\Delta)}$ complexes to check in total.  \endproof

The bound on how many complexes we have to check is a worst-case scenario, and the number we have to actually construct by this procedure is usually much smaller than the given bound.  One reason for this is that quite often, we construct a complex from the same set of parameters at multiple steps.  This commonly happens for complexes that have $g_i(\Gamma) = b_i(\Delta)$ for more than one value of $i$.  If this happens, we could see that we have already computed the number of facets of a given complex once, and not bother to compute it again the second or third time it shows up.

The other reason why this is an overestimate is that in Lemma~\ref{g2bruteforce}, we implicitly assumed that the bound of Lemma~\ref{b1quadr2} is as bad as it can possibly be at every single step until the last possible moment, at which point we suddenly have $j_1(\Gamma) = 0$ and hit the bounds of the lemma exactly.  If ${\sqrt{f_{12}(\Delta)f_{23}(\Delta)} \over f_{13}(\Delta)}$ is large, this is quite improbable.  The usual scenario is that as $g_2(\Gamma)$ varies, it won't take very many values to stumble upon a case where $j_1(\Gamma)$ is either near $g_2(\Gamma)$ or else very small compared to $g_2(\Gamma)$.  This causes us to nearly hit the bound of Lemma~\ref{b1quadr2}, which greatly restricts how many additional cases we have to check, rather than relying on Lemma~\ref{b1nearquad}.

It is also worthwhile to note that the quantity ${\sqrt{f_{12}(\Delta)f_{23}(\Delta)} \over f_{13}(\Delta)}$ is rarely large.  For a large integer $n$, if one picks $f_{12}(\Delta)$, $f_{13}(\Delta)$, and $f_{23}(\Delta)$ uniformly at random from $[n]$ and then sorts them to make $f_{12}(\Delta) \leq f_{13}(\Delta) \leq f_{23}(\Delta)$, an easy triple integral approximation finds that the expected value of ${\sqrt{f_{12}(\Delta)f_{23}(\Delta)} \over f_{13}(\Delta)}$ is essentially ${8 \over 9}$.  Therefore, the expected number of complexes that one must check by the method of Theorem~\ref{maintheorem} is less than 18.

We can actually do better than that.  If we use the line $\sqrt{24 + 8{f_{12}(\Delta)f_{23}(\Delta) \over f_{13}(\Delta)^2}}$ from the proof of Lemma~\ref{g2bruteforce}, this has an average value of $4\sqrt{2} < 6$.  If we use this rather than 9 as the approximation for the average upper bound on the number of complexes to check in steps 10 through 13, then on average, we have to check fewer than 15 complexes.  In practice, it tends to be a lot less than this, even.

\section{Some examples}

In the previous section, Theorem~\ref{maintheorem} explained how to compute $m(\Delta)$.  In this section, we give some examples of how the procedure works, with both some typical cases and some extremal ones to argue that it would likely be impractical to greatly improve upon Theorem~\ref{maintheorem}, so it a satisfactory solution to the problem.  We start with a few trivial examples.

\begin{example}
\textup{Let $f_1(\Delta) = 3$, $f_2(\Delta) = 5$, $f_3(\Delta) = 7$, $f_{12}(\Delta) = 23$, $f_{13}(\Delta) = 14$, and $f_{23}(\Delta) = 18$.  We compute $f_{12}(\Delta) = 23 > 15 = f_1(\Delta)f_2(\Delta)$, so there is no 3-colored complex $\Delta$ having the given face numbers, and we stop.}
\end{example}

\begin{example}
\textup{Let $f_1(\Delta) = 3$, $f_2(\Delta) = 5$, $f_3(\Delta) = 7$, $f_{12}(\Delta) = 13$, $f_{13}(\Delta) = 16$, and $f_{23}(\Delta) = 18$.  The inequalities of part 1 of Theorem~\ref{maintheorem} hold, so we move on.  In part 2, we compute}
$$\bigg\lfloor {f_{12}(\Delta) \over f_1(\Delta)} \bigg\rfloor \bigg\lfloor {f_{13}(\Delta) \over f_1(\Delta)} \bigg\rfloor = \bigg\lfloor {13 \over 3} \bigg\rfloor \bigg\lfloor {16 \over 3} \bigg\rfloor = 20 \geq 18 = f_{23}(\Delta).$$
\textup{Thus, Lemma~\ref{vertedge} asserts that $m(\Delta) = f_1(\Delta)f_{23}(\Delta) = (3)(18) = 54$.  The lemma also explains how to find a complex $\Gamma$ with the desired flag f-numbers and 54 facets, if desired.}
\end{example}

\begin{example}
\textup{Let $f_1(\Delta) = 17$, $f_2(\Delta) = 31$, $f_3(\Delta) = 25$, $f_{12}(\Delta) = 15$, $f_{13}(\Delta) = 12$, and $f_{23}(\Delta) = 279$.  The inequalities of point 1 hold and those of point 2 fail, so neither settles the problem and we move on.  Step 3 advises us to ensure that $f_{12}(\Delta) \leq f_{13}(\Delta) \leq f_{23}(\Delta)$.  This does not hold with the numbers as given, as $15 > 12$.  We want to rearrange the colors such that $f_{12}(\Delta) = 12$, $f_{13}(\Delta) = 15$, and $f_{23}(\Delta) = 279$.  This can be done by swapping colors 2 and 3, which also gives us $f_2(\Delta) = 25$ and $f_3(\Delta) = 31$.}

\textup{Step 4 starts by computing}
$$b_1(\Delta) = \Bigg\lfloor \sqrt{{f_{12}(\Delta)f_{13}(\Delta) \over f_{23}(\Delta)}} \Bigg\rfloor = \Bigg\lfloor \sqrt{(12)(15) \over 279} \Bigg\rfloor \approx \lfloor .803 \rfloor = 0.$$
\textup{Since $b_1(\Delta) = 0$, Lemma~\ref{b10} tells us that} $$m(\Delta) = f_{12}(\Delta)f_{13}(\Delta) = (12)(15) = 180.$$
\end{example}

The next example is a typical use of the full Theorem~\ref{maintheorem}.  It has few enough complexes in $\mathcal{F}(\Delta)$ that it is easy to compute them all, so that Lemma~\ref{g2bruteforce} doesn't particularly matter.

\begin{example}
\textup{Let $f_1(\Delta) = 533$, $f_2(\Delta) = 471$, $f_3(\Delta) = 818$, $f_{12}(\Delta) = 4972$, $f_{13}(\Delta) = 5311$, and $f_{23}(\Delta) = 5630$.  We can quickly compute that steps 1 and 2 do not solve the problem, and the numbers of edges are already sorted as step 3 dictates.  Step 4 asks us to compute $b_1(\Delta) = 68$, $b_2(\Delta) = 72$, and $b_3(\Delta) = 77$.}

\textup{The remaining steps essentially ask us to brute force the various complexes in $\mathcal{F}(\Delta)$.  We list the step at which we construct each complex, the parameters of the complex, and the number of facets.  When we hit on parameters used earlier, we note it and do not reconstruct a complex that we have already used.}

$$\begin{array}{cccccc}
\textup{step} & g_1(\Gamma) & g_2(\Gamma) & g_3(\Gamma) & r(\Gamma) & f_{123}(\Gamma) \\ 5 & 68 & 73 & 78 & 1 & \textup{undefined} \\ 5 & 69 & 72 & 78 & 2 & \textup{undefined} \\ 6 & 68 & 73 & 77 & 2 & 382896 \\ 8 & 69 & 72 & 76 & 1 & 382736 \\ 9 & 6 & 6 & 818 & 3 & \textup{not in } \mathcal{F}(\Delta) \\ 11 & 68 & 73 & 77 & 2 & \textup{previous} \end{array}$$

\textup{There is no complex for step 7 because the condition of Lemma~\ref{bibr} is violated.  We could have quickly discarded the two undefined complexes of step 5 on the basis that it has $g_i(\Gamma) > b_i(\Delta)$ for two values of $i$.  We do not bother to invoke Lemma~\ref{g2bruteforce} for steps 11-13, as there are few enough complexes that we can find them all by brute force.  By inspection, $m(\Delta) = 382896$.}
\end{example}

The next example gives a typical demonstration of the power of Lemma~\ref{g2bruteforce}.  The class $\mathcal{F}(\Delta)$ is huge, but this lemma lets us compute few enough complexes that we can list them all here.

\begin{example}
\textup{Let $f_1(\Delta) = 13$, $f_2(\Delta) = 5471$, $f_3(\Delta) = 3818$, $f_{12}(\Delta) = 1843$, $f_{13}(\Delta) = 2157$, and $f_{23}(\Delta) = 3150248$.  We can quickly compute that steps 1 and 2 do not solve the problem, and the numbers of edges are already sorted as step 3 dictates.  Step 4 asks us to compute $b_1(\Delta) = 1$, $b_2(\Delta) = 1640$, and $b_3(\Delta) = 1920$.}

\textup{This time, there aren't very many possible complexes outside of steps 11-13, but in these final steps, we get complexes in $\mathcal{F}(\Delta)$ with $g_2(\Gamma)$ ranging from 1460 to 1843.  A direct brute force approach would require checking several hundred complexes.  Fortunately, Lemma~\ref{b1quadr2} immediately allows us to limit the computations to values of $g_2(\Delta)$ ranging from 1637 to 1644.}

$$\begin{array}{cccccc}
\textup{step} & g_1(\Gamma) & g_2(\Gamma) & g_3(\Gamma) & r(\Gamma) & f_{123}(\Gamma) \\ 5 & 1 & 1842 & 2156 & 1 & \textup{undefined} \\ 5 & 1 & 1640 & 1920 & 2 & 3198156 \\ 6 & 1 & 1640 & 1920 & 2 & \textup{previous} \\ 9 & 0 & 825 & 3818 & 3 & \textup{not in } \mathcal{F}(\Delta) \\ 11 & 1 & 1640 & 1920 & 2 & \textup{previous} \\ 11 & 1 & 1641 & 1919 & 2 & 3198122 \\ 11 & 1 & 1642 & 1918 & 2 & 3198086 \\ 11 & 1 & 1643 & 1917 & 2 & 3198048 \\ 11 & 1 & 1644 & 1916 & 2 & 3198008 \\ 11 & 1 & 1639 & 1922 & 2 & 3198098 \\ 11 & 1 & 1638 & 1923 & 2 & 3198013 \\ 11 & 1 & 1637 & 1924 & 2 & 3198040
\end{array}$$

\textup{If one of the complexes computed later had more facets than the ones we computed before reaching step 11, that could have further restricted how many complexes we would have to compute in step 11.  Regardless, this is still far more efficient than having to compute the number of facets of every single complex in $\mathcal{F}(\Delta)$.  Note that it was sufficient to try 8 complexes.  For comparison, Theorem~\ref{maintheorem} said that we would need to do the computations for at most 114 complexes.}
\end{example}

Finally, we wish to demonstrate that finding the complex with the maximal number of vertices can force $g_2(\Delta)$ to be arbitrarily far away from $b_2(\Delta)$.  More precisely, the difference can be on the order of ${\sqrt{f_{12}(\Delta)f_{23}(\Delta)} \over f_{13}(\Delta)}$ even as this quantity becomes arbirarily large.

\begin{example}
\textup{Pick any real number $t$ and let $f_1(\Delta) = 2$, $f_2(\Delta) = \lfloor 100^t \rfloor$, $f_3(\Delta) = \lfloor 100^t \rfloor$, $f_{12}(\Delta) = \lfloor 100^t \rfloor$, $f_{13}(\Delta) = \lfloor 100^t + 2(10)^t \rfloor$, and $f_{23}(\Delta) = \Big\lfloor \big\lfloor {2 \over 3}f_{12}(\Delta) \big\rfloor \big( \big\lfloor {2 \over 3}f_{13}(\Delta)\big\rfloor + c \big) \Big\rfloor$, for some real number $c$ near ${1 \over 2}$.  We can compute $b_1(\Delta) = 1$, $b_2(\Delta) \approx {2 \over 3}f_{12}(\Delta)$, and $b_3(\Delta) \approx {2 \over 3}f_{13}(\Delta)$.  Furthermore, for a suitable choice of $c$, the complex that maximizes $f_{123}(\Gamma)$ has}
$$g_2(\Gamma) \approx b_2(\Gamma) - 2.4(10)^t \approx b_2(\Gamma) - .36{\sqrt{f_{12}(\Delta)f_{23}(\Delta)} \over f_{13}(\Delta)}.$$

\textup{For example, if $t = 2$ and $c = .443$, then we get $f_{12}(\Delta) = 10000$, $f_{13}(\Delta) = 10200$, $f_{23}(\Delta) = 45331753$, $b_2(\Delta) = 6666$, and $b_3(\Delta) = 6799$.  After running through the various possibilities of Theorem~\ref{maintheorem}, we see that the complex that maximizes $f_{123}(\Gamma)$ has $g_2(\Gamma) = 6643$, which differs from $b_2(\Gamma)$ by 23.  For comparison, ${\sqrt{f_{12}(\Delta)f_{23}(\Delta)} \over f_{13}(\Delta)} \approx 66$.}

\textup{Larger values of $t$ also let us take $c$ to be a little smaller.  For example, if $t = 4$ and $c = .416$, then we get $f_{12}(\Delta) = 100000000$, $f_{13}(\Delta) = 100020000$, $f_{23}(\Delta) = 4445333316613330$, $b_2(\Delta) = 66666666$, and $b_3(\Delta) = 66679999$.  Letting a computer run the necessary computations yields that $m(\Delta) = 5556666649191260$, and the complex that produces this many facets has $r_2(\Gamma) = 2$ and $g_2(\Gamma) = 66664202$.  This differs from $b_2(\Delta)$ by 2464, and for comparison, ${\sqrt{f_{12}(\Delta)f_{23}(\Delta)} \over f_{13}(\Delta)} \approx 6666$.}
\end{example}

In the previous example, $f_{12}(\Delta)$ was very close to $f_{13}(\Delta)$, which means that ${\sqrt{f_{12}(\Delta)f_{23}(\Delta)} \over f_{13}(\Delta)} \approx \sqrt{f_{23}(\Delta) \over f_{13}(\Delta)}$.  The next example sets $c = .45$ and then generalizes the previous example and shows that the number of complexes required can still be on the order of ${\sqrt{f_{12}(\Delta)f_{23}(\Delta)} \over f_{13}(\Delta)}$ even as ${f_{13}(\Delta) \over f_{12}(\Delta)}$ is arbitrarily large.

\begin{example}
\textup{Pick any positive real number $t$ and any integer $w$ and let $f_1(\Delta) = 2$, $f_2(\Delta) = \lfloor 100^t \rfloor$, $f_3(\Delta) = \lfloor w100^t \rfloor$, $f_{12}(\Delta) = \lfloor 100^t \rfloor$, $f_{13}(\Delta) = \lfloor w100^t + 2\sqrt{w}(10)^t \rfloor$, and $f_{23}(\Delta) = \Big\lfloor \big\lfloor {2 \over 3}f_{12}(\Delta) \big\rfloor \big( \big\lfloor {2 \over 3}f_{13}(\Delta)\big\rfloor + .45 \big) \Big\rfloor$.  We can compute $b_1(\Delta) = 1$, $b_2(\Delta) \approx {2 \over 3}f_{12}(\Delta)$, and $b_3(\Delta) \approx {2 \over 3}f_{13}(\Delta)$.  Furthermore, if $t$ is large enough that ${f_{13}(\Delta) \over f_{12}(\Delta)} \approx w$, the complex that maximizes $f_{123}(\Gamma)$ has}
$$g_2(\Gamma) \approx b_2(\Gamma) - .23{10^t \over \sqrt{w}} \approx b_2(\Gamma) - .35{\sqrt{f_{12}(\Delta)f_{23}(\Delta)} \over f_{13}(\Delta)}.$$

\textup{For example, if we set $t = 3$ and $w = 100$, we get $f_{13}(\Delta) = 100020000$, $f_{23}(\Delta) = 44453289179999$, $b_2(\Delta) = 666666$, and $b_3(\Delta) = 66679966$.  A computer search finds that $m(\Delta) = 55566644505542$ and the complex that attains this bound has $g_2(\Gamma) = 666643$.  This differs from $b_2(\Delta)$ by 23; for comparison, ${\sqrt{f_{12}(\Delta)f_{23}(\Delta)} \over f_{13}(\Delta)} \approx 67$.}
\end{example}

In these examples, in order for the complex $\Gamma$ that maximizes the number of facets to have $g_2(\Gamma)$ far away from $b_2(\Delta)$, it is necessary that many consecutive possible values of $g_2(\Gamma)$ have $j_1(\Gamma)$ much larger than 0 and much smaller than $g_2(\Gamma)$.  If this happens, then decreasing $g_2(\Gamma)$ by 1 increases $g_3(\Gamma)$ by the same amount ($w$ in the above example) many consecutive times.  This additional structure makes it easy to get a formula for $f_{123}(\Gamma)$ as a function of $g_2(\Gamma)$ that holds for many consecutive values of $g_2(\Gamma)$, which can greatly reduce the computations needed to find $m(\Delta)$ in the particularly bad cases where Theorem~\ref{maintheorem} calls for constructing a large number of simplicial complexes. Thus, even the worst cases are not nearly so bad as they seem.

Of course, one could still hope for a quick and clever solution to this problem as has happened with some previous characterizations of f-vectors of various classes of complexes. The next example explains why an easy characterization is improbable, as adding one extra vertex or edge can dramatically change the complex that maximizes the number of facets.

\begin{example}
\textup{Let $f_1(\Delta) = 2$, $f_2(\Delta) = 6683$, $f_3(\Delta) = 7000$, $f_{12}(\Delta) = 10000$, $f_{13}(\Delta) = 10200$, and $f_{23}(\Delta) = 45331745$.  We can compute that $m(\Delta) = 56664978$.  Furthermore, there is only one complex $\Gamma \in \mathcal{F}(\Delta)$ such that $f_{123}(\Gamma) = 56664978$, and it has $p(\Gamma) = 1$, $q(\Gamma) = 3$, $g_1(\Gamma) = 1$, $g_2(\Gamma) = 6683$, and $g_3(\Gamma) = 6783$.}

\textup{If we set $f_2(\Delta) = 6682$ and leave the rest of the flag f-numbers unchanged, this obviously excludes the previously optimal complex.  This time, we get $m(\Delta) = 56664977$, which corresponds to two complexes $\Gamma_1, \Gamma_2 \in \mathcal{F}(\Delta)$.  The two complexes are defined by $p(\Gamma_1) = 3$, $q(\Gamma_1) = 1$, $g_1(\Gamma_1) = 1$, $g_2(\Gamma_1) = 6643$, $g_3(\Gamma_1) = 6823$, $p(\Gamma_2) = 2$, $q(\Gamma_2) = 1$, $g_1(\Gamma_2) = 1$, $g_2(\Gamma_2) = 6642$, and $g_3(\Gamma_2) = 6824$.  What happened in this example is that $g_2(\Gamma)$ for the unique $\Gamma \in \mathcal{D}(\Delta)$ was quite far to one side of $b_2(\Delta) = 6666$, and changing the number of allowed vertices of one color by 1 made it so that there were two complexes $\Gamma_1, \Gamma_2 \in \mathcal{D}(\Delta)$, both of which $g_2(\Gamma_1)$ and $g_2(\Gamma_2)$ quite far on the other side of $b_2(\Delta)$.}

\textup{Furthermore, we can get similar results by adding one edge.  Let $f_{13}(\Delta) = 10201$ and leave the rest of the flag f-numbers the same as in the original example.  This time, we get $m(\Delta) = 56668334$, and there are again two complexes $\Gamma_1, \Gamma_2 \in \mathcal{F}(\Delta)$ such that $f_{123}(\Gamma_1) = f_{123}(\Gamma_2) = 56668334$.  These two complexes are defined by exactly the same parameters as $\Gamma_1$ and $\Gamma_2$ had in the previous paragraph; the extra edge merely adds some extra facets.  This time, the big change in the structure of the complex is not due to a cap on the number of vertices; the same complexes would still be the only ones in $\mathcal{D}(\Delta)$ even if $f_2(\Delta)$ were greatly increased.  One can still define $\Gamma$ by the same parameters as before, but this time, $f_{123}(\Gamma) = 56668295 < m(\Delta)$.}
\end{example}

This same behavior also occurs with smaller numbers, but if $g_2(\Gamma)$ differs from $b_2(\Delta)$ by only 1 or 2, it is much less clear what happened.

\section{More colors}

Having characterized the flag f-vectors of 3-colored complexes, it is natural to ask whether the characterization carries over to more colors.  Unfortunately, even the case of four colors is dramatically more complicated than that of three.

The basic approach of the three color case does carry over, however.  Recall that we started by ignoring the discreteness of faces and allowing non-integer numbers of vertices.  The same scheme can be done with more colors, and is along the lines of what Walker did in \cite{walker}.

If given a proposed flag f-vector on n colors $\{f_S\}_{S \subset [n]}$, one can propose that the faces of color set $S$ be a complete $|S|$-partite complex on some vertices of each color of $S$.  That is, if $S = \{i_1, i_2, \dots, i_n\}$, we can suppose that the faces of color set $S$ consist of all ways to choose one vertex out of $f_S^{i_1}$ of color set $i_1$, one vertex out of $f_S^{i_2}$ of color set $i_2$, and so forth, with the restriction that $f_S = f_S^{i_1}f_S^{i_2} \dots f_S^{i_n}$.  The simplicial complex restriction that any subface of a face must itself be a face corresponds to the requirement $f_T^i \geq f_S^i$ for every $i \in T$ and $T \subset S$.

As Walker did, we can take the logarithms of both sides and get $\log(f_S) = \log(f_S^{i_1}) + \log(f_S^{i_2}) + \dots + \log(f_S^{i_n})$.  This turns the problem into a linear programming problem of maximizing $\log(f_{[n]})$ subject to the known values of $\log(f_S)$ and the inequalities $f_T^i \geq f_S^i$.  If one can find the optimal solution in the continuous case, one could hope that the optimal solution in the discrete case would be nearby.

Unfortunately, not only is it unclear how to find an efficient solution in the discrete case, but with four or more colors, having a solution in the continuous case doesn't even guarantee that there is a solution in the discrete case.  As we saw earlier, if we set $f_{12}(\Delta) = f_{13}(\Delta) = f_{23}(\Delta) = 5$, the optimal solution in the continuous case is $f_{123} = 5\sqrt{5} > 11$, but the discrete case only allows 9 facets.  If we use these same numbers as part of a flag f-vector for a four-colored complex and try to require $f_{123}(\Delta) = 11$, we may well find solutions in the continuous case, but there will be no solution in the discrete case.  Unlike the case of three colors, faces of dimension two are no longer facets, and cannot be ignored simply by posing the problem as one of maximizing the number of facets.

Regardless of whether this method can be extended to higher dimensions, it does provide a non-trivial class of examples where the exact characterization is known.  Any proposed theorem toward characterizing the flag f-vectors of colored complexes or the flag h-vectors of balanced Cohen-Macaulay complexes or balanced shellable complexes can now be checked against the known, exact result in the case of three colors.

\end{document}